\documentclass[onefignum,onetabnum]{siamonline250211}

\usepackage{lipsum}
\usepackage{amsfonts}
\usepackage{graphicx}
\usepackage{epstopdf}
\ifpdf
  \DeclareGraphicsExtensions{.eps,.pdf,.png,.jpg}
\else
  \DeclareGraphicsExtensions{.eps}
\fi

\usepackage{enumitem}
\setlist[enumerate]{leftmargin=.5in}
\setlist[itemize]{leftmargin=.5in}

\newsiamremark{remark}{Remark}
\newsiamremark{hypothesis}{Hypothesis}
\crefname{hypothesis}{Hypothesis}{Hypotheses}
\newsiamthm{claim}{Claim}
\newsiamremark{fact}{Fact}
\crefname{fact}{Fact}{Facts}

\headers{Greedy learning to optimize with convergence guarantees}{P. Fahy, M. Golbabaee, M. J. Ehrhardt}

\title{Greedy Learning to Optimize with Convergence Guarantees\thanks{Submitted to the editors June 25, 2025.
\funding{Patrick Fahy is supported by a scholarship from the EPSRC center for Doctoral Training in Statistical Applied Mathematics at Bath (SAMBa), under the project EP/S022945/1. MJE acknowledges support from the EPSRC (EP/T026693/1, EP/V026259/1, EP/Y037286/1) and the European Union Horizon 2020 research and innovation programme under the Marie Skodowska-Curie grant agreement REMODEL. MG acknowledges support from EPSRC grant EP/X001091/1.
}}}

\author{Patrick Fahy\thanks{Department of Mathematical Sciences, University of Bath, Bath, BA2 7AY, UK (\email{pf341@bath.ac.uk}).}
\and Mohammad Golbabaee\thanks{Department of Engineering Mathematics, University of Bristol, Bristol, BS8 1TW, UK (\email{m.golbabaee@bristol.ac.uk}).}
\and Matthias J. Ehrhardt\thanks{Department of Mathematical Sciences, University of Bath, Bath, BA2 7AY, UK
  (\email{m.ehrhardt@bath.ac.uk}).}}

\usepackage{amsopn}

\ifpdf
\hypersetup{
  pdftitle={Greedy learning to optimize with convergence guarantees},
  pdfauthor={P. Fahy, M. Golbabaee, M. J. Ehrhardt}
}
\fi

\headers{Greedy learning to optimize with convergence}{P. Fahy, M. Golbabaee, M. J. Ehrhardt}
\usepackage{subcaption}

\usepackage[utf8]{inputenc} % allow utf-8 input
\usepackage[T1]{fontenc}    % use 8-bit T1 fonts
\usepackage{amsfonts}       % blackboard math symbols
\usepackage{nicefrac}       % compact symbols for 1/2, etc.
\usepackage{microtype}      % microtypography
\usepackage{xcolor}         % colors

\usepackage{amsmath} % for \arg\min and other math commands

\newcommand{\ft}{\mathcal{F}}

\usepackage{array}
\usepackage{amssymb}

\usepackage[numbers]{natbib}
\usepackage{hyperref}

\usepackage{siunitx}
\usepackage{graphicx}
\usepackage{subcaption}

\newtheorem{assumption}{Assumption}

\usepackage{lipsum}
\usepackage{amsfonts}
\usepackage{graphicx}
\usepackage{epstopdf}

\usepackage{algorithmicx}

\usepackage{amsmath} % for mathematical symbols
\usepackage{algorithm} % for the algorithm environment
\usepackage{algpseudocode} % for the pseudocode environment (algorithmic replacement)

\begin{document}

\maketitle

% REQUIRED
\begin{abstract}
   Learning to optimize (L2O) is an approach that leverages training data to accelerate the solution of optimization problems. Many approaches use unrolling to parametrize the update step and learn optimal parameters. This requires choosing an unroll length in advance and typically incurs memory costs that grow with the number of unrolled iterations. Moreover, learned optimizers often lack convergence guarantees. We propose a greedy L2O framework that learns iteration-dependent parameters sequentially, by minimizing the average objective value after one additional optimization step. This converts the finite-unroll training problem into a sequence of one-step learning problems, allowing training to proceed until a stopping criterion is met while maintaining constant device-memory requirements with respect to the number of training iterations. We apply this framework to a preconditioned variant of Polyak's heavy-ball method, using several parameterizations of the learned preconditioners, including a convolutional parameterization for imaging problems. For linear parameterizations, the greedy learning subproblem is convex whenever the training objectives are convex, and admits closed-form solutions for least-squares objectives. We prove convergence in the training set even when the preconditioners are not necessarily symmetric nor positive definite. Convergence on a class of unseen functions is also obtained under certain assumptions. We test our learned algorithms on two inverse problems, image deblurring and Computed Tomography, on which learned convolutional preconditioners demonstrate improved empirical performance over classical optimization algorithms such as Nesterov's Accelerated Gradient Method and the quasi-Newton method L-BFGS.
\end{abstract}

% REQUIRED

% REQUIRED

\section{Introduction}

We consider the optimization problem
\begin{equation} \label{eq--orig-opt}
    \min_{x \in \mathcal X} f(x),
\end{equation}
where \(f:\mathcal X \to \mathbb R\) is \(L\)-smooth and has a minimizer $f^*$, and \(\mathcal X\) is a Hilbert space. Classical optimization algorithms are built in a theoretically justified manner, with guarantees on their performance and convergence properties. However, practitioners often concentrate on problems within a much smaller class. This is often the case for linear inverse problems, where an observation $y \in \mathcal{Y}$ is generated from a ground-truth $x_{\text{true}}$ via some linear forward operator $A: \mathcal{X} \to \mathcal{Y}$, such that $y = Ax_{\text{true}} + \varepsilon$, where $\varepsilon \in \mathcal{Y}$ is some random noise, and the goal is to recover $x_{\text{true}}$. For example, in reconstructing images from blurred observations $y$ generated by a blurring operator $A$, one might minimize a function from the class:
\begin{equation} \label{eq--example-function-class}
    \mathcal{F} = \left\{ f: \mathcal{X} \to \mathbb{R} : f(x) = \frac12 \| Ax - y \|^2 + \mathcal{S}(x), y \sim \mathcal{P}(\mathcal{Y}) \right\},
\end{equation}
where $\mathcal{S}: \mathcal{X} \to \mathbb{R}$ is a chosen regularizer and $\mathcal{P}(\mathcal{Y})$ is some probability distribution on $\mathcal{Y}$ detailing the observations $y$ of interest. 

\paragraph{A preconditioned heavy-ball algorithm} Many classical algorithms can be seen as modifying the gradient direction. Newton's method accelerates convergence by applying the inverse Hessian to the gradient, which can be costly in practice. Quasi-Newton methods like BFGS \cite{nocedal2006quasi} approximate the Hessian, and L-BFGS \cite{liu1989limited} is used when BFGS is too memory-intensive. Methods that exploit momentum, such as Polyak's Heavy Ball (HB) method \cite{polyak1964} and Nesterov's Accelerated Gradient Method (NAG) \cite{nesterov1983}, use information from the past two iterates. In this work, we aim to accelerate optimization by learning preconditioners $G_t$ and $H_t$ in the update 
\begin{equation}
     x^{t+1} = x^t - G_t \nabla f\left(x^t\right) + H_t (x^t - x^{t-1}).
\end{equation}
% $x^{t+1} = x^t - G_t \nabla f\left(x^t\right) + H_t (x^t - x^{t-1})$.

\paragraph{Learning an optimization algorithm} Learning to optimize (L2O) aims to use training data from a class \(\mathcal F\) to learn an update rule that performs well on new objectives from the same class. The solution at each iteration $t$ is updated by a parametrized function, for example $R_\theta : \mathcal{X}^3 \to \mathcal{X} $ (i.e. the update rule) as dependent on parameters $\theta_t$ at iteration $t$ as $x^{t+1} = x^t - R_{\theta_t}\left(x^t, x^{t-1}, \nabla f\left(x^t\right)\right)$.
Ideally, one would like to learn parameters that give good performance over the whole optimization trajectory, for example by using the infinite-horizon loss
\begin{equation} \label{eq--infinite-horizon-l2o}
    \mathcal L_\infty((\theta_t)_{t\geq 0})
    =
    \mathop{\mathbb{E}}_{f \in \mathcal{F}, x^0}
    \left[
        \sum_{t=1}^{\infty}
        \omega_t \bigl(f(x^t)-f^\star\bigr)
    \right],
\end{equation}
where  \((\omega_t)_{t\geq 0}\) are non-negative
weights. In practice, however, standard unrolling methods \cite{monga2021algorithm} replace
\cref{eq--infinite-horizon-l2o} by a finite-horizon problem and often consider fixed parameters $\theta_t = \theta$ for all $t$. That is, for some handpicked \(T>0\), they train parameters by minimizing a loss of the form
\begin{equation} \label{eq--finite-unroll-l2o}
    \mathcal L_T(\theta)
    =
    \mathop{\mathbb{E}}_{f \in \mathcal{F}, x^0}
    \left[
        \sum_{t=0}^{T}
        \omega_t f\left(x^t\right)
    \right].
\end{equation}
Thus, the unroll length \(T\) must be chosen before training. For large \(T\), training is typically expensive both in terms of memory and computational cost. The value of an L2O method should be understood in an amortized sense, which is natural in inverse problems where many observations are reconstructed using the same forward model. Let \(C_{\rm train}\) be the cost of learning an optimizer, \(C_{\rm base}\) the average cost of solving one problem with a baseline method, and \(C_{\rm learned}\) the corresponding cost using the learned optimizer. The learned algorithm is then computationally preferable when applied to \(M\) unseen problems if $M >C_{\rm train}/({C_{\rm base}-C_{\rm learned}})$.

\subsection{Existing works}
% optimization problems are widely used in imaging applications, such as in image denoising and deblurring, and medical imaging with CT, and MRI. For this reason, $\mathcal{X}$ may not necessarily equal $\mathbb{R}^n$.
Learned optimization algorithms often lack convergence guarantees, including many that use RNNs \cite{andrychowicz2016learning, metz2019understanding} or Reinforcement learning \cite{2016RLL2O}. \cite{mathstructures} propose a learned optimization algorithm in which, at each iteration, a neural network outputs both a diagonal coordinate-wise scaling of the current gradient and an additional additive displacement term. This parametrization is motivated by the structure of classical optimization algorithms with known convergence properties. However, the learned update is not itself provably convergent.

Other approaches achieve provable convergence, which can be enforced with safeguarding \cite{heaton2023safeguarded}, or constructing convergent algorithms by learning parameters within a provably convergent set \cite{banert2020, banert2024accelerated}. \cite{tan2023boosting, tan2023mirror} learn mirror maps using input-convex neural networks within the mirror descent optimization algorithm such that the algorithm is provably convergent. \cite{learningalgohyperparams} learn scalar hyperparameters for first-order methods applied to convex optimization problems. Their framework uses time-varying hyperparameters during an initial regime and then switches to fixed hyperparameters up to a prescribed maximum unroll length. They also derive closed-form optimal step sizes for gradient descent with finite lookahead: one step for general objectives, and two- and three-step lookahead for least-squares problems. Lastly, \cite{l2OPACSucker} and \cite{PACOther} apply the PAC-Bayes framework to L2O, giving probabilistic generalization bounds for learned optimizers to address how performance on a finite training set transfers to unseen problem instances.

While L2O often seeks to exploit training data from a targeted family of problems, there is also a line of work that aims to design first-order methods with optimal worst-case guarantees over broad classes, such as $L$-smooth convex or $(L,\mu)$-smooth strongly convex functions. The Performance Estimation Problem (PEP) \cite{drori2014performance,taylor2017smooth,taylor2017exact} casts the worst-case performance of a first-order method after a prescribed number of iterations as an optimization problem over all objective functions. Reformulations yield tractable semidefinite programs and, in several important settings, exact worst-case bounds. Rather than learning from data on a narrow problem family, it seeks methods that are optimal over an entire function class. 

Lastly, adaptive algorithms improve optimization during use. For example, Armijo line-search \cite{armijo1966minimization} seeks to find a good step size at each iteration, while methods like AdaGrad \cite{duchi2011adaptive} and optimal diagonal preconditioners \cite{qu2024optimal} adapt preconditioners throughout the optimization process. Online optimization \cite{hazan2016introduction}, with methods such as Coin Betting \cite{orabona2016coin} and Adaptive Bound Optimization \cite{mcmahan2010adaptive}, offer a game-theoretic perspective.

% Adaptive algorithms focus on improving optimization algorithms during employment, for example, Armijo line-search \cite{armijo1966minimization} finds a suitable step size at each iteration. Algorithms such as  AdaGrad \cite{duchi2011adaptive} and optimal diagonal preconditioners \cite{qu2024optimal} search for better preconditioners during optimization. Similarly, online optimization \cite{hazan2016introduction} is a game-theoretic interpretation of optimization, with methods such as  Coin Betting \cite{orabona2016coin} and Adaptive Bound Optimization \cite{mcmahan2010adaptive}.

% Similarly, online optimization \cite{hazan2016introduction} focuses on improving optimization algorithms during employment. Techniques like Armijo line-search \cite{armijo1966minimization} and optimal diagonal preconditioners for quadratics \cite{qu2024optimal} adjust updates to accelerate empirical convergence. Adaptive methods such as AdaGrad \cite{duchi2011adaptive} and bound optimization \cite{mcmahan2010adaptive} address the non-smooth setting. Multidimensional Armijo focuses on smooth strongly convex functions \cite{kunstner2023searching}, and methods such as Coin Betting \cite{orabona2016coin}, optimize learning rates without prior knowledge.

Our work is also related to unrolled and layer-wise optimization-inspired architectures for inverse problems. Approximate-message-passing (AMP) methods \cite{donoho2009message} are iterative algorithms for large-scale linear inverse problems, and AMP-inspired networks such as LAMP/LVAMP \cite{borgerding2017amp} and LDAMP \cite{metzler2017learned} unfold these iterations into neural networks with layer-dependent parameters. Similarly, methods based on Learned ISTA \cite{LISTA} learn parameters in finite-depth networks, and \cite{chen2018theoretical} provides convergence theory for a structured unfolded ISTA model. These unfolded architectures learn a finite-depth solver tailored to a fixed number of layers, whereas our method learns iteration-dependent parameters greedily and is designed to continue until a stopping criterion is met while retaining convergence guarantees.

\subsection{Contributions and Outline}

Our paper contributes in the following ways: 
 \begin{itemize}
    \item We introduce a greedy learning-to-optimize framework to convert the infinite-horizon problem \cref{eq--infinite-horizon-l2o} into an infinite sequence of one-step learning problems. To do this, we learn the parameters of a preconditioned heavy-ball method sequentially to minimize the average next-step objective value on the training set, rather than by backpropagating through a fixed unrolled trajectory. Thus, the unroll length does not need to be chosen before training; instead, training may continue until a prescribed stopping criterion is met. Moreover, the method has constant device-memory requirements with respect to the number of training iterations: \cref{sec--greedy-L2O}. 
     \item Convergence in the training set is proved even when learned preconditioners are neither symmetric nor positive definite: \cref{section--convergence}.  Furthermore, convergence is proved on a class of unseen functions using soft constraints for parameter learning.
     \item Learning parameters is no more difficult than solving the initial optimization problem. For example, when objective functions $f$ are convex, parameter learning is a convex optimization problem for `linear parametrizations' of $G_t$ and $H_t$, enabling training that is significantly faster, with closed-form solutions for least-squares functions: \cref{section--linear-params}.
     \item A novel parametrization of $G_t$ and $H_t$ as convolutions. This provides a structured and computationally efficient parameterization that is well suited to imaging problems.
     \item We test our learned algorithms on two inverse problems, image deblurring and Computed Tomography, on which we demonstrate improved empirical performance over classical optimization algorithms such as Nesterov's Accelerated Gradient Method and the quasi-Newton method L-BFGS: \cref{sec--experiments}. Furthermore, we compare our learned algorithm to the deviation approach by \cite{banert2024accelerated}.
 \end{itemize}

\section{Notation}
Let $\mathcal{X}$ be a Hilbert space with corresponding field $\mathbb{R}$ and norm $\|\cdot \|$. 
% \comment{Throughout this paper, we take $\mathcal{X} = \mathbb{R}^{h_1 \times h_2}$ for $h_1, h_2 \in \mathbb{N}$. But then I would need to define $x_j$ somehow. If I use $x_{ij}$ then for the linear operators on $\mathcal{X}$, I would need to use four indices instead of the normal two}. 
A function $f: \mathcal{X} \rightarrow \mathbb{R}$ is $L$-smooth with parameter $L > 0$ if its gradient is Lipschitz continuous, i.e., if for all $x, y \in \mathcal{X}$, $\|\nabla f(x) - \nabla f(y)\| \leq L \|x - y\|$. A function $f: \mathcal{X} \rightarrow \mathbb{R}$ is bounded below if there exists some $M \in \mathbb{R}$ such that $f(x) \geq M$ for all $x \in \mathcal{X}$. We say that $f \in \mathcal{F}_L$ if $f$ is continuously differentiable, $L$-smooth, and has a minimizer. A function $f: \mathcal{X} \rightarrow \mathbb{R}$ is convex if for all $x, y \in \mathcal{X}$ and for all $\alpha \in [0, 1]$, $f(\alpha x + (1 - \alpha) y) \leq \alpha f(x) + (1 - \alpha) f(y)$. Furthermore, a function $f: \mathcal{X} \rightarrow \mathbb{R}$ is strongly convex with parameter $\mu > 0$ if $f - \mu\| \cdot \|^2/2$ is convex. If  $f \in \mathcal{F}_L$ is $\mu$-strongly convex, we say $f \in \mathcal{F}_{L, \mu}$.

% \begin{itemize}
%     \item $f \in \mathcal{F}_L$ if $f$ is convex, continuously differentiable and $L$-smooth. 
%     \item $f \in \mathcal{F}_{L, \mu}$ if, in addition, $f$ is $\mu$-strongly convex.
% \end{itemize}
We assume that the Hilbert space $\mathcal{X}$ has dimension $\operatorname{dim}(\mathcal{X}) = n$ and, therefore, admits a finite orthonormal basis $\{e_1, \dots, e_n \}$. For $x,y \in \mathcal{X}$ and $j  \in \{1, \dots, n\}$, define $x_j := \langle x, e_j \rangle$ and the pointwise product $x \odot y$ by $[x \odot y]_j = x_jy_j$, and the pointwise division by $[x\oslash y]_i = x_i/y_i$. For Hilbert spaces $\mathcal{X}$ and $\mathcal{Y}$, denote the space of linear operators from $\mathcal{X}$ to $\mathcal{Y}$ by $\mathcal{L}(\mathcal{X}, \mathcal{Y})$. If $\mathcal{Y} = \mathcal{X}$, we write $\mathcal{L}(\mathcal{X})$. For example, if $\mathcal{X} = \mathbb{R}^n$, $\mathcal{L}(\mathcal{X})$ is the space of $n \times n$ matrices.  Denote the adjoint of $A \in \mathcal{L}(\mathcal{X}, \mathcal{Y})$ by $A^*$, meaning that for $x \in \mathcal{X}, y \in \mathcal{Y}$, $\langle Ax, y \rangle = \langle x, A^*y \rangle$. We take $I \in \mathcal{L}(\mathcal{X})$ as the identity operator: $I(x) = x$ for all $x \in \mathcal{X}$, and assume there exists $\mathbf{1} \in \mathcal{X}$ such that $\mathbf{1} \odot x = x$ for all $x \in \mathcal{X}$.
    
\section{Greedy learning of preconditioned heavy-ball} \label{sec--greedy-L2O}

This section introduces our proposed method: greedy learning to optimize. We first motivate why learned preconditioners can, in principle, produce very large one-step improvements on a fixed training set. Next, define a parametrized preconditioned heavy-ball method and formulate the greedy parameter-learning problem. Finally, we present the training algorithm and describe how the learned parameters are used to optimize unseen objectives.

\subsection{Motivation: one-step improvement by learned preconditioning}

At each iteration $t \in \{0,1,2,\dots\}$, we first consider the parametrized preconditioned gradient descent update
\begin{equation} \label{eq--general-update}
    x^{t+1} = x^t - G_{\theta_t} \nabla f\left(x^t\right),
\end{equation}
where \(G_{\theta_t}\in \mathcal L(\mathcal X)\) is a learned linear operator, parametrized by \(\theta_t\in\Theta\). The following propositions show that it is possible to obtain convergence after just one iteration of the update \cref{eq--general-update}. This highlights that one can do much better than gradient descent on small datasets. Firstly, we show that it is possible to even when $G$ is a pointwise operator, i.e. $Gx:= p \odot x$ for some $p \in \mathcal{X}$.

\begin{proposition} \label{prop--imeediate-diag-conv}
    Assume that $f: \mathcal{X} \to \mathbb{R}$ is continuously differentiable and strongly convex and denote its unique global minimum by $x^*$. Then for any initial point $x^0 \in \mathcal{X}$ such that $\left[\nabla f\left(x^0\right)\right]_i \neq  0$ for all $i \in \{1, \dots, n\}$, there exists a pointwise preconditioner, such that gradient descent reaches the minimizer in one iteration.
\end{proposition}

\begin{proof}%[\textbf{Proof of \cref{prop--imeediate-diag-conv}.}]
Choose the vector $p \in \mathcal{X}$ such that $p = (x^0 - x^*) \oslash \nabla f\left(x^0\right)$. Then as $x^1 = x^0 - p \odot \nabla f\left(x^0\right)$, we have
\begin{align*}
    x^1 =  x^0 - \left(x^0 - x^*\right) \odot \nabla f\left(x^0\right) \oslash \nabla f\left(x^0\right) =   x^*.
\end{align*}
\end{proof}

% \begin{proposition} \label{prop--imeediate-diag-conv}
%     Assume that $f: \mathcal{X} \to \mathbb{R}$ is continuously differentiable and has a global minimum, and take an initial point $x^0 \in \mathcal{X}$ such that $[\nabla f\left(x^0\right)]_i \neq 0$ for all $i \in \{1, \dots, n\}$. Then there exists $p \in \mathcal{X}$ such that, $x^0 - p \odot  \nabla f\left(x^0\right) \in \arg\min_x f(x)$.
% \end{proposition}

% \begin{proof}%[\textbf{Proof of \cref{prop--imeediate-diag-conv}.}]
% Take any $x^* \in \arg\min_x f(x)$. Choose the vector $p \in \mathcal{X}$ such that
% \begin{equation}
%     p_i = \frac{[x^0 - x_*]_i}{[\nabla f\left(x^0\right)]_i},
% \end{equation}
% then for each $i \in \{1, \dots, n\}$, we have
% \begin{align*}
%     [x^0 - p \odot \nabla f\left(x^0\right)]_i =  [x^0]_i - \frac{[x^0 - x^*]_i}{[\nabla f\left(x^0\right)]_i} [\nabla f\left(x^0\right)]_i =   x^*_i.
% \end{align*}
% \end{proof}

% \begin{remark}
%     Therefore, learning an update $x_k^1  = x_k^0 - p_{\theta}( w_k^0) \odot \nabla f_k(x_k^0)$ can be seen as equivalent to learning a solution map $x_k^1 = H_{\theta}( w_k^0)$.
% \end{remark}

While the pointwise parametrization may obtain convergence after one iteration for one function, for an arbitrary linear operator $G \in \mathcal{L}(\mathcal{X})$, under certain conditions, one can obtain convergence after one iteration for multiple functions.
\begin{proposition} \label{prop--immediate-full-conv}
    For $k \in \{1,\dots, N\}$, assume that $f_k: \mathcal{X} \to \mathbb{R}$ is continuously differentiable, and has a global minimum, with any initial point $x_k^0 \in \mathcal{X}$. Assume that the set of gradients $\{ \nabla f_1(x_1^0), \dots,  \nabla f_N(x_N^0)\}$ is linearly independent. Then if $N \leq n$, there exists an operator $P \in \mathcal{L}(\mathcal{X})$ such that $x_k^0 - P  \nabla f_k(x_k^0) \in \arg\min_x f_k(x)$, for all $k \in \{1,\dots,N\}$. 
\end{proposition}

\begin{proof}%[\textbf{Proof of \cref{prop--immediate-full-conv}.}]
    Take any $x^* \in \arg\min_x f(x)$. We wish to find a linear operator $P \in \mathcal{L}(\mathcal{X})$ such that $x_k^* = x_k^0 - P \nabla f_k(x_k^0)$ for $k \in \{1, \dots, N\}$.
% \begin{align*}
% \begin{cases}
%     x_1^* = x_1^0 - P \nabla f_1(x_1^0), \\
%     \vdots \\
%     x_N^* = x_N^0 - P \nabla f_N(x_N^0).
% \end{cases}
% \end{align*}
This gives $nN = N\dim (\mathcal{X})$ linear equations in $n^2$ unknowns. Rewritten, these read
\begin{equation}
    P \begin{bmatrix}
         \nabla f_1(x_1^0) | \dots | \nabla f_N(x_N^0)
    \end{bmatrix} = \begin{bmatrix}
         x_1^0 - x_1^* | \dots | x_N^0 - x_N^*
    \end{bmatrix}.
\end{equation}
As the columns of $\begin{bmatrix}
         \nabla f_1(x_1^0) | \dots | \nabla f_N(x_N^0)
    \end{bmatrix}$ are linearly independent, such a $P$ exists if $nN \leq n^2$, which is equivalent to $N \leq n$.
% Note in the case that $N = n$, we have a unique choice of $P$:
% \begin{equation}
%     P = \begin{bmatrix}
%          \nabla f_1(x_1^0) | \dots | \nabla f_N(x_N^0)
%     \end{bmatrix}^{-1} \begin{bmatrix}
%          x_1^0 - x_1^* | \dots | x_N^0 - x_N^*
%     \end{bmatrix}
% \end{equation}
% \comment{Not complete}
\end{proof}

\cref{prop--imeediate-diag-conv} and \cref{prop--immediate-full-conv} motivate learning $\theta_t$ by considering the function values only at the next iteration, due to the possibility of convergence after one iteration. 

\subsection{Parametrized preconditioned heavy-ball}

We extend \cref{eq--general-update} by adding a learned momentum preconditioner.  Let \(\Theta\) and \(\Phi\) be Hilbert spaces of parameters. We learn parameters $\theta_t \in \Theta$ and $\phi_t \in \Phi$ in the update 
\begin{equation} \label{eq--general-update-hb}
    x^{t+1} = x^t - G_{\theta_t} \nabla f\left(x^t\right) + H_{\phi_t} (x^t - x^{t-1}),
\end{equation}
where $G_{\theta_t}, H_{\phi_t} \in \mathcal{L}(\mathcal{X})$. In order to learn the parameters $\theta_t, \phi_t$ for $t \in \{0,1,2,\dots\}$, we introduce a training set of objective functions $\mathcal{T}:= \{f_1, \dots, f_N\}$, with $f_k \in \mathcal{F}_{L_k}$ for $k \in \{1, \dots, N\}$, with corresponding initial points $\mathcal{X}_0:= \{x_1^0, \dots, x_N^0\}$, and define $f_k^* = \min_x f_k(x)$. The following definition provides a condition for which the update rule \cref{eq--general-update} generalizes gradient descent (GD).
\begin{definition} \label{eq--GGD}
    We say that the family $(G_\theta, H_{\phi})$ is \textit{GGD} if the family generalizes gradient descent, meaning there exist parameters $\tilde{\phi}$ such that $H_{\tilde{\phi}} = 0$, and for all $\alpha \ge 0$, there exist parameters $\theta$ such that $G_{\theta} = \alpha I$.
\end{definition}

Examples of parametrizations that satisfy the GGD property are shown in \cref{section--linear-params}. Let $\tau = 1 / {L_{\text{train}}}$, where ${L_{\text{train}}} \geq \max \{ L_1, \dots, L_N \}$ upper bounds the largest smoothness coefficient in the training data. This choice of step size in gradient descent ensures convergence for all functions $f_k \in \mathcal{T}$. From this point forward, we assume $(G_\theta, H_{\phi})$ is GGD, meaning in particular that there exists $\tilde{\theta}, \tilde{\phi}$ such that $G_{\tilde{\theta}} = \tau I, H_{\tilde{\phi}} = 0$.  Furthermore, the GGD property can be leveraged to establish provable convergence for a set of unseen functions by introducing a penalty when the parameters $(\theta_t, \phi_t)$ deviate significantly from $(\tilde{\theta}, \tilde{\phi})$. With this purpose, at iteration $t$, defining $x_k^{t+1} := x_k^t - G_{\theta}\nabla f_k(x_k^t) + H_{\phi} (x_k^t - x_k^{t-1})$, which depends on parameters $\theta$ and $\phi$, we solve the optimization problem
\begin{align} \label{eq--opt-problem}
    (\theta_{t}, \phi_t) \in \arg\min_{\theta, \phi} \left \{ g_{t, \lambda_t, \mu_t}(\theta, \phi) := \frac1N \sum_{k=1}^N f_k(x_k^{t+1}) + \frac {\lambda_t} 2 \| \theta - \tilde{\theta} \|^2 + \frac {\mu_t} 2 \| \phi - \tilde{\phi} \|^2  \right\} ,
\end{align}
for some regularization parameters $\lambda_t, \mu_t \geq 0$, which are used to balance the importance of the regularizers. Such a strategy is greedy, as learning refers to tuning the parameters $\theta_t$ and $\phi_t$ considering only the function values at the next iteration. The sequential training procedure for parameter learning is detailed in \cref{alg--train-algo}. For unrolling with a standard implementation of backpropagation, device memory requirements scale linearly with the number of training iterations. However, with our greedy method, once the parameters $\theta_t, \phi_t$ and the next iterates $x_k^{t+1}$ for $k \in \{1, \dots, N\}$ have been calculated, $\theta_t$ and $\phi_t$ are no longer required to be stored on the device. Therefore, device memory is constant with increasing training iterations for our greedy method. In addition, the unroll length $T$ is not required to be chosen apriori. For example one may terminate the training procedure if the condition $\| \nabla f(x_k^t) \|_2 < \xi$ is satisfied for all $k \in \{1, \dots, N\}$ for some chosen tolerance $\xi$. To minimize an unseen function $f$ with initial point $x^0$, we propose \cref{alg--test-algo}.

% Unlike unrolling, where memory requirements scale linearly with the number of training iterations, greedy training is a memory-friendly method, with memory that is constant with increasing training iterations.

\begin{algorithm}[H]
\caption{Training algorithm for greedy parameter learning}\label{alg--train-algo}
\begin{algorithmic}[1]
\State \textbf{Input:} Functions $f_1, \dots, f_N$, initial points $x_1^0, \dots, x_N^0$, $x_1^{-1}:=x_1^{0}, \dots, x_N^{-1}:=x_N^{0}$, tolerance $\xi$, regularization parameter sequences $(\lambda_t)_{t\geq 0}$ and $(\mu_t)_{t\geq 0}$ with $\lambda_t,\mu_t\geq 0$.
\For{$t = 0, 1, 2, \dots$}
    \State Find $\theta_t, \phi_t$ \text{ s.t. } $g_{t, \lambda_t, \mu_t}(\theta_t, \phi_t) < g_{t, \lambda_t, \mu_t}(\tilde{\theta}, \tilde{\phi})$
    \For{$k = 1, 2, \dots, N$}
        \State $x_k^{t+1} = x_k^t - G_{\theta_t}\nabla f_k(x_k^t) + H_{\phi_t}(x_k^t - x_k^{t-1})$
    \EndFor
     \If{$ \|\nabla f_k(x_k^t)\|_2 < \xi$ \text{ for all } $k \in \{1, \dots, N\}$}
    \State $T = t$
        \State \textbf{break}
    \EndIf
\EndFor
\State \textbf{Output:} Learned parameters $(\theta_0, \phi_0), \dots, (\theta_{T}, \phi_{T})$.
\end{algorithmic}
\end{algorithm}

\begin{algorithm}[H]

\caption{Learned algorithm to minimize a function $f$}\label{alg--test-algo}
\begin{algorithmic}[1]
\State \textbf{Input:} Function $f$ with initial point $x^0$, $x^{-1}:=x^0$, tolerance $\xi$, learned parameters $(\theta_0, \phi_0), \dots, (\theta_{T}, \phi_{T})$.
\For{$t = 0, 1, 2, \dots$}
    \If{$t \leq T$}
    \State $x^{t+1} = x^t - G_{\theta_t}\nabla f\left(x^t\right) + H_{\phi_t}(x^t - x^{t-1})$
    \Else
    \State $x^{t+1} = x^t - G_{\theta_{T}}\nabla f\left(x^t\right) + H_{\phi_T}(x^t - x^{t-1})$
    \EndIf
    \If{$ \|\nabla f\left(x^t\right)\|_2 < \xi$}
        \State \textbf{break}
    \EndIf
\EndFor
\State \textbf{Output:} $x^{t+1}$.
\end{algorithmic}
\end{algorithm}

\begin{remark}[Lyapunov-based greedy training]
The greedy training strategy proposed in this paper is aligned with algorithms for which one can prove a one-step decrease inequality in a quantity that is directly computable from the training data. In the present work this quantity is the next-step objective value, which motivates minimizing $f_k(x_k^{t+1})$. More generally, one could consider greedy training based on a discrete Lyapunov (or potential) function. Such arguments are already standard in the analysis of several discrete optimization algorithms, including Nesterov-type accelerated methods \cite{wilson2021lyapunov}, accelerated proximal-gradient methods like FISTA \cite{beck2009fast}, and accelerated primal-dual algorithms \cite{condat2026nesterov}. This suggests that, in principle, one could replace the objective-based greedy loss by the next-step Lyapunov function value. We leave this direction for future work.

%However, unlike the loss used in the present paper, such Lyapunov functions typically depend on quantities such as $x^\star$, which are not exactly available in practice. One could instead use high-accuracy approximate reference solutions, but the resulting training objective and convergence analysis would then require explicit control of the approximation error.  We therefore leave this direction for future work.

% A further issue is convexity of the parameter-learning problem. For example, in NAG, if one learns preconditioners for both the gradient step and the extrapolation step the parameters enter the look-ahead point and hence the gradient evaluation through composition. In that case, even when the training functions $f_k$ are convex, the resulting parameter-learning problem is not necessarily convex. However, if the momentum coefficients $(\beta_t)$ are kept as a handcrafted sequence and one learns only the gradient-step operators through linear parametrizations, then the next iterate remains affine in the learned parameters and hence the greedy learning subproblem remains convex for convex training functions. However, such theory is beyond the scope of this paper.
\end{remark}

\section{Convergence results} \label{section--convergence}

This section contains convergence results for \cref{alg--test-algo}. Firstly, we show that the learned method converges on the training set whenever each greedy update performs at least as well as a reference gradient-descent step. This result does not require regularization, nor for the learned preconditioners to be necessarily symmetric or positive definite. Following this, we show that convergence can be guaranteed on unseen functions by adding soft constraints that keep the learned method close to a safe gradient-descent reference method.

Before we present the convergence results, we require the following definition, which ensures the parametrized algorithm adopts the convergence properties of gradient descent. 
\begin{definition}
    We say that $(\theta_t, \phi_t)$ is \textit{BGD} (better than gradient descent) with regularization parameters ${\lambda_t, \mu_t}$ if
    \begin{equation} \label{def--BGD}
        g_{t, {\lambda_t}, \mu_t}(\theta_t, \phi_t) \leq g_{t, {\lambda_t}, \mu_t}(\tilde{\theta}, \tilde{\phi}) = \frac1N \sum_{k=1}^N f_k \left(x_k^t - \tau \nabla f_k (x_k^t)\right).
    \end{equation}
\end{definition}
In \cref{section--linear-params} we introduce parameterizations $G_{\theta}$ for which the BGD property is easily obtained during training.

\subsection{Convergence on training data}

The following lemma is required to prove the convergence of our learned method, with a proof provided in \cref{app:convergence}.
\begin{lemma} \label{prop:F-properties}
    Define $F:  \mathcal{X}^N \to \mathbb{R} $ by $F(\mathbf{x}) = 1/N \sum_{k=1}^N f_k(x_k),
        \mathbf{x} = (x_1, x_2, \dots, x_N)\in \mathcal{X}^N$. If each $f_k \in \mathcal{F}_{L_k}$ then $F \in \mathcal{F}_{L_F}$, with $L_F = {L_{\text{train}}}/N$.
\end{lemma}
The function $F$ denotes the average of the training functions $f_k$. The following theorem proves convergence of $\nabla F$, and therefore each $\nabla f_k$.
\begin{theorem}[Convergence on training data] \label{prop:conv}
    Suppose that ${\lambda_t}, \mu_t \geq 0$ and  $(\theta_t, \phi_t)_{t=0}^\infty$ is a BGD sequence of parameters. Then with \cref{alg--train-algo}, we have $\nabla f_k (x_k^t) \to 0$ as $t \to \infty$ for all $k \in \{1, \dots, N\}$.
\end{theorem}

\begin{proof}%[\textbf{Proof of \cref{prop:conv}}]
As $\theta_t$ is BGD, we have that
\begin{align*}
    F(x^{t+1}) \le g_{t, {\lambda_t, \mu_t}}(\theta_t, \phi_t) &\leq g_{t, {\lambda_t}, \mu_t} (\tilde{\theta}, \tilde{\phi}) 
    =\frac1N \sum_{k=1}^N f_k\left( x_k^t - \tau \nabla f_k(x_k^t)\right) 
    %&= F\left(x^t - \tau \left(\nabla f_1(x_1^t), \dots, \nabla f_N(x_N^t)\right)\right)\\
    = F\left(x^t - \tau N \nabla F\left(x^t\right)\right).
    %&= F\left(x^t - \tau_F \nabla F\left(x^t\right)\right),
\end{align*}
    Note that $\tau N =  1/{L_F}$, then using the $L_F$-smoothness of $F$ we have sufficient descent:
    \begin{align*}
        F(x^{t+1}) & \leq F\left(x^t - \frac1{L_{F}} \nabla F\left(x^t\right)\right)\leq F\left(x^t\right) - \frac1{2L_F} \| \nabla F\left(x^t\right) \|^2.
    \end{align*}
    Now, using \cref{lemma::convergence_grad_fk}, we have that $\nabla f_k (x_k^t) \to 0$ for $t \to \infty$ for all $k \in \{1, \dots, N\}$
\end{proof}
Note that in particular, this means that convergence in training is obtained even when $\lambda_t, \mu_t = 0$ for all $t$. In particular, the learned preconditioners $G_t$ are never necessarily positive-definite.

\begin{remark}
    In the case when $\lambda_t = \mu_t = 0$, if the Lipschitz constants $L_k$ are unknown, then instead of comparing with the objective decrease of $F$ using a step size of $1/L_F$, one may prove convergence on training data by obtaining sufficient descent instead by comparing with any step size $\tau_t$ leading to a sufficient decrease in $F$ for convergence, for example, one found using backtracking line search on the function $F$. Details are omitted for brevity. In this case, as $\lambda_t = \mu_t = 0$, convergence on test data is not proved.
\end{remark}

% \rcomment{You might want to say more about this if you consider it worth a remark. It seems out of place.}
% \begin{remark}
%     This result can be generalized for a more general optimization algorithm. Often, people use neural networks to parametrize optimization algrotihms, for example, they consider the update $x^{t+1} = K_{\theta_t} (w_t)$, for a parametrized function $K_{\theta_t}: \mathcal{W} \to \mathcal{X}$. Such a function often contains of a large number of parameters, and therefore memory issues quickly become an issue in training. This can be remedied using checkpointing \cite{checkpointing2016}, however if unrolling to a large number of iterations, this makes training prohibitively expensive. \rcomment{Say time-wise. Maybe make this a bigger point saying it the first time memory is discussed?}
% \end{remark}

\subsection{Convergence on test data} \label{subsec:conv-test-data}

We now show convergence on test data. Firstly, the following assumption is required.
\begin{assumption} \label{assumption::liminfbgdcts}
     Let $\lim \inf_{t \to \infty} \lambda_t > 0$, $\lim \inf_{t \to \infty} \mu_t > 0$, $(\theta_t, \phi_t)_{t=0}^\infty$ be BGD, and $G : \Theta \to \mathcal{L}(\mathcal{X})$ and $H : \Phi \to \mathcal{L}(\mathcal{X})$ be continuous.
\end{assumption}

In practice, the simplest choice is to take $\lambda_t=\mu_t=0$, since this removes the regularizers in parameter learning and typically gives the best empirical performance. The role of regularization is instead theoretical: choosing $\lambda_t$ and $\mu_t$ to be small but eventually bounded away from zero ensures that the learned parameters remain close to the safe reference values $(\tilde{\theta},\tilde{\phi})$, which is what allows us to prove convergence on unseen functions. Thus, there is a trade-off between empirical flexibility and provable safety. Designing a practically optimal schedule for $(\lambda_t,\mu_t)$ is beyond the scope of this paper.

We now show that if the regularization parameters $(\lambda_t, \mu_t)$ are eventually non-vanishing, then the learned parameters tend towards $(\tilde{\theta}, \tilde{\phi})$.
\begin{lemma} [Convergence of parameters] \label{lemma--theta-t-conv} 
    Let \cref{assumption::liminfbgdcts} hold. Then $(\theta_t, \phi_t) \to (\tilde{\theta}, \tilde{\phi})$ as $t \to \infty$, and therefore $G_{\theta_t} \to \tau I$ and $H_{\phi_t} \to 0 \text{ as } t \to \infty$.
\end{lemma}

\begin{proof}
First, note that
\begin{align*}
    \frac{\lambda_t}2 \| \theta_t - \tilde{\theta} \|^2 + \frac{\mu_t}2 \| \phi_t - \tilde{\phi}\|^2 + F(x^{t+1}) = g_{t, \lambda_t, \mu_t}(\theta_t, \phi_t) \le g_{t, 0, 0}(\tilde{\theta}, \tilde{\phi}) \leq  F\left(x^t\right),
\end{align*}
by the BGD and GGD properties. The assumptions $\lim\inf_{t \to \infty} \lambda_t, \mu_t > 0$ means that there exist $\lambda, \mu, K_1 > 0$ such that $\lambda_t \geq \lambda$ and $\mu_t \geq \mu$ for all $t \ge K_1$. Then for all $t \ge K_1$, we have $\lambda \| \theta_t - \tilde{\theta} \|^2 + \mu \| \phi_t - \tilde{\phi}\|^2\leq 2(F\left(x^t\right) - F(x^{t+1}))$. Taking a summation up to $K_2 > K_1$ gives
\begin{align}
   \sum_{t=K_1}^{K_2}\lambda \| \theta_t - \tilde{\theta} \|^2 + \mu \| \phi_t - \tilde{\phi}\|^2\leq 2(F(x_{K_1}) - F(x_{K_2+1})) \leq 2(F(x_0) - F^*),
\end{align}
where the final inequality is due to $F(x_{K_1}) \le F(x_0)$, and $F(x_{K_2+1}) \ge F^*$. In particular, the left-hand summations are bounded above by a constant in $t$. Therefore, as the series converge, taking $K_2 \to \infty$ requires $\theta_t \to \tilde{\theta}$ and $\phi_t \to \tilde{\phi}$ as $t \to \infty$. By continuity of $G$ and $H$, $G_{\theta_t} \to \tau I$ and $H_{\phi_t} \to 0 \text{ as } t \to \infty$.
\end{proof}

The idea is that we start with a method that fits the data very well, leading to quick initial convergence, but in the interest of safety, over time, we become closer to an algorithm with proven convergence, in particular, with $G_{\theta_t}$ positive-definite eventually. Note that one could alternatively switch to plain gradient descent after some iteration $T$, which would trivially recover the standard convergence rates. We instead keep $(\theta_T,\phi_T)$ fixed after iteration $T$, since this may retain some advantage over exact gradient descent while still remaining sufficiently close to the safe limit for the convergence analysis to apply. The following lemma gives an upper bound on the descent of a test function $f$ using \cref{alg--test-algo}. This result is used to prove convergence on test data.

\begin{lemma}[Descent of perturbed heavy-ball] \label{lemma:thm-calcs}
    Suppose $\| G_{\theta_T} - \tau I \| \le \varepsilon_1$, and $\| H_{\phi_T} \| \le \varepsilon_2$ for some $\varepsilon_1, \varepsilon_2 \ge 0$. Define constants $C_1, C_2$ by
    \begin{align}
        C_1 &= \tau \left( 1 - \frac{\tau L}{2} \right) - \left(\varepsilon_1 + \frac {\varepsilon_2} 2\right) |\tau L - 1| -   \frac L 2 \varepsilon_1 \left(  \varepsilon_1 + \varepsilon_2  \right) \\
        C_2 &= \frac{\varepsilon_2} 2\left( |\tau L - 1| + L (\varepsilon_1  + \varepsilon_2) \right).
    \end{align}
    Then using \cref{alg--test-algo} for  $f \in \mathcal{F}_L$ for some $L > 0$, for all $t > T$, 
    \begin{equation} \label{eq::general-ineq-errors}
    \begin{aligned}
        f\left(x^{t+1}\right) 
        %& \le f\left(x^t\right) - \left( \tau \left( 1 - \frac{\tau L}{2} \right) - \varepsilon_1 |\tau L - 1| -  \frac {L \varepsilon_1^2} 2 - \frac {\varepsilon_2 \left( |\tau L - 1| + L \varepsilon_1 \right)} 2  \right) \|r^t\|^2  \\
        %&+ \frac{\varepsilon_2} 2\left( |\tau L - 1| + L (\varepsilon_1  + \varepsilon_2) \right) \| p^t \|^2 \\
        & \le f\left(x^t\right) - C_1 \|\nabla f\left(x^t\right)\|^2  + C_2 \| x^t - x^{t-1} \|^2.
    \end{aligned}
    \end{equation}
\end{lemma}
\begin{proof}
 In the following define $r^t := \nabla f\left(x^t\right)$, $p^t := x^t - x^{t-1}$,  and $M : = G_{\theta_T} - \tau I$, then $\| M \| \leq \varepsilon_1$. Firstly, by $L$-smoothness of $f$, we have for $t > T$,
\begin{align*}
    f\left(x^{t+1}\right) &\le f\left(x^t\right) + \langle r^t, p^{t+1} \rangle + \frac {L} 2 \| p^{t+1} \|^2 =  f\left(x^t\right) + \left\langle p^{t+1}, r^t + \frac {L} 2 p^{t+1} \right\rangle \\
    &=  f\left(x^t\right) + \left\langle -(\tau I + M)r^t + H_{\phi_T} p^t, r^t + \frac {L} 2 (-(\tau I + M)r^t + H_{\phi_T} p^t) \right\rangle \\
    %&=  f\left(x^t\right) + \left\langle -\tau r^t - Mr^t + H_{\phi_T} p^t, \left( 1- \frac{\tau L} 2 \right) r^t - \frac L 2 M r^t + \frac {L} 2H_{\phi_T} p^t \right\rangle \\
    &=  f\left(x^t\right) - \tau \left( 1 - \frac{\tau L}{2} \right) \|r^t\|^2 + (\tau L - 1) \langle Mr^t, r^t \rangle + (1 - \tau L) \langle r^t, H_{\phi_T}p^t \rangle \\
    &+ \frac L 2 \| M r^t \|^2 - L \langle Mr^t, H_{\phi_T}p^t \rangle + \frac L 2 \| H_{\phi_T}p^t \|^2 \\
    & \leq f\left(x^t\right) - \left( \tau \left( 1 - \frac{\tau L}{2} \right) - \varepsilon_1 |\tau L - 1| -  \frac {L \varepsilon_1^2} 2 \right) \|r^t\|^2  \\
    &+ \varepsilon_2 \left( |\tau L - 1| + L \varepsilon_1 \right) \| r^t \| \|p^t \| + \frac {L\varepsilon_2^2} 2 \| p^t \|^2 .
\end{align*}
Young's inequality states that for any $a, b \in \mathbb{R}$ that $2ab \le  a^2 + b^2$. Using this gives \\ $\varepsilon_2 \left( |\tau L - 1| + L \varepsilon_1 \right) \| r^t \| \|p^t \| \leq \frac{\varepsilon_2 \left( |\tau L - 1| + L \varepsilon_1 \right)} 2 \|r^t\|^2 + \frac {\varepsilon_2 \left( |\tau L - 1| + L \varepsilon_1 \right)} 2 \| p^t \|^2$ for $t > T$, and \cref{eq::general-ineq-errors} follows.
% Now, by $L$-smoothness of $f$ and Young's inequality, we have
% \begin{align}
%     f\left(x^{t+1}\right) &\le f\left(x^t\right) + \langle r^t, p^{t+1} \rangle + \frac L 2 \| p^{t+1} \|^2 \\
%     &\le  f\left(x^t\right) + \left(- \tau + \varepsilon_1 + \frac {\varepsilon_2} 2 \right) \|r^t \|^2  + \frac {\varepsilon_2} 2 \|  p^t \|^2 \\
%     &+ \frac L 2 \left( (\tau + \varepsilon_1)\left(\varepsilon_1 +  \varepsilon_2 + \tau\right) \| r^t \|^2 + \varepsilon_2\left( \varepsilon_1 + \varepsilon_2 + \tau \right) \|p^t \|^2  \right)\\
%     &= f\left(x^t\right) + \left( - \tau + \varepsilon_1 + \frac {\varepsilon_2} 2 + \frac L 2  (\tau + \varepsilon_1)\left(\varepsilon_1 +  \varepsilon_2 + \tau \right) \right) \|r^t \|^2 \\
%     & + \frac {\varepsilon_2} 2 \left( 1 + L( \varepsilon_1 + \varepsilon_2 + \tau ) \right) \|  p^t \|^2.
% \end{align}
\end{proof}

Using this lemma, we present our main Theorem, which provides convergence of test functions.
\begin{theorem}[Convergence on test data] \label{thm--reg-conv}
Let \cref{assumption::liminfbgdcts} hold. Then, for any $L < 2 L_{\text{train}}$, there exists a final training iteration $T$ such that for all $f \in \mathcal{F}_{L }$, and any starting point $x^0$, using \cref{alg--test-algo} we have $\nabla f\left(x^t\right) \rightarrow 0$ as $t \rightarrow \infty$. Furthermore, there exists a constant $C>0$, such that
    \begin{equation} \label{eq:nonconvex-conv-rate}
        \min_{s \in  \{0, 1, \dots, t\}} \| \nabla f\left(x^{s}\right)\|^2 \leq  \frac C {t}.
    \end{equation}
\end{theorem}
\begin{proof}
Fix $\varepsilon_1, \varepsilon_2 \ge 0$. $G_{\theta_t} \to \tau I$ and $H_{\phi_t} \to 0$ as $t \to \infty$ implies that there exists a final training iteration $T$ such that $\| G_{\theta_T} - \tau I \| \le \varepsilon_1 $  and $\| H_{\phi_T} \| \le \varepsilon_2 $. Again define $r^t := \nabla f\left(x^t\right)$, $p^t := x^t - x^{t-1}$, and $M : = G_{\theta_T} - \tau I$. For some $\gamma > 0$, and for $t > T$, define
\begin{equation} \label{eq::psi-definition}
    \Psi_t := f\left(x^t\right) - f^* + \gamma  \| p^t \|^2.
\end{equation}
Note that $p^{t+1} = x^{t+1} - x^t = -(\tau I + M)r^t + H_{\phi_T} p^t$, and so for $t > T$,
\begin{align*}
    \| p^{t+1} \|^2 &= \| -(\tau I + M)r^t + H_{\phi_T} p^t \|^2 
    = \| (\tau I + M)r^t \|^2 + \|H_{\phi_T} p^t \|^2 - 2 \langle (\tau I + M)r^t, H_{\phi_T}p^t \rangle \\
    &\le (\tau + \varepsilon_1)^2 \| r^t \|^2 + \varepsilon_2^2 \|p^t \|^2 + 2(\tau + \varepsilon_1) \varepsilon_2 \| r^t \| \|p^t \|.
\end{align*}
Young's inequality gives $2(\tau + \varepsilon_1) \varepsilon_2 \| r^t \| \|p^t \| \le (\tau + \varepsilon_1) \varepsilon_2 \| r^t \|^2 + (\tau + \varepsilon_1) \varepsilon_2 \|p^t \|^2$, and so 
\begin{equation} \label{eq:ptp1-calc}
    \| p^{t+1} \|^2 \le  (\tau + \varepsilon_1)\left(\varepsilon_1 +  \varepsilon_2 + \tau\right) \| r^t \|^2 + \varepsilon_2\left( \varepsilon_1 + \varepsilon_2 + \tau \right) \|p^t \|^2, 
\end{equation}
for $t > T$. Then, using \cref{lemma:thm-calcs} and \cref{eq:ptp1-calc},
\begin{align*}
    &\Psi_{t+1} - \Psi_t = f\left(x^{t+1}\right) - f\left(x^t\right) + \gamma \left[ \| p^{t+1} \|^2 - \| p^t \|^2 \right] \\
    & \leq  - \left( \tau \left( 1 - \frac{\tau L}{2} \right) - \left(\varepsilon_1 + \frac {\varepsilon_2} 2\right) |\tau L - 1| -   \frac L 2 \varepsilon_1 \left(  \varepsilon_1 + \varepsilon_2  \right) - \gamma (\tau + \varepsilon_1)\left(\varepsilon_1 +  \varepsilon_2 + \tau\right)  \right) \|r^t\|^2  \\
    &-  \left(\gamma \left(1 -  \varepsilon_2\left( \varepsilon_1 + \varepsilon_2 + \tau \right)\right) - \frac{\varepsilon_2} 2\left( |\tau L - 1| + L (\varepsilon_1  + \varepsilon_2) \right) \right) \| p^t \|^2.
\end{align*}
% Then 
% \begin{align*}
%     \Psi_{t+1} - \Psi_t &= f\left(x^{t+1}\right) - f\left(x^t\right) + \gamma \left[ \| p^{t+1} \|^2 - \| p^t \|^2 \right] \\
%     & \leq \left( - \tau + \varepsilon_1 + \frac {\varepsilon_2} 2 + \frac L 2  (\tau + \varepsilon_1)\left(\varepsilon_1 +  \varepsilon_2 + \tau \right) + \gamma  (\tau + \varepsilon_1)\left(\varepsilon_1 +  \varepsilon_2 + \tau\right) \right) \|r^t \|^2 \\
%     & + \left (\frac {\varepsilon_2} 2 ( 1 + L( \varepsilon_1 + \varepsilon_2 + \tau)) - \gamma \left( 1 -  \varepsilon_2\left( \varepsilon_1 + \varepsilon_2 + \tau \right) \right) \right) \|  p^t \|^2 \\
% \end{align*}
Define
% \begin{align*}
%     C_1(\varepsilon_1, \varepsilon_2, \gamma) &=  \tau - \varepsilon_1 - \frac {\varepsilon_2} 2 - \frac L 2  (\tau + \varepsilon_1)\left(\varepsilon_1 +  \varepsilon_2 + \tau \right) - \gamma  (\tau + \varepsilon_1)\left(\varepsilon_1 +  \varepsilon_2 + \tau \right)\\
%     C_2(\varepsilon_1, \varepsilon_2, \gamma) &=  \gamma \left( 1 -  \varepsilon_2\left( \varepsilon_1 + \varepsilon_2 + \tau \right) \right) -\frac {\varepsilon_2} 2 ( 1 + L( \varepsilon_1 + \varepsilon_2 + \tau)),
% \end{align*}
\begin{align}
    C_1(\varepsilon_1, \varepsilon_2, \gamma) &:=  \tau \left( 1 - \frac{\tau L}{2} \right) - \left(\varepsilon_1 + \frac {\varepsilon_2} 2\right) |\tau L - 1| -   \frac L 2 \varepsilon_1 \left(  \varepsilon_1 + \varepsilon_2  \right) - \gamma (\tau + \varepsilon_1)\left(\varepsilon_1 +  \varepsilon_2 + \tau\right) \label{eq::C1}, \\
    C_2(\varepsilon_1, \varepsilon_2, \gamma) &:= \gamma \left(1 -  \varepsilon_2\left( \varepsilon_1 + \varepsilon_2 + \tau \right)\right) - \frac{\varepsilon_2} 2\left( |\tau L - 1| + L (\varepsilon_1  + \varepsilon_2) \right) \label{eq::C2},
\end{align}
then
\begin{equation} \label{eq::psi-decrease}
    \Psi_{t+1} - \Psi_t \leq - C_1(\varepsilon_1, \varepsilon_2, \gamma) \| r^t \|^2 -  C_2(\varepsilon_1, \varepsilon_2, \gamma) \|p^t \|^2.
\end{equation}
We would like to choose $\varepsilon_1, \varepsilon_2, \gamma$ such that $C_1(\varepsilon_1, \varepsilon_2, \gamma), C_2(\varepsilon_1, \varepsilon_2, \gamma) > 0$. This means that 
% \begin{equation}
%    \frac{\varepsilon_2  ( 1 + L( \varepsilon_1 + \varepsilon_2 + \tau))}{2\left( 1 -  \varepsilon_2\left( \varepsilon_1 + \varepsilon_2 + \tau \right) \right)} \le \gamma < \frac{2\tau - 2\varepsilon_1 - \varepsilon_2 - L  (\tau + \varepsilon_1)\left(\varepsilon_1 +  \varepsilon_2 + \tau \right)}{2(\tau + \varepsilon_1)(\tau + \varepsilon_1 + \varepsilon_2)}
% \end{equation}
\begin{equation*} \label{eq::gamma}
   \frac{ \frac{\varepsilon_2} 2\left( |\tau L - 1| + L (\varepsilon_1  + \varepsilon_2) \right)}{1 -  \varepsilon_2\left( \varepsilon_1 + \varepsilon_2 + \tau \right)} < \gamma < \frac{\tau \left( 1 - \frac{\tau L}{2} \right) - \left(\varepsilon_1 + \frac {\varepsilon_2} 2\right) |\tau L - 1| -   \frac L 2 \varepsilon_1 \left(  \varepsilon_1 + \varepsilon_2  \right)}{ (\tau + \varepsilon_1)(\tau + \varepsilon_1 + \varepsilon_2)}.
\end{equation*}
This is equivalent to requiring
\begin{equation} \label{eq::P}
\begin{aligned}
    P(\varepsilon_1, \varepsilon_2) &:=  \frac{\varepsilon_2} 2\left( |\tau L - 1| + L (\varepsilon_1  + \varepsilon_2) \right) (\tau + \varepsilon_1)(\tau + \varepsilon_1 + \varepsilon_2)  \\
    &- \left(\tau \left( 1 - \frac{\tau L}{2} \right) - \left(\varepsilon_1 + \frac {\varepsilon_2} 2\right) |\tau L - 1| -   \frac L 2 \varepsilon_1 \left(  \varepsilon_1 + \varepsilon_2  \right)\right)\left(1 -  \varepsilon_2\left( \varepsilon_1 + \varepsilon_2 + \tau \right)\right)<0.
\end{aligned}
\end{equation}
Note that $P(0, 0) = -\tau(1- \frac{\tau L}{2}) < 0$ if $L < 2L_{\text{train}}$. As $P$ is continuous in $\varepsilon_1$ and $\varepsilon_2$, there exist $(\tilde{\varepsilon_1}, \tilde{\varepsilon_2})$ such that if $L < 2L_{\text{train}}$ then $P(\varepsilon_1, \varepsilon_2) < 0$ for all $\varepsilon_i \in [0, \tilde{\varepsilon_i})$, $i \in \{1, 2\}$. Therefore, for $\varepsilon_1, \varepsilon_2$ sufficiently small, the interval for $\gamma$ such that $C_1, C_2 > 0$ is nonempty. Now, note that $\Psi_t \geq 0$ for all $t$ as $\gamma \geq 0$, and $f\left(x^t\right) \ge f^*$ for all $t$. Therefore,
\begin{align*}
    \Psi_{t+1} - \Psi_t \leq - C_1(\varepsilon_1, \varepsilon_2, \gamma) \| r^t \|^2 -  C_2(\varepsilon_1, \varepsilon_2, \gamma) \|p^t \|^2 \leq - C_1(\varepsilon_1, \varepsilon_2, \gamma) \| r^t \|^2,
\end{align*}
and so by taking a summation, 
\begin{align*}
    &- \Psi_T \le\Psi_{t+1} - \Psi_T = \sum_{s=T}^t (\Psi_{s+1} - \Psi_s) \leq - C_1 \sum_{s=T}^t \| r^s \|^2 \implies \sum_{s=T}^t \| r^s \|^2 \leq \frac{\Psi_T}{C_1},
\end{align*}
and by taking the limit $t \to \infty$ we see that $r^s \to 0$ as $s \to \infty$. Finally, let $C(\varepsilon_1, \varepsilon_2, \gamma) := \Psi_T/C_1(\varepsilon_1, \varepsilon_2, \gamma)$. As $\nabla f$ is Lipschitz, there exists a constant $C_3>0$ such that $\| \nabla f\left(x^{s}\right) \|^2 \leq C_3$ for all $s \in \{0, \dots, T-1\}$, and so
\begin{align*}
    \min_{s=0, \dots, t}\|\nabla f\left(x^s\right)\|^2 &\leq \frac1 t \sum_{s=0}^t \| \nabla f\left(x^{s}\right) \|^2  
    %\leq \frac1 t \sum_{s=T}^t \| \nabla f\left(x^{s}\right) \|^2 + \frac1 t \sum_{s=0}^{T-1} \| \nabla f\left(x^{s}\right) \|^2 
    \leq \frac{C(\varepsilon_1, \varepsilon_2, \gamma) + C_3 T}{t}.
\end{align*}

\end{proof}

\begin{remark} \label{remark-L-range}
The previous Theorem says that for any $L < 2L_{\text{train}}$, there exists a final training iteration $T$ such that convergence holds using \cref{alg--test-algo}. Conversely, if $\| G_{\theta_T} - \tau I \| \le \varepsilon_1$, and $\| H_{\phi_T} \| \le \varepsilon_2$, then the convergence results are obtained for all $f \in \mathcal{F}_L$, for $1/\tau \le L < h_1(\varepsilon_1, \varepsilon_2)$ or $L < \min(1/\tau, h_2(\varepsilon_1, \varepsilon_2))$, for 
\begin{equation*}
    h_1(\varepsilon_1, \varepsilon_2) := \frac{ (2\tau + 2\varepsilon_1 + \varepsilon_2) - \varepsilon_2 (\tau + \varepsilon_1)^2 - 2 \varepsilon_2^2(\tau + \varepsilon_1) - \varepsilon_2^3 }{ (\tau + \varepsilon_1)(\tau + \varepsilon_1 + \varepsilon_2) },
\end{equation*}
and $h_2(\varepsilon_1, \varepsilon_2)$ is given by
\begin{equation*}
\frac{\left(2\tau - 2\varepsilon_1 - \varepsilon_2\right)\left(1 - \varepsilon_2(\varepsilon_1 + \varepsilon_2 + \tau)\right) - \varepsilon_2 (\tau + \varepsilon_1)(\tau + \varepsilon_1 + \varepsilon_2)}{\varepsilon_2 (\varepsilon_1 + \varepsilon_2 - \tau) (\tau + \varepsilon_1)(\tau + \varepsilon_1 + \varepsilon_2) + \left(\tau^2 - 2\varepsilon_1\tau - \varepsilon_2\tau + \varepsilon_1^2 + \varepsilon_1\varepsilon_2\right) \left(1 - \varepsilon_2(\varepsilon_1 + \varepsilon_2 + \tau)\right)}.
\end{equation*}
Furthermore, $h_1$ and $h_2$ satisfy
\begin{equation*}
    \max_{\varepsilon_1, \varepsilon_2 \ge 0}h_1(\varepsilon_1, \varepsilon_2) = h_1(0,0) = 2L_{\text{train}}, \quad \inf_{\varepsilon_1, \varepsilon_2 \ge 0}h_1(\varepsilon_1, \varepsilon_2) = - \infty, \quad \inf_{\varepsilon_1, \varepsilon_2 \ge 0}h_2(\varepsilon_1, \varepsilon_2) = - \infty.
\end{equation*}
This means that given error bounds $\varepsilon_1, \varepsilon_2$, we can use $h_1$ and $h_2$ to determine on which functions our learned algorithm achieves convergence. \cref{app:convergence} provides the proof. In particular, the upper bound for the $L$-smoothness constant of test functions is maximised when $h_1$ is maximised, which occurs when $\varepsilon_1=\varepsilon_2=0$, in which case $h_1(0,0) = 2/\tau = 2 L_{\text{train}}$.

\cref{cor::ex-error-bounds} contains explicit bounds for $\varepsilon_1$ and $\varepsilon_2$ in terms of $\tau$ only for convergence to be guaranteed when using \cref{alg--test-algo}.
\end{remark}

\subsection{Rates under additional structure}

The previous Theorem ensures a convergence rate for the squared gradient norm throughout optimization for non-convex functions. The following corollary provides pointwise convergence and linear convergence rates in the Polyak--{\L}ojasiewicz case. 

\begin{definition}[Polyak--{\L}ojasiewicz inequality]
\label{def:pl-inequality}
Let $f:\mathcal{X}\to \mathbb R$ be Fr\'echet differentiable and bounded below, with
$f^\ast:=\inf_{x\in \mathcal{X}} f(x)$. $f$ satisfies the
Polyak--{\L}ojasiewicz inequality with parameter $\mu>0$ if
\[
    \frac12\|\nabla f(x)\|^2 \geq \mu \bigl(f(x)-f^\ast\bigr)
    \qquad \text{for all } x\in \mathcal{X} .
\]
\end{definition}

Note that a $\mu$-strongly convex function satisfies the
Polyak--{\L}ojasiewicz inequality with parameter $\mu>0$.
\begin{corollary}[Point convergence and P{\L} convergence rate]
\label{cor:point-convergence}
Let \cref{assumption::liminfbgdcts} hold, and let
$f\in\mathcal F_L$ for $L<2L_{\mathrm{train}}$ satisfy the P{\L} inequality
from \cref{def:pl-inequality} with parameter $\mu>0$. Suppose that, in the
notation of \cref{thm--reg-conv},
\[
    \|G_{\theta_T}-\tau I\|\leq \varepsilon_1,
    \quad
    \|H_{\phi_T}\|\leq \varepsilon_2,  \quad C_1>0, \quad
    C_2>0 .
\]
Define
\[
    \rho
    :=
    \min\left\{
        2\mu C_1,\,
        \frac{C_2}{\gamma}
    \right\}, \qquad \Psi_t:=f\left(x^t\right)-f^\ast+\gamma\|x^t-x^{t-1}\|^2.
\]
Then $\rho \in (0, 1]$, and the iterates produced by \cref{alg--test-algo} satisfy for $t \ge T$
\[
    f\left(x^t\right)-f^\ast
    \leq
    (1-\rho)^{t-T}\Psi_T, \quad \text{ and } x^t\to \bar x\in\arg\min f .
\] 
\end{corollary}
\begin{proof}
By \cref{lemma--theta-t-conv}, for some final training iteration $T$,
\[
    \|G_{\theta_T}-\tau I\|\le \varepsilon_1,
    \qquad
    \|H_{\phi_T}\|\le \varepsilon_2 .
\]
For $t\ge T$, let $p_t:=x^t-x^{t-1}$ and
$\Psi_t:=f\left(x^t\right)-f^\ast+\gamma\|p_t\|^2$. By \cref{thm--reg-conv} and the P{\L} inequality,
\[
\begin{aligned}
    \Psi_{t+1}-\Psi_t
    &\le -C_1\|\nabla f\left(x^t\right)\|^2-C_2\|p_t\|^2  \\
    &\le -2\mu C_1\bigl(f\left(x^t\right)-f^\ast\bigr)-C_2\|p_t\|^2 .
\end{aligned}
\]
Since $\rho\le 2\mu C_1$ and $\rho\le C_2/\gamma$,
\[
    2\mu C_1\bigl(f\left(x^t\right)-f^\ast\bigr)+C_2\|p_t\|^2
    \ge \rho\bigl(f\left(x^t\right)-f^\ast+\gamma\|p_t\|^2\bigr)
    =\rho\Psi_t .
\]
Thus $\Psi_{t+1}\le(1-\rho)\Psi_t$, and consequently
\[
    f\left(x^t\right)-f^\ast \le \Psi_t \le (1-\rho)^{t-T}\Psi_T .
\]
Note that by the definition of $C_2$, $0 \le C_2 / \gamma \le 1$ It remains to prove point convergence. Since $\gamma>0$,
\[
    \|p_t\|
    \le \sqrt{\frac{\Psi_t}{\gamma}}
    \le \sqrt{\frac{\Psi_{T}}{\gamma}}\,(1-\rho)^{{(t-T)}/2}.
\]
As $(1-\rho)^{1/2}<1$, we have $\sum_{t=T}^\infty \|p_t\|<\infty$. Hence,
for $m>n \ge T$,
\[
    \|x^m-x^n\|
    \le \sum_{t=n+1}^m \|p_t\|
    \to0
    \text{ as } n,m\to\infty,
\]
so $(x^t)$ is Cauchy. Since $X$ is complete, $x^t\to\bar x$ for some
$\bar x\in X$. Moreover, $\nabla f\left(x^t\right)\to0$ by \cref{thm--reg-conv};
therefore, by continuity, $\nabla f(\bar x)=0$. Applying the P{\L} inequality at
$\bar x$ yields
\[
    0=\frac12\|\nabla f(\bar x)\|^2
    \ge \mu\bigl(f(\bar x)-f^\ast\bigr).
\]
Since $f(\bar x)\ge f^\ast$, we have $f(\bar x)=f^\ast$, and hence
$\bar x\in\arg\min f$.
\end{proof}

%\begin{remark}[K{\L} convergence theory]
The proof of \cref{thm--reg-conv} shows that, under certain regularization assumptions, the learned method admits the Lyapunov function $\Psi_t$, which satisfies a sufficient-decrease estimate whenever $C_1,C_2>0$. This sufficient decrease appears in convergence analyses based on the Kurdyka--{\L}ojasiewicz inequality. One therefore may expect to yield finite length and convergence of the whole sequence to a critical point. We do not pursue this extension here, but see for instance \cite{attouch2013convergence, ochs2014ipiano}.
%\end{remark}

The previous results consider generalization guarantees when both $\theta_t$ and $\phi_t$ are learned at each iteration $t$. The following result considers generalization when one only learns a preconditioner on the gradient, i.e. $H_{\phi_t} = 0$ for all $t$. In particular, while the previous results can be applied to this case, this restriction enables a sublinear convergence rate in function optimality for convex functions, under certain conditions.  

\begin{theorem} [Convex convergence rate] \label{cor:only-g-learned}
    Suppose that only $G_{\theta_t}$ is learned at each iteration (i.e. $H_{\phi_t} = 0$ for all $t$). Assume that $\lim \inf_{t \to \infty} \lambda_t > 0$, $(\theta_t)_{t=0}^\infty$ is BGD, $G : \Theta \to \mathcal{L}(\mathcal{X})$ is continuous, and at iteration $T$, for some $0 \le \varepsilon < \tau$, we have $\| G_{\theta_T} - \tau I \| \le \varepsilon$. Then, for all convex $f \in \mathcal{F}_{L}$ with $L < 2L_{\text{train}}/ (1 + \varepsilon L_{\text{train}}) $ and any starting point $x^0$,  using \cref{alg--test-algo} we have $\nabla f\left(x^t\right) \to 0$ as $t \to \infty$. Furthermore, suppose that $(x^t)_{t=1}^\infty$ is a bounded sequence. Then there exists a constant $C>0$, such that
    \begin{equation}
        f(x^{t}) - f(x^*) \leq \frac C t.
    \end{equation}
    
    % In particular, if $\varepsilon = 0$, then convergence holds for all $L < 2 L_{\text{train}}$, and for all $\varepsilon \in [0, \tau)$ such that $1<\frac{2}{1 + \varepsilon L_{\text{train}}} \leq 2$, meaning that convergence is guaranteed for all $f \in \mathcal{F}_{L_{\text{train}}}$. 
\end{theorem}

\begin{proof}
By \cref{lemma:thm-calcs}, we have $f\left(x^{t+1}\right) \le f\left(x^t\right) - D \|\nabla f\left(x^t\right)\|^2$, where $D$ is given by 
\begin{equation}
    D := \tau \left( 1 - \frac{\tau L}{2} \right) - \varepsilon |\tau L - 1| -  \frac {L \varepsilon^2} 2.
\end{equation}
\textbf{Case 1} - $L \ge L_{\text{train}}$, then $D = (\tau  + \varepsilon) \left( 1 -  \frac L 2 \left( \tau + \varepsilon \right) \right)$, which is positive if and only if 
\begin{equation}
    L < \frac 2 {\tau + \varepsilon} = \frac 2 {1 + \varepsilon L_{\text{train}}} L_{\text{train}}.
\end{equation}
Note that $2 / ({1 + \tau L_{\text{train}}}) = 1 $, therefore in case 1 as $L \ge L_{\text{train}}$ we require $\varepsilon < \tau$.

\textbf{Case 2} - $L < L_{\text{train}}$, then $D = ( \tau - \varepsilon ) \left( 1 - \frac L 2 (\tau - \varepsilon) \right)$, which is positive if $\varepsilon < \tau$ and
\begin{equation}
    L < \frac 2 {\tau - \varepsilon} = \frac 2 {1 - \varepsilon L_{\text{train}}} L_{\text{train}},
\end{equation}
however, this is always the case when $L < L_{\text{train}}$ and $\varepsilon < \tau$.

For the convergence rate, define $D = \sup_{t=0, 1, \dots} \{\| x^t - x^* \|\}$, which is finite as $(x^t)$ is bounded. Due to the convexity of $f$ and the Cauchy-Schwarz inequality, we have that 
    \begin{align*}
        f(x^{t}) - f(x^*) &\leq \langle \nabla f\left(x^t\right) , x^t - x^* \rangle 
        \leq \| \nabla f\left(x^t\right)\| \| x^t - x^*\| 
         \leq D \| \nabla f\left(x^t\right)\|.
    \end{align*}
    Therefore for $t \geq T$ we have
    \begin{align*}
        f\left(x^{t+1}\right) \leq f\left(x^t\right) - c(\nu, \tau) \|\nabla f\left(x^t\right)\|^2 
         \leq f\left(x^t\right) - \frac{c(\nu, \tau)}{D^2}(f(x^{t}) - f(x^*))^2.
    \end{align*}
    Define $\Delta_t = f(x_t) - f^* \ge 0$ for $t > T$, then by standard arguments for gradient descent, see e.g. \cite{nesterov2018lectures}, shown in \cref{lemma-nesterov-lemma} in the Appendix, we have 
    \begin{align*}
        f\left(x^{t+T}\right) - f(x^*) \leq \frac{D^2/c(\nu, \tau)}{t}.
    \end{align*}
\end{proof}

All the results presented in this section give the worst-case convergence rate of the learned algorithm. In \cref{sec--experiments} we will see that the empirical performance of the learned algorithms may exceed that of NAG and L-BFGS.

\section{Linear parametrizations} \label{section--linear-params}

In this section, we consider `linear parametrizations' of $G$ and $H$, defined below.
\begin{definition} \label{def:line-param}
    % We call $G$ a linear parameterization if for all $x \in \mathcal{X}$, there exists a linear operator $B_x  \in \mathcal{L}(\Theta, \mathcal{X})$ such that $G_{\theta} x = B_x \theta$.
    We say a parametrization $P_{\theta}$ with parameters $\theta \in \Theta$ a linear parameterization if for all $x \in \mathcal{X}$, there exists a linear operator $B_x  \in \mathcal{L}(\Theta, \mathcal{X})$ such that $P_{\theta} x = B_x \theta$.
    % $G: \Theta \to \mathcal{L}(\mathcal{X})$ is a linear map. 
    % % $G: \Theta \to \mathcal{X}$ is linear. Then 
%     This means there exists $B_k^t \in \mathcal{L}(\Theta, \mathcal{X})$ such that
% \begin{equation} \label{eq--linear}
%     G_{\theta} \nabla f_k (x_k^t) = B_k^t \theta.
% \end{equation}
\end{definition}
Hereafter, we assume $G_{\theta}$ and $H_{\phi}$ are linear parameterizations. Therefore, there exist $B_k^t \in \mathcal{L}(\Theta, \mathcal{X}), C_k^t \in \mathcal{L}(\Phi, \mathcal{X})$ such that
\begin{equation} \label{eq--linear}
    G_{\theta} \nabla f_k (x_k^t) = B_k^t \theta, \quad H_{\phi} (x_k^t - x_k^{t-1}) = C_k^t \phi.
\end{equation}
% Similarly, we call $H$ a linear parametrization if there exists $C_k^t \in \mathcal{L}(\Phi, \mathcal{X})$ such that $H_{\phi} (x_k^t - x_k^{t-1}) = C_k^t \phi$. 
The motivation is that when $G$ and $H$ are linear parametrizations and the functions $f_k$ are convex, each optimization problem \cref{eq--opt-problem} is convex (as it is the composition of a convex function with a linear function \cite{Beck2014}). Therefore, learning comprises solving a sequence of convex optimization problems. In this case, there exist fast, provably convergent algorithms to find global solutions. In the nonconvex case, learning parameters requires solving an optimization problem that is no more difficult than minimizing the functions $f_k$. Due to the speed of training enabled by linear parameterizations, we are able to learn algorithms up to large $T$. In \cref{sec--experiments}, we see this enables algorithms to be learned up to iterations where a pre-selected tolerance has been satisfied. Four examples of linear parametrizations are provided in \cref{tab:params} for $G$. The analogue for $H$ simply requires replacing $\nabla f_k (x_k^t)$ with $(x_k^t - x_k^{t-1})$. These parametrizations are used for the numerical experiments in \cref{sec--experiments}. 

 % For example, when objective functions $f_k$ are convex, each optimization problem \cref{eq--opt-problem} is convex (as it is the composition of a convex function with a linear function \citep{Beck2014}). Therefore, learning comprises solving a sequence of convex optimization problems.

 % \rcomment{Also check boundary conditions on convolution for the adjoint - check my calculations}

\begin{table}[ht]

    \caption{Examples of linear parametrizations}
    
    \label{tab:params}
    \centering
    \begin{tabular}{c c c c}
        \hline
        \textbf{Label} & \textbf{Description} & \textbf{Parametrization $B_k^t\theta = $} & \textbf{Adjoint} $(B_k^t)^*w =$ \\
        \hline
        PS & Scalar step size & $\theta \nabla f_k(x_k^t), \theta \in \mathbb{R}$ & $\langle w,  \nabla f_k (x_k^t)\rangle \in \mathbb{R}$  \\
        PP & Pointwise operator & $\theta \odot \nabla f_k(x_k^t), \theta \in \mathcal{X}$  & $w \odot \nabla f_k (x_k^t) \in \mathcal{X}$ \\
        PC & Convolution & $\theta \ast \nabla f_k(x_k^t), \theta \in \mathcal{X}$ & $ w \ast \overline{\nabla f_k (x_k^t)} \in \mathcal{X}$  \\
        PF & Full linear operator & $\theta \nabla f_k(x_k^t), \theta \in \mathcal{L}(\mathcal{X})$ &  $w \otimes \nabla f_k (x_k^t) \in \mathcal{L}(\mathcal{X})$ \\
        \hline
    \end{tabular}
\end{table}

\begin{lemma} \label{lemma--linear-properties}
    All parametrizations in \cref{tab:params} satisfy the GGD property \cref{eq--GGD}, and are all continuous with respect to their parameters.
\end{lemma}
    
    \begin{proof}%[\textbf{Proof of \cref{lemma--linear-properties}.}]
        \begin{enumerate}
            \item For scalar step sizes, $G_{\theta} = \theta I$, take $\tilde{\theta} = \tau$. $H_{\phi} = \phi I$, take $\tilde{\phi} = 0$.
            \item For a pointwise parametrization, $G_{\theta}x = \theta \odot x$, take $\tilde{\theta} = \tau \mathbf{1}$. $H_{\phi} = \phi \odot I$, take $\tilde{\phi} = 0 \mathbf{1}$.
            \item For full operator parametrization, $G_{\theta} = \theta \in \mathcal{L}(\mathcal{X})$, take $\tilde{\theta} = \tau I$. $H_{\phi} = \phi $, take $\tilde{\phi} = 0$.
            \item For the convolutional parametrization, $H_{\phi}x = \phi \ast x$, take $\tilde{\phi} = 0$, and $G_{\theta} x = \theta \ast x$, take
            \begin{equation*}
                \tilde{\theta} (i,j) = \begin{cases}
                    \tau, \text{ if } i=j=0,\\
                    0, \text{ otherwise}.
                \end{cases}
            \end{equation*}
        \end{enumerate}
        $G_{\theta}, H_{\phi}$ are clearly continuous in $\theta$ for all listed parametrizations.
    \end{proof}

Then linear parameterizations satisfy the only assumption on $G$ and $H$ and therefore the convergence results hold.
% \begin{proposition}
%     If the assumptions from \cref{thm--reg-conv} are satisfied, then for linear parametrizations in \cref{tab:params}, we obtain the convergence results.
% \end{proposition} \label{cor--linear}

% \begin{proof} %[\textbf{Proof of \cref{cor--linear}.}]
%      With any parametrization in \cref{tab:params}, $G: \Theta \rightarrow \mathcal{L}(\mathcal{X})$ is continuous by \cref{lemma--linear-properties}. For \cref{thm--reg-conv} to hold, we then need $\theta_t$ is BGD and $\lim \inf_{t \to \infty} \lambda_t > 0$, which are both assumed. 
% \end{proof}

\subsection{Closed-form solutions}
\label{section--closedform}

If each function $f_k \in \mathcal{T}$ can be written as a least-squares function, then the parameters $\theta_t, \phi_t$ at iteration $t$ have a closed-form solution. 

\begin{proposition}[Closed-form solution for least-squares objectives]
\label{prop--general-closed-form}
Let
\[
    f_k(x)=\frac12\|A_kx-y_k\|^2,\qquad k=1,\dots,N,
\]
and suppose that \(G,H\) are linear parametrisations, so that $x_k^{t+1}(\theta,\phi) = x_k^t-B_k^t\theta+C_k^t\phi$. Set \(u:=(\theta,\phi)\), \(\tilde u:=(\tilde\theta,\tilde\phi)\), and define $r_k^t:=A_kx_k^t-y_k,R_t:=\operatorname{diag}(\lambda_t I_\Theta,\mu_t I_\Phi)$,
and \(S_k^t:\Theta\times\Phi\to\mathcal Y\) by
\[
    S_k^t u := A_kB_k^t\theta-A_kC_k^t\phi .
\]
Then a solution of \cref{eq--opt-problem} is given by
\[
    u_t
    =
    \left(
        R_t+\frac1N\sum_{k=1}^N (S_k^t)^*S_k^t
    \right)^\dagger
    \left(
        R_t\tilde u
        +
        \frac1N\sum_{k=1}^N (S_k^t)^*r_k^t
    \right),
\]
where \(M^\dagger\) denotes the Moore--Penrose pseudoinverse of $M$.
\end{proposition}

\begin{proof}
Since
\[
    A_kx_k^{t+1}(\theta,\phi)-y_k
    =
    r_k^t-S_k^tu,
\]
the objective is the convex quadratic
\[
    g_{t,\lambda_t,\mu_t}(u)
    =
    \frac1{2N}\sum_{k=1}^N\|r_k^t-S_k^tu\|^2
    +
    \frac12\langle R_t(u-\tilde u),u-\tilde u\rangle .
\]
The first-order optimality condition is therefore
\[
    \left(
        R_t+\frac1N\sum_{k=1}^N (S_k^t)^*S_k^t
    \right)u
    =
    R_t\tilde u
    +
    \frac1N\sum_{k=1}^N (S_k^t)^*r_k^t .
\]
\end{proof}

\begin{remark}[Scalar step size]
If \(B_k^t\theta=\theta\nabla f_k(x_k^t)\) and \(H_{\phi_t}=0\) for $\theta \in \mathbb{R}$, the previous corollary reduces to
\[
    \theta_t
    =
    \frac{
        \lambda_t \tau
        +
        \frac1N\sum_{k=1}^N \|\nabla f_k(x_k^t)\|^2
    }{
        \lambda_t
        +
        \frac1N\sum_{k=1}^N \|A_k\nabla f_k(x_k^t)\|^2
    } .
\]
For \(\lambda_t=0\) and \(N=1\), this recovers the exact line-search step size for least-squares objectives.
\end{remark}
Calculations for the closed-form solutions for parametrizations in \cref{tab:params} are given in \cref{section--linear-params}. 

\subsection{Learning parameters using optimization} \label{subsec:optimization}

For general functions $f_k$, a closed-form solution does not exist, and we instead require an optimization algorithm. With information of $\nabla_{\theta} g_{t, {\lambda_t}, \mu_t}(\theta, \phi)$, $\nabla_{\phi} g_{t, {\lambda_t}, \mu_t}(\theta, \phi)$, and $L_{\nabla_{\theta}g_{t, {\lambda_t}, \mu_t}}$, $L_{\nabla_{\phi}g_{t, {\lambda_t}, \mu_t}}$, the Lipschitz constants of $\nabla_{\theta} g_{t, {\lambda_t}, \mu_t}(\theta, \phi)$ and $\nabla_{\theta} g_{t, {\lambda_t}, \mu_t}(\theta, \phi)$, respectively, one can use first-order optimization algorithms for learning parameters $\theta_t, \phi_t$ without requiring step size tuning. Examples include gradient descent, NAG, L-BFGS, or stochastic optimization methods such as SGD, SVRG \cite{gower2020variancereduced} or Adam \cite{kingma2014adam} (especially for large $N$, due to both speed and memory considerations). 
% For example, one can start at an initial point $\theta_t^0$ at iteration $t$ and update via gradient descent
% \begin{equation} \label{eq--theta-t-w}
%     \theta_t^{w+1} = \theta_t^w - \frac1{L_{g_{t, {\lambda_t}, \mu_t}}}\nabla g_{t, {\lambda_t}}(\theta_t^w).
% \end{equation}
The following result details important calculations regarding linear parametrizations.

\begin{proposition} \label{prop:approx-general}
    For linear parametrizations $G$ and $H$, $\nabla _{\theta} g_{t, {\lambda_t}, \mu_t}(\theta, \phi)$ and $\nabla _{\phi} g_{t, {\lambda_t}, \mu_t}(\theta, \phi)$ can be calculated as 
\begin{align} \label{eq--grad-g-app}
    \nabla_{\theta} g_{t, {\lambda_t}, \mu_t}(\theta, \phi) &= {\lambda_t} (\theta - \tilde{\theta}) -\frac1N \sum_{k=1}^N  (B_k^t)^*\nabla f_k(x_k^t - G_{\theta}\nabla f_k(x_k^t) + H_{\phi}(x_k^t - x_k^{t-1})),\\
    \nabla_{\phi} g_{t, {\lambda_t}, \mu_t}(\theta, \phi) &= {\mu_t} (\phi - \tilde{\phi}) +\frac1N \sum_{k=1}^N  (C_k^t)^*\nabla f_k(x_k^t - G_{\theta}\nabla f_k(x_k^t) + H_{\phi}(x_k^t - x_k^{t-1})).
\end{align}
Furthermore, the corresponding Lipschitz constants of $\nabla_{\theta} g_{t, {\lambda_t}, \mu_t}(\theta, \phi)$ and $\nabla_{\phi} g_{t, {\lambda_t}, \mu_t}(\theta, \phi)$ are given by
\begin{align}
    L_{\nabla_{\theta} g_{t, {\lambda_t}, \mu_t}}&= {\lambda_t} + \frac1N \sum_{k=1}^N L_k \|B_k^t \|^2,\quad
    L_{\nabla_{\phi} g_{t, {\lambda_t}, \mu_t}} = {\mu_t} + \frac1N \sum_{k=1}^N L_k \|C_k^t \|^2,
\end{align}
respectively.
\end{proposition}

\begin{proof}
From the definition of $g_{t, \lambda_t, \mu_t}$ in \cref{eq--opt-problem}, applying the chain rule achieves the desired result. To calculate the Lipschitz constants, 
\begin{align*}
    &\| \nabla_ \theta g_{t, \lambda_t, \mu_t}(\theta_1, \phi) - \nabla _ \theta g_{t, \lambda_t, \mu_t}(\theta_2, \phi) \| 
    %&= \left\| {\lambda_t} (\theta_1 - \theta_2) +  \frac1N \sum_{k=1}^N (A_kB_k^t)^* (\nabla f_k(x_k^t - B_k^t \theta_2) - \nabla f_k(x_k^t - B_k^t \theta_2)) \right\|\\
    %&\leq {\lambda_t} \| \theta_1 - \theta_2 \| +  \frac1N \sum_{k=1}^N  \left\|(A_kB_k^t)^* (\nabla f_k(x_k^t - B_k^t \theta_2) - \nabla f_k(x_k^t - B_k^t \theta_1)) \right\|\\
    \\
    &\leq  {\lambda_t} \| \theta_1 - \theta_2 \| + \frac1N \sum_{k=1}^N \| B_k^t \| \left\| \nabla f_k(x_k^t - B_k^t \theta_2 + C_k^t \phi) - \nabla f_k(x_k^t - B_k^t \theta_1 + C_k^t \phi) \right\|\\
    &\leq  {\lambda_t} \| \theta_1 - \theta_2 \| +  \frac1N \sum_{k=1}^N L_k \| B_k^t \| \left\| B_k^t( \theta_1  - \theta_2)\right\| \leq  ({\lambda_t}  + \frac1N \sum_{k=1}^N L_k \| B_k^t \|^2 )\left\| \theta_1  - \theta_2 \right\|,
\end{align*}
and, similarly for $\phi$,
\begin{align*}
    &\| \nabla_ \phi g_{t, \lambda_t, \mu_t}(\theta, \phi_1) - \nabla _ \phi g_{t, \lambda_t, \mu_t}(\theta, \phi_2) \| \leq  \left({\mu_t}  + \frac1N \sum_{k=1}^N L_k \| C_k^t \|^2 \right)\left\| \phi_1  - \phi_2 \right\|.
\end{align*}
% Due to the properties of the triangle inequality, the Cauchy-Schwarz inequality, and the operator norm, this bound is tight. Therefore the Lipschitz constant of $\nabla_{\theta} g_{t, {\lambda_t}, \mu_t}(\theta, \phi)$ is given by
% \begin{equation*}
%     {\lambda_t} + \frac1N \sum_{k=1}^N L_k \| B_k^t \|^2.
% \end{equation*}
\end{proof}
We provide parametrization-specific calculations in \cref{tab:params} in \cref{appendix:approx-linear-params}.

The BGD property is easily verified during training for each parametrization by checking $g_{t, \lambda_t, \mu_t}(\theta_t, \phi_t) \leq g_{t, \lambda_t, \mu_t}(\tilde{\theta}, \tilde{\phi})$, and can be ensured by initializing $(\theta_t, \phi_t) = (\tilde{\theta}, \tilde{\phi})$ and any descent in $g_{t, \lambda_t}$ would ensure the BGD property.

\section{Numerical Experiments} \label{sec--experiments}

\subsection{Setup}
%\textbf{The optimization problem.} 
In this section, we test the linear parametrizations in \cref{tab:params} on two inverse problems in imaging: image deblurring and CT. We create observations from given ground-truth data according to the model $y = Ax_{\text{true}} + \varepsilon$ using the specified forward operator $A$ and noise distribution. Once these observations have been created, the ground-truth data are no longer used. For both inverse problems, $\mathcal{X} = \mathbb{R}^{h_1 \times h_2}, \mathcal{Y} = \mathbb{R}^{h_3 \times h_4}$ for $h_1, h_2, h_3, h_4 \in \mathbb{N}$, and $\varepsilon$ is noise drawn from iid zero-mean Gaussian distributions. To approximate $x_{\text{true}}$ from $y$, we solve
\begin{equation} \label{eq--objective-function}
    \min_x \left\{ f(x) := \frac12 \|Ax - y\|^2 + \alpha \mathcal{S}(x)\right\},
\end{equation}
for a regularizer $\mathcal{S}: \mathcal{X} \to \mathbb{R}$ and a fixed regularization parameter $\alpha$.

%We implement L-BFGS using Wolfe line search and NAG with backtracking line search \comment{fix if I change this}. 

\paragraph{Problem 1: Huber TV-regularized Image Deblurring}
In \cref{experiment::blur} we consider an Image deblurring problem. The forward operator $A$ in \cref{eq--objective-function} is a Gaussian blur with a $5 \times 5$ kernel size and a standard deviation $\sigma = 1.5$, and we normalize the forward operator so that $\|A\| = 1$. We use the STL-10 dataset \cite{STL10} with greyscale images of size $96 \times 96$ as $\mathcal{X}$. The noise $\varepsilon$ is modeled with a standard deviation of $0.01$, and we set $\alpha = 2 \times 10^{-4}$ and the Huber parameter $\epsilon = 0.005$, resulting in $L = 1.32$. The initial point $x^0 = y \in \mathcal{Y} = \mathcal{X}$ is chosen equal to the observation. 

For the regularizer, we use the Huber Total Variation \cite{rudin1992TV, huberloss1992}, defined as
    \begin{equation} \label{eq:huber-tv}
        \mathcal{S}(x)= \sum_{i,j=1}^{h_1, h_2} h_{\epsilon}\left(\sqrt{(\nabla x)_{i, j, 1}^2+(\nabla x)_{i, j, 2}^2}\right), \quad
        h_{\varepsilon}(s) = \begin{cases}
            \frac{1}{2\epsilon}{s^2}, \text{ if } |s| \leq \epsilon \\
            |s| - \frac{\epsilon}{2}, \text{ otherwise},
        \end{cases}
    \end{equation}
    where $\epsilon > 0$ is a hyperparameter and the finite difference operator $\nabla: \mathbb{R}^{h_1 \times h_2} \rightarrow \mathbb{R}^{h_1 \times h_2 \times 2}$ is defined in \cite{chambolle2016introduction}. Note that this choice of regularizer makes the function $f$ non-quadratic. Then, each function $f$ is $L$-smooth, where $L = 1 + 8\alpha/\epsilon$ \cite{chambolle2016introduction}. 

\paragraph{Problem 2: Huber TV-regularized Computed Tomography}
In \cref{sec::ct}, we consider Computed Tomography. The forward operator $A$ in \cref{eq--objective-function} is the discrete Radon transform in 2D, and we simulate CT measurements using ODL with ASTRA toolbox as the backend \cite{ODL, ASTRA} with a parallel-beam geometry and $180$ projection angles evenly distributed over a $180$-degree range. We normalize the forward operator so that $\|A\| = 1$. For the dataset, we use ground-truth images in the SARS-CoV-2 CT-scan dataset \cite{SARS-CoV-2CT-scan}, and in optimization take the initial point $x^0 = 0 \in \mathcal{X}$. We take $\mathcal{X} = \mathbb{R}^{256 \times 256}$ and $\mathcal{Y} = \mathbb{R}^{180 \times 363}$, where $363 = \left\lceil{256 \sqrt{2}}\right\rceil$. The noise $\varepsilon$ is modeled with standard deviation $0.01$. The regularizer is chosen as \cref{eq:huber-tv}, and we fix the regularization parameter $\alpha = 3 \times 10^{-4}$ and the Huber parameter $\epsilon = 0.01$, resulting in $L = 1.24$. 
    
\paragraph{Problem 3: Image Deblurring with a Nonconvex Regularizer} In \cref{sec::nonconvex-exp}, we consider the image deblurring problem with he forward operator $A$, the image dataset, the choice of $x^0$, and the standard deviation of the noise $\varepsilon$ are chosen as in \cref{experiment::blur}. We use the Weakly Convex Ridge Regularizer, with learned parameters as in \cite{WCRR}. This regularizer can be written as
    \begin{equation} \label{eq::wcrr}
        \mathcal{S}(x) =  \sum_{i=1}^{N_C} \sum_{k \in \mathbb{R}^2} \psi_i\left(\left(h_i * x\right)[k]\right),
    \end{equation}
    where $\left(h_i\right)_{i=1}^{N_C}$ are learned filters, and $\left(\psi_i\right)_{i=1}^{N_C}$ are potential functions with Lipschitz continuous derivative. The regularizer depends on a hyperparameter $\sigma$, and we take $\sigma = 18$. Then the regularizer satisfies $\operatorname{Lip}(\nabla \mathcal{S}) \leq 48.032$, and we take $\alpha = 3 \times 10^{-3}$, so that each function $f$ is $L$-smooth, where $L \leq 1.1441$. The parameters $\sigma, \alpha$ were selected using a grid search for optimizing reconstruction PSNR to the ground truth images in the training set.

\paragraph{Learning parameters} For a parametrization in \cref{tab:params}, to learn parameters $(\theta_t, \phi_t)$, we initialize $(\theta_t^0, \phi_t^0) = (\tilde{\theta}, \tilde{\phi})$ to ensure that any descent in the loss function during parameter learning ensures that the learned parameters are BGD. We apply NAG for solving the optimization problem \cref{eq--opt-problem}.
% , with early stopping using a validation set of $25$ functions, in order to improve generalization performance of the learned algorithm. 
When only $\theta_t$ is learned, we stop NAG after $500$ iterations for the scalar parametrization and $5000$ iterations otherwise. When both $\theta_t$ and $\phi_t$ are learned, for the scalar parametrization, we take $100$ iterations of NAG in $\theta_t$, then $100$ iterations of NAG in $\phi_t$ (as associated step-sizes can vary by a factor of $100$), and repeat $10$ times. For other parametrizations, we use the same procedure but use $1000$ iterations instead of $100$. We use a training set of $25$ functions. Unless otherwise stated, the learned convolutional kernels have dimensions $h_1 \times h_2$, matching the size of the images in $\mathcal{X}$. 
%See \cref{fig:blur-kernels} for a comparison with different kernel sizes, for the imaging deblurring problem seen below in section. 

\paragraph{Training details}
The following outlines final training iterations $T$ and training durations for the parametrizations in \cref{tab:params} with $\lambda_t = \mu_t = 0$ for all iterations. $T$ is not a user-chosen hyperparameter, and training is stopped once the learned method reaches a satisfactory performance for the problem at hand. For Problems 1 and 3, $\| \nabla f\left(x^t\right)\|_2^2 < \| \nabla f\left(x^0\right)\|_2^2 / 10^7$, and for Problem 2, $\| \nabla f\left(x^t\right)\|_2^2 < \| \nabla f\left(x^0\right)\|_2^2 / 10^{12}$. For learned algorithms with momentum, i.e. learning parameters $\phi_t$ in addition to $\theta_t$, we prefix the parametrization labels in \cref{tab:params} with ``M-''.  

% Problem 1: For the M-PS parametrization, parameters were learned up to iteration $T=100$, and training took $52$ minutes. M-PC, $T=80$, $16.5$ hours. M-PP, $T=100$, $130$ minutes. PS, $T=400$, $52$ minutes. PC, $T=400$, $20.5$ hours. PP, $T=400$, $130$ minutes. 

% Problem 2: M-PC $T=50$, $33$ hours. M-PS, $T=100$, $5$ hours. 

% Problem 3: M-PC, $T=100$, $31$ hours. M-PS, $T=150$, $4.5$ hours.
For Problem 1, the M-PS parametrization parameters were learned up to iteration $T=100$, which took $52$ minutes. The M-PC parametrization was trained for $T=80$ iterations and required $16.5$ hours, while the M-PP parametrization was trained for $T=100$ iterations and took $130$ minutes. The PS, PC, and PP parametrizations were trained for $T=400$ iterations. This took $52$ minutes for PS, $20.5$ hours for PC, and $130$ minutes for PP. For Problem 2, the M-PC parametrization with $T=50$ required $33$ hours for training, whereas the M-PS parametrization with $T=100$ took $5$ hours. For Problem 3, the M-PC parametrization was trained for $T=100$ iterations, taking $31$ hours, while the M-PS parametrization was trained for $T=150$ iterations, which took $4.5$ hours.

In our experiments, the training and test problems have the same fixed dimension. Nevertheless, the PS and PC parametrizations naturally extend across resolutions, although no guarantees about performance can be made. Whereas the PP and PF parametrizations depend on the discretized dimension and would typically require retraining if the resolution changed.

\paragraph{Evaluation} Given a dataset of functions ${f_1, \dots, f_N}$, the mean value at iteration $t$ is defined as $F\left(x^t\right) = 1/N \sum_{k=1}^N f_k(x_k^t)$. Furthermore, we define ``function optimality'' for a function $f$ with minimizer $x_f^*$ at iteration $t$ by $(f\left(x^t\right) - f(x_f^*))/(f\left(x^0\right) - f(x_f^*))$. For a function $f$, its approximate minimizer $x_{f}^* \in \mathcal{X}$ is calculated using NAG for $2000$ iterations. For a test set of $100$ functions, we visualize the maximum and minimum function optimality over the dataset and the function optimality for $F$.

\paragraph{Benchmark Algorithms} The learned algorithms are compared to NAG with backtracking \cite{beck2009fast} and L-BFGS with the Wolfe conditions \cite{wolfe1969convergence}. For Problem 1, we compare to a handcrafted preconditioner of the form $(\delta I + A^*A)^{-1}$ due to ill-conditioning, and use backtracking at each step so that the update algorithm reads $x^{t+1} = x^t - \alpha_t P \nabla f\left(x^t\right)$, which we denote by PGD. We take $\delta \approx 0.032$, which is found using the Nelder--Mead algorithm to minimize the objective value after $100$ iterations of the PGD algorithm. For Problem 1 and Problem 2, we compare the performance of our proposed algorithm to the deviation-based approach \cite{banert2024accelerated}, with update rule given by $x^{t+1} = x^t - \tfrac1L (\nabla f\left(x^t\right) + \Delta_t)$, where $\Delta_t$ is equal to $\Delta_t := \alpha h_t \| \nabla f\left(x^t\right) \|/\sqrt{1+\|h_t\|^2}$ for $\alpha = 0.9$ and $h_t = h_t(x^t, \nabla f\left(x^t\right), \Delta_{t-1})$ as the output of the Neural Network with the same architecture as in \cite{banert2024accelerated}. We use the same training set as for learning our proposed algorithm and use Adam optimizer \cite{kingma2014adam} with learning rate $10^{-3}$ and a batch size of $1$ due to memory constraints. Computations were performed in single precision on an Nvidia RTX 3600 12GB GPU.

\subsection{Results for Problem 1} \label{experiment::blur}

% \textbf{Training details.} Training with the greedy method was performed with $\lambda_t, \mu_t = 0$ for all $t$ for the parameterizations M-PS, M-PP, M-PC, PS, PP, PC. The following lists details of training times and final training iterations $T$ for the different parametrizations. M-PS, $T=100$, $52$ minutes. M-PC, $T=80$, $16.5$ hours. M-PP, $T=100$, $130$ minutes. PS, $T=400$, $52$ minutes. PC, $T=400$, $20.5$ hours. PP, $T=400$, $130$ minutes.

\begin{figure}[h!]
    \centering
    \begin{subfigure}[t]{0.33\textwidth}
        \centering
        \includegraphics[width=\textwidth]{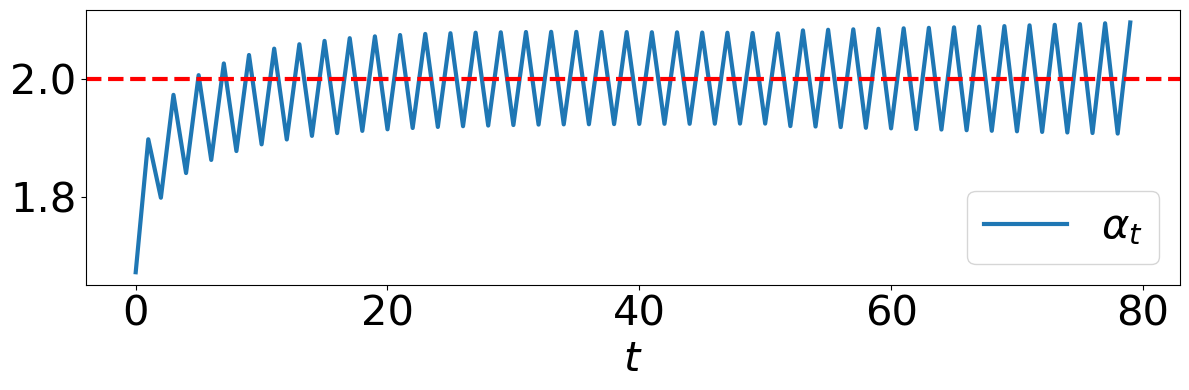} % Replace with your image file
        
        \caption{Learned scalar parameters PS.}
        \label{fig:blur-scalars}
    \end{subfigure}
    \begin{subfigure}[t]{.64\textwidth}
        \centering
        \includegraphics[width=\textwidth]{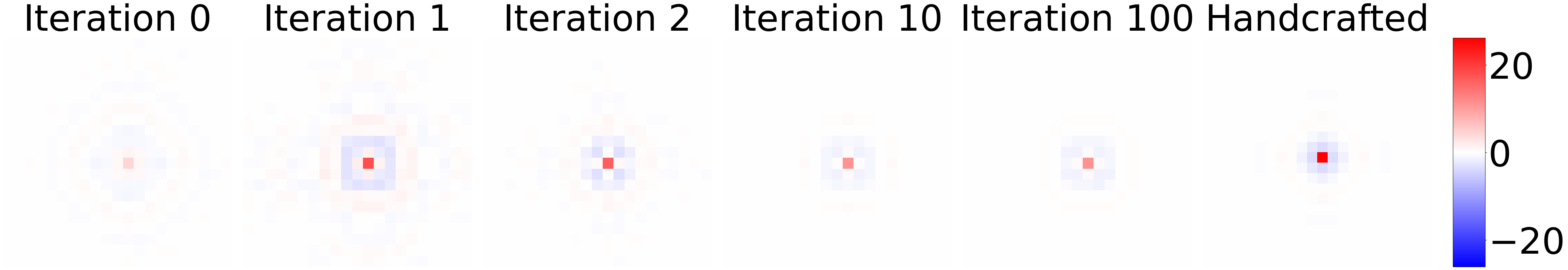}
        
        \caption{$21 \times 21$ center crop of PC kernels, and the handcrafted kernel.}
        \label{fig:blur-kernels}
    \end{subfigure}
    \newline
    % \begin{subfigure}[t]{.45\textwidth}
    %     \centering
    %     \includegraphics[width=\textwidth]{images/blur/learned_pw.png}
        
    %     \caption{}
    %     \label{fig:blur-diags}
    % \end{subfigure}
    % \begin{subfigure}[t]{0.33\textwidth}
    %     \centering
    %     \includegraphics[width=\textwidth]{images/blur/learned_alpha_betas.png} % Replace with your image file
        
    %     \caption{}
    %     \label{fig:blur-scalars-mom}
    % \end{subfigure} 

    \begin{subfigure}[t]{0.4\textwidth}
        \centering
        \includegraphics[width=\textwidth]{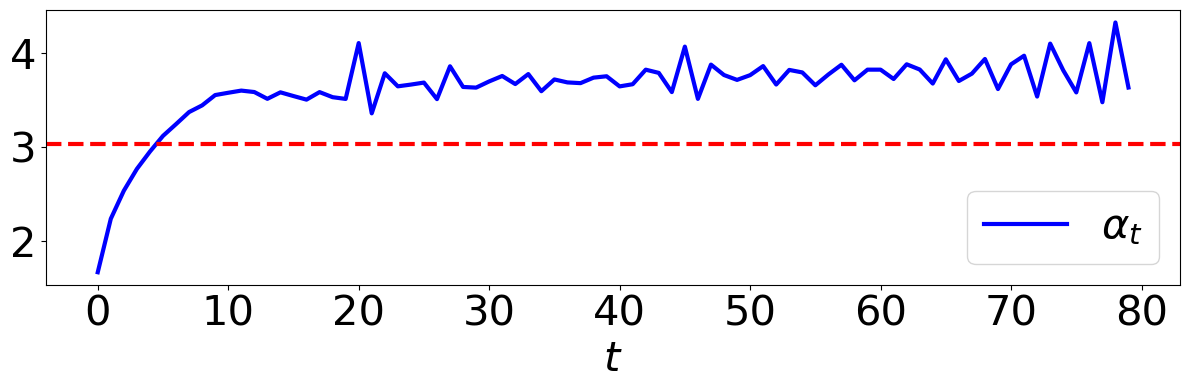} % Replace with your image file
        
        \caption{Learned $\alpha_t$ M-PS scalars and the line $y=4/L$.}
        \label{fig:blur-alpha-p}
    \end{subfigure} 
    \begin{subfigure}[t]{0.4\textwidth}
        \centering
        \includegraphics[width=\textwidth]{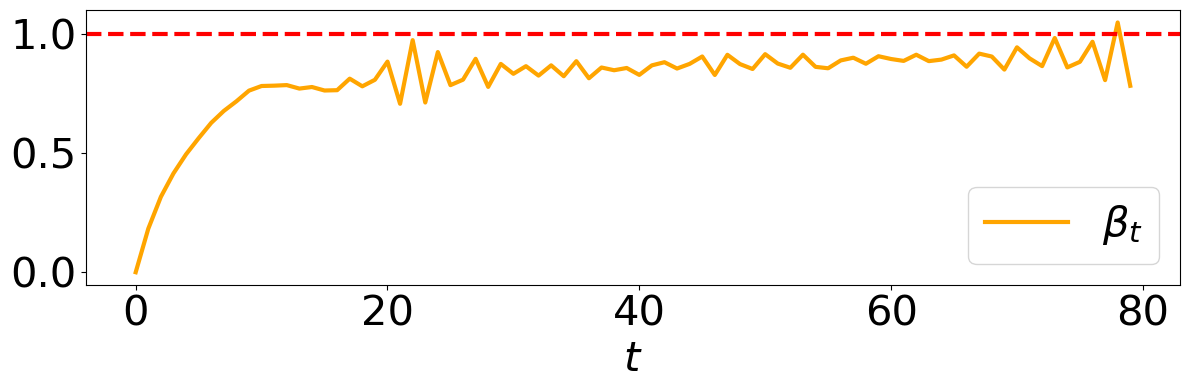} % Replace with your image file
        
        \caption{Learned $\beta_t$ M-PS scalars and the line $y=1$.}
        \label{fig:blur-alpha-q}
    \end{subfigure} 
    \newline
    \begin{subfigure}[t]{.49\textwidth}
        \centering
        \includegraphics[width=\textwidth]{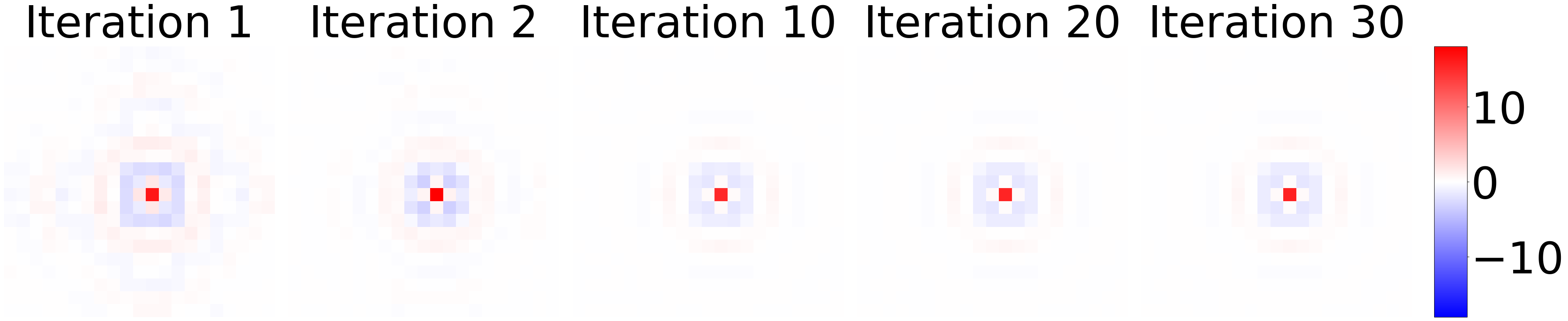}
        
        \caption{$21 \times 21$ center crop of $\theta_t$ M-PC kernels.}
        \label{fig:blur-p-kernels}
    \end{subfigure}
    \begin{subfigure}[t]{.49\textwidth}
        \centering
        \includegraphics[width=\textwidth]{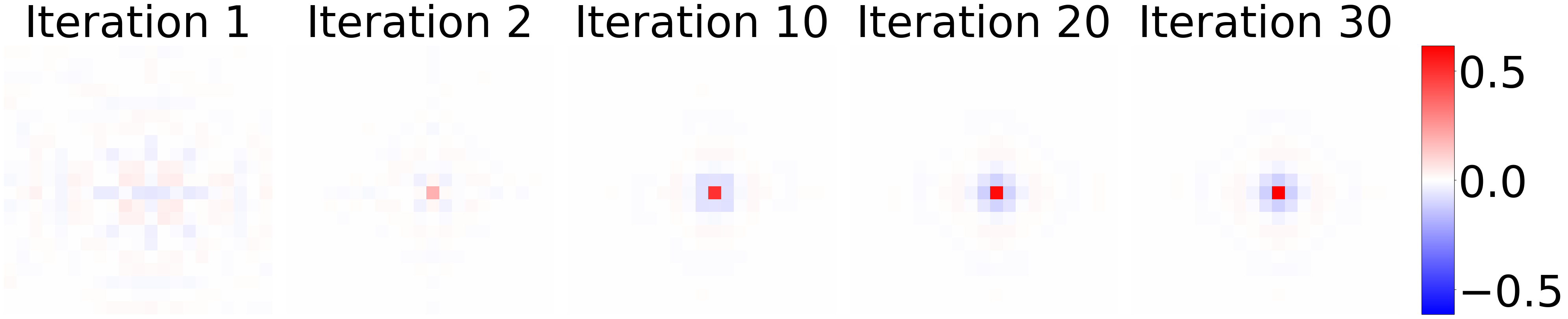}
        
        \caption{$21 \times 21$ center crop of $\phi_t$ M-PC kernels.}
        \label{fig:blur-q-kernels}
    \end{subfigure}
    \caption{Learned preconditioners for Problem 1.}
    \label{fig:blur-diags-scalars-kernels}
\end{figure}

\paragraph{Visualizing learned preconditioners} \cref{fig:blur-scalars} shows that the learned scalar parameters PS eventually fluctuate around $2/L$, which is outside of the range of provable convergence of gradient descent with a constant step size. Despite this, the learned algorithm leads to convergence on training data as $t \to \infty$ by \cref{prop:conv}. Step size sequences which contain steps larger than $2/L$ and which guarantee convergence with an accelerated convergence rate are considered in \cite{grimmer2024longsteps, silversteps}. Learned step sizes above $2/L$ are also encountered in \cite{learningalgohyperparams}. The learned convolutional kernels PC in \cref{fig:blur-kernels} contain positive and negative values and are predominantly weighted towards the center, suggesting that information from neighboring pixels is prioritized over more distant ones.  As the number of iterations increases, the kernels exhibit increasing similarity, though no formal convergence result has been established. The handcrafted preconditioner $P$ corresponds to convolution with the kernel shown in \cref{fig:blur-kernels}. \cref{fig:blur-alpha-p} and \cref{fig:blur-alpha-q} show the learned scalars and compare their values with the asymptotic optimal Heavy Ball values for $\mu$-strongly convex and $L$-smooth functions \cite{polyak1987}, for $\mu \approx 0$, $\alpha^* =  4/ {(\sqrt{L} + \sqrt{\mu})^2}\approx 4/L$ and $\beta^* = ((\sqrt{L} - \sqrt{\mu})/(\sqrt{L} + \sqrt{\mu}))^2 \approx 1$. The learned $\alpha_t$ approach $4$ instead of $4/L$, possibly as the data-fit term is $1$-smooth. \cref{fig:blur-p-kernels} and \cref{fig:blur-q-kernels} consider the M-PC parametrization, where the kernels are weighted heavily towards the center, as seen in the learned kernels in the PC parametrization \cref{fig:blur-kernels}, and with $\phi_t$ having pixel values of magnitude less than $1$. 

\begin{figure}[h!]
    \centering
    \begin{subfigure}[t]{0.38\textwidth}
        \centering
        \includegraphics[width=\textwidth]{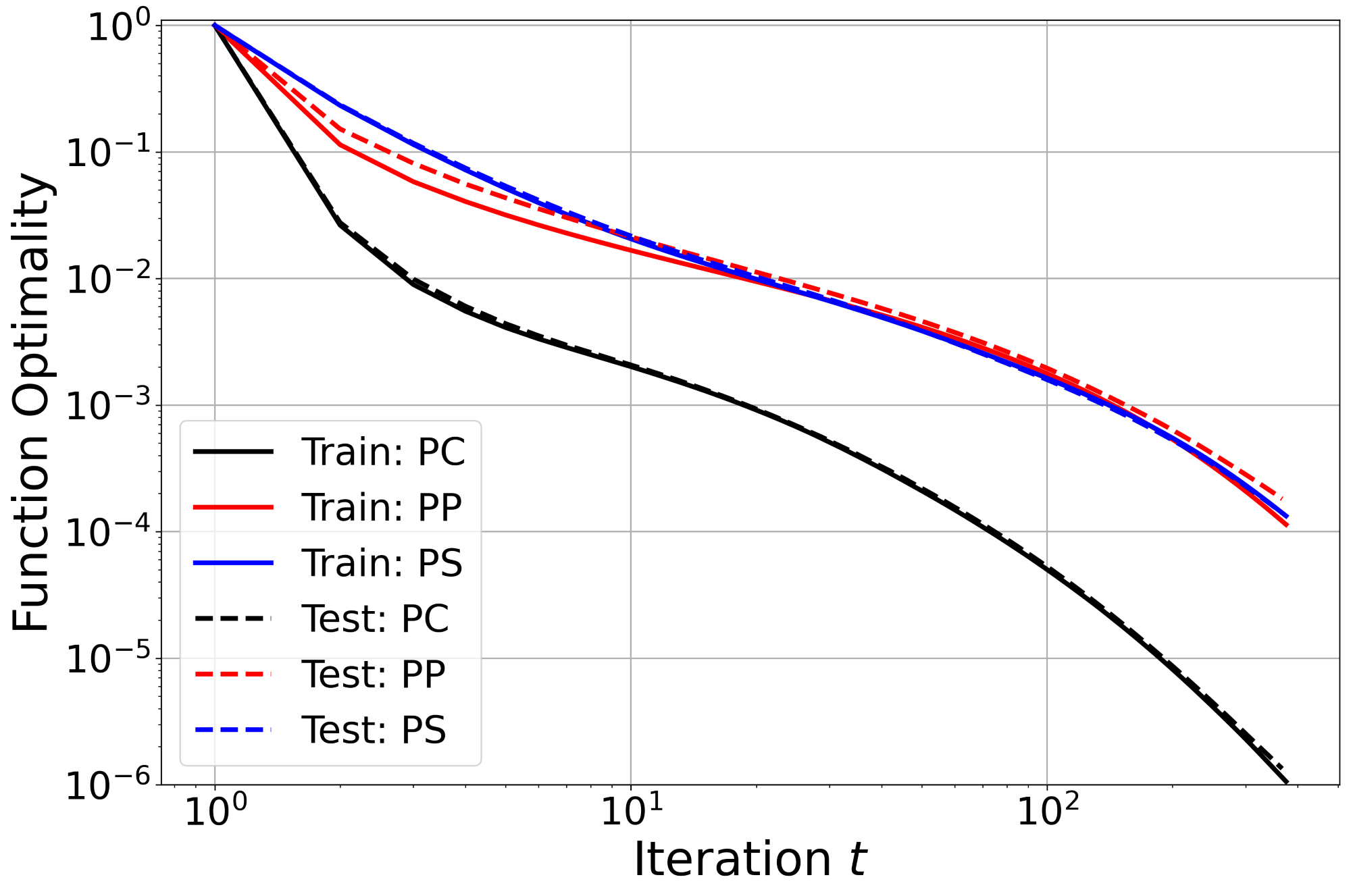} % Replace with your image file
        %\caption{}
        \label{fig::gen-orig}
    \end{subfigure}
    \hspace{15px}
    \begin{subfigure}[t]{0.38\textwidth}
        \centering
        \includegraphics[width=\textwidth]{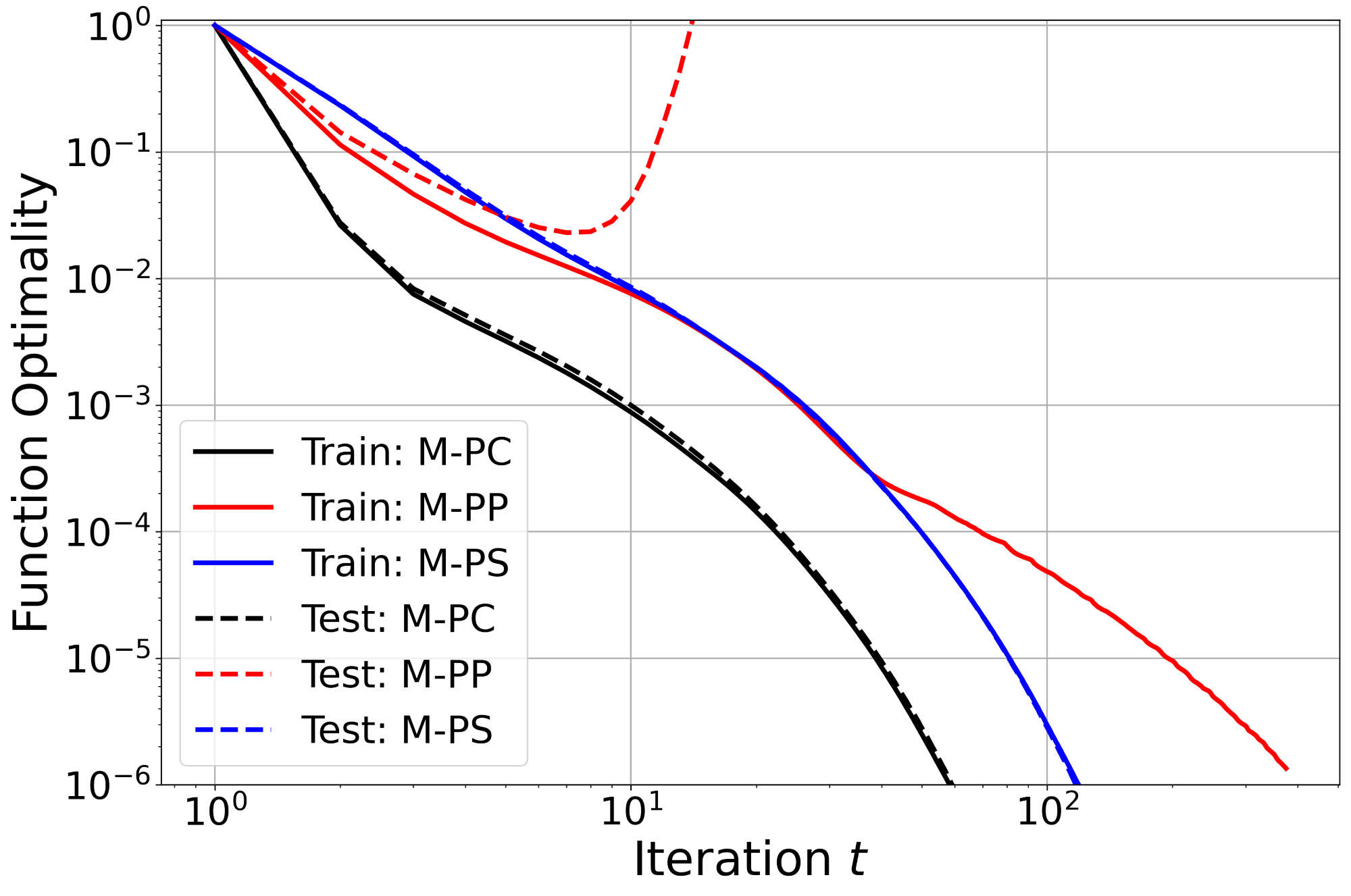} % Replace with your image file
        %\caption{}
        \label{fig::gen-mom}
    \end{subfigure}
    
    \caption{Generalization performance of learned preconditioners with and without momentum for Problem 1. Left: Learned methods without momentum generalize well to the test data. Right: Methods with momentum. M-PP does not generalize well, but M-PS and M-PC do.}
    \label{fig--blur-gen-test}
\end{figure}

\begin{figure}[h!]
    \centering
    \begin{subfigure}[t]{0.39\textwidth}
        \centering
        \includegraphics[width=\textwidth]{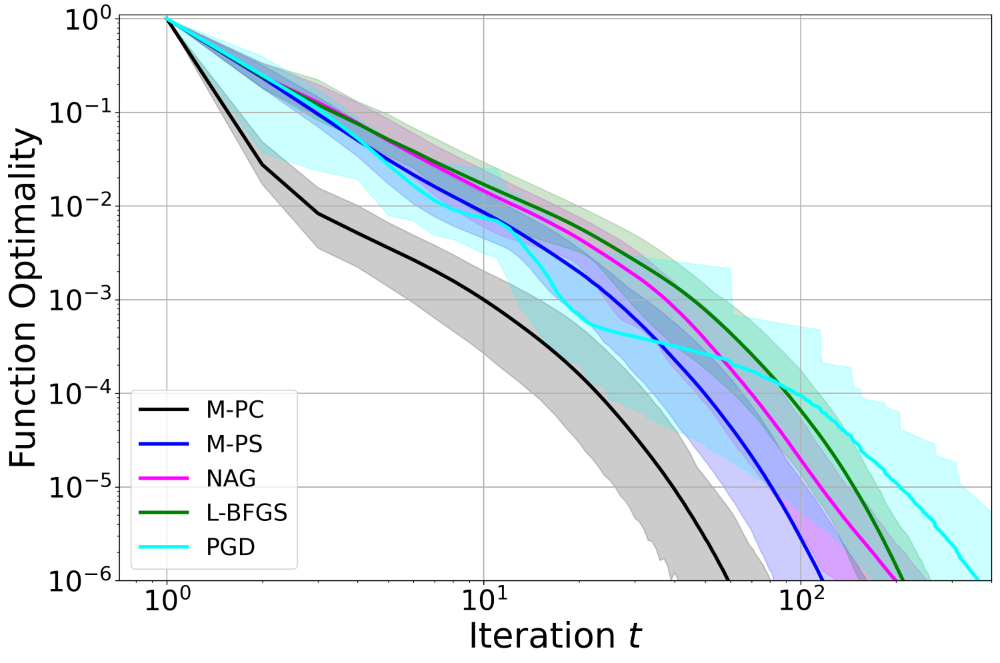} % Replace with your image file
        %\caption{}
        \label{blur-test-iter}
    \end{subfigure}
    \begin{subfigure}[t]{0.39\textwidth}
        \centering
        \includegraphics[width=\textwidth]{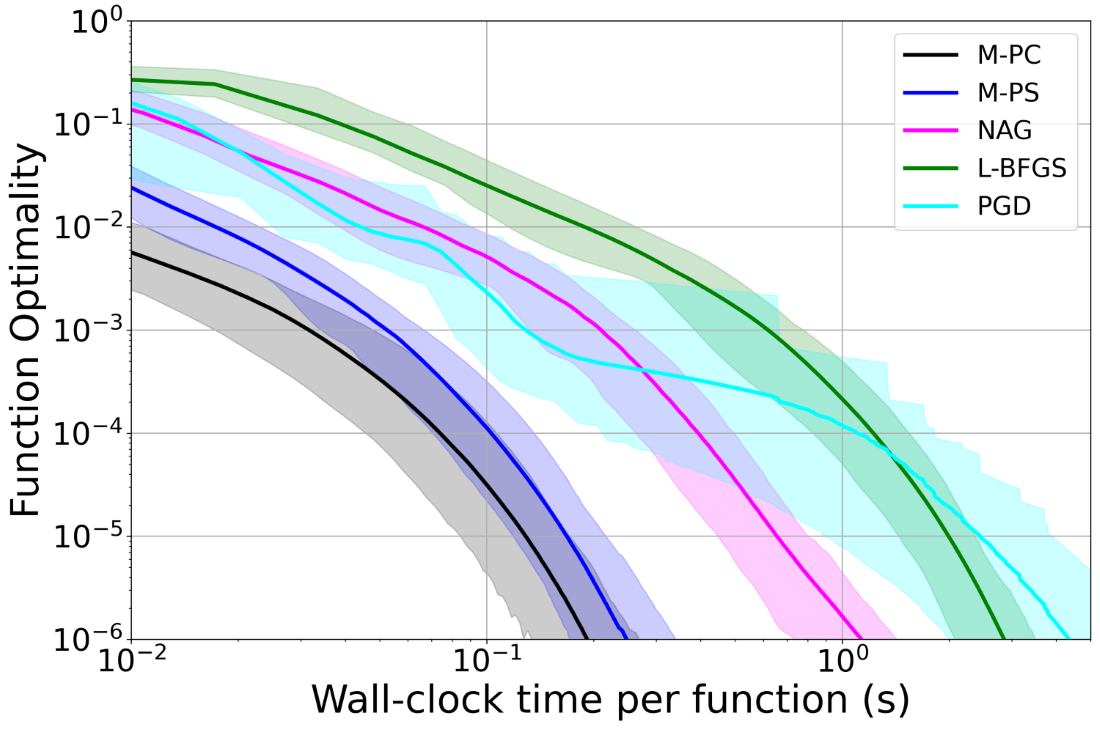} % Replace with your image file
        %\caption{}
        \label{blur-test-time}
    \end{subfigure}
    \hspace{10px}
    % \begin{subfigure}[t]{0.49\textwidth}
    %     \centering
    %     \includegraphics[width=\textwidth]{blur_images/out-of-iteration/OOI-scalar-pointwise.png} % Replace with your image file
    %     \caption{}
    %     \label{fig--2d}
    % \end{subfigure}
    
    \caption{Test performance of the proposed method with M-PC and M-PS parametrizations versus benchmark non-learned algorithms for Problem 1. We see that the learned algorithms significantly outperform the benchmark algorithms, both in terms of iterations (left) and wall-clock time (right). Intervals around each mean represent maximum and minimum values over the dataset. }
    \label{fig--blur-test}
\end{figure}

\paragraph{Generalization of learned preconditioners} \cref{fig--blur-gen-test} shows that the learned param\-etrizations PS, PP, and PC generalize well to test data. PP performs comparably to PS for this example, despite having an equal number of parameters as PC, which captures global information of the image, rather than only pixel-level details. Furthermore, we see that the performance in the training and test sets for the M-PC and M-PS parameterizations are similar. However, for the M-PP parametrization, the test curve diverges, meaning that the learned parameters overfit to the training data. When comparing both plots, we also see the improved performance that comes as a result of learning of parameters $\phi_t$. For example, the M-PC parametrization reaches a tolerance of $10^{-6}$ in under $60$ iterations, whereas the PC parametrization reaches this tolerance around iteration $400$. Furthermore, the PS parametrization does not reach a tolerance of $10^{-4}$ after $400$ iterations, but the M-PS parametrization reaches $10^{-6}$ after just over $100$ iterations. Therefore, we only compare the learned M-PC and M-PS preconditioners against benchmark algorithms in \cref{fig--blur-test}.

\paragraph{Learned algorithm performance}  \cref{fig--blur-test} shows that the learned M-PC parametrization outperforms NAG and L-BFGS on the test data. We also see that the performance with respect to wall-clock time of the learned algorithms greatly outperforms NAG and L-BFGS. Due to the increased cost of convolution, the learned M-PS algorithm performs more similarly to M-PC, but the learned kernels still outperform the other algorithms. \cref{fig--blur-test} shows that the learned algorithms with M-PC and M-PS parametrizations significantly outperform the handcrafted PGD algorithm.

% \textbf{Comparison to a hand-crafted convolutional preconditioner.}  \cref{fig--blur-test} shows that the learned algorithms with M-PC and M-PS parametrizations significantly outperform PGD.

\begin{figure}[h!]
    \centering
    \begin{subfigure}[t]{0.24\textwidth}
        \centering
        \includegraphics[width=.8\textwidth]{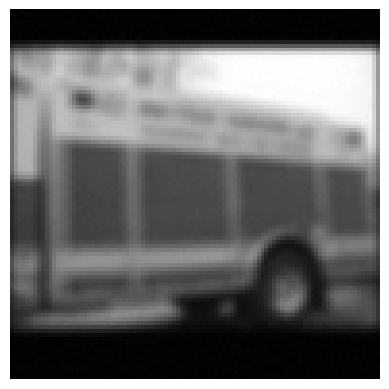}
        
        \caption{Observation}
    \end{subfigure}
    \begin{subfigure}[t]{0.24\textwidth}
        \centering
        \includegraphics[width=.8\textwidth]{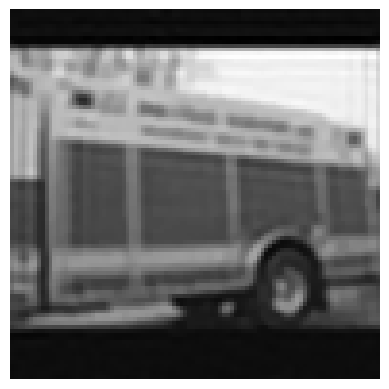}
        
        \caption{NAG iteration $15$}
    \end{subfigure}
    \begin{subfigure}[t]{0.24\textwidth}
        \centering
        \includegraphics[width=.8\textwidth]{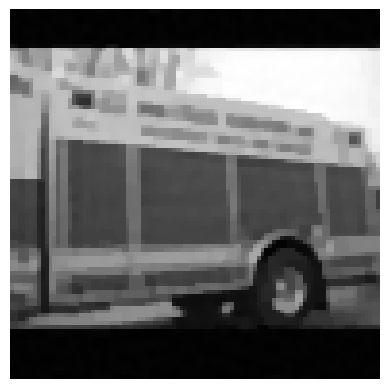}
        
        \caption{M-PC iteration $15$}
    \end{subfigure}
    \begin{subfigure}[t]{0.24\textwidth}
        \centering
        \includegraphics[width=.8\textwidth]{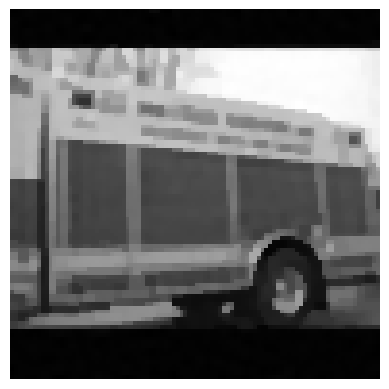}
        
        \caption{Final reconstruction.}
    \end{subfigure}
    
    \caption{A Comparison of reconstructions for Problem 1. The learned M-PC algorithm provides significantly better reconstruction at iteration $15$.}
    \label{fig:blur-comparisons}
\end{figure}

\paragraph{Reconstruction comparison} \cref{fig:blur-comparisons} demonstrates qualitatively that the learned convolutional algorithm with momentum achieves a high-quality image reconstruction in only $15$ iterations, at which point NAG produces a significantly lower-quality reconstruction, indicating the effect learning to optimize can have to speed up reconstruction.

% For the operator $A$ given by a Gaussian blur with standard deviation $1.5$ and kernel size $5 \times 5$, and a constant $\delta = 0.2$, the operator $(\delta I + A^*A)^{-1}$ corresponds to a convolution with the kernel gives as in \cref{fig:blur-handcrafted-kernel}.%\cref{app-fig--inverse-kernel}.
% \begin{figure}[h!]
%     \centering
%     \includegraphics[width=0.3\textwidth]{blur_images/inv_kernel.png} % Replace with your image file
    
%     \caption{The kernel corresponding to the operator $(\delta I + A^*A)^{-1}$. \rcomment{How is this similar to what you have learned?}}
%     \label{app-fig--inverse-kernel}
% \end{figure}

\begin{figure}
    \centering
    \centering
    \begin{subfigure}[t]{0.38\textwidth}
        \centering
        \includegraphics[width=\textwidth]{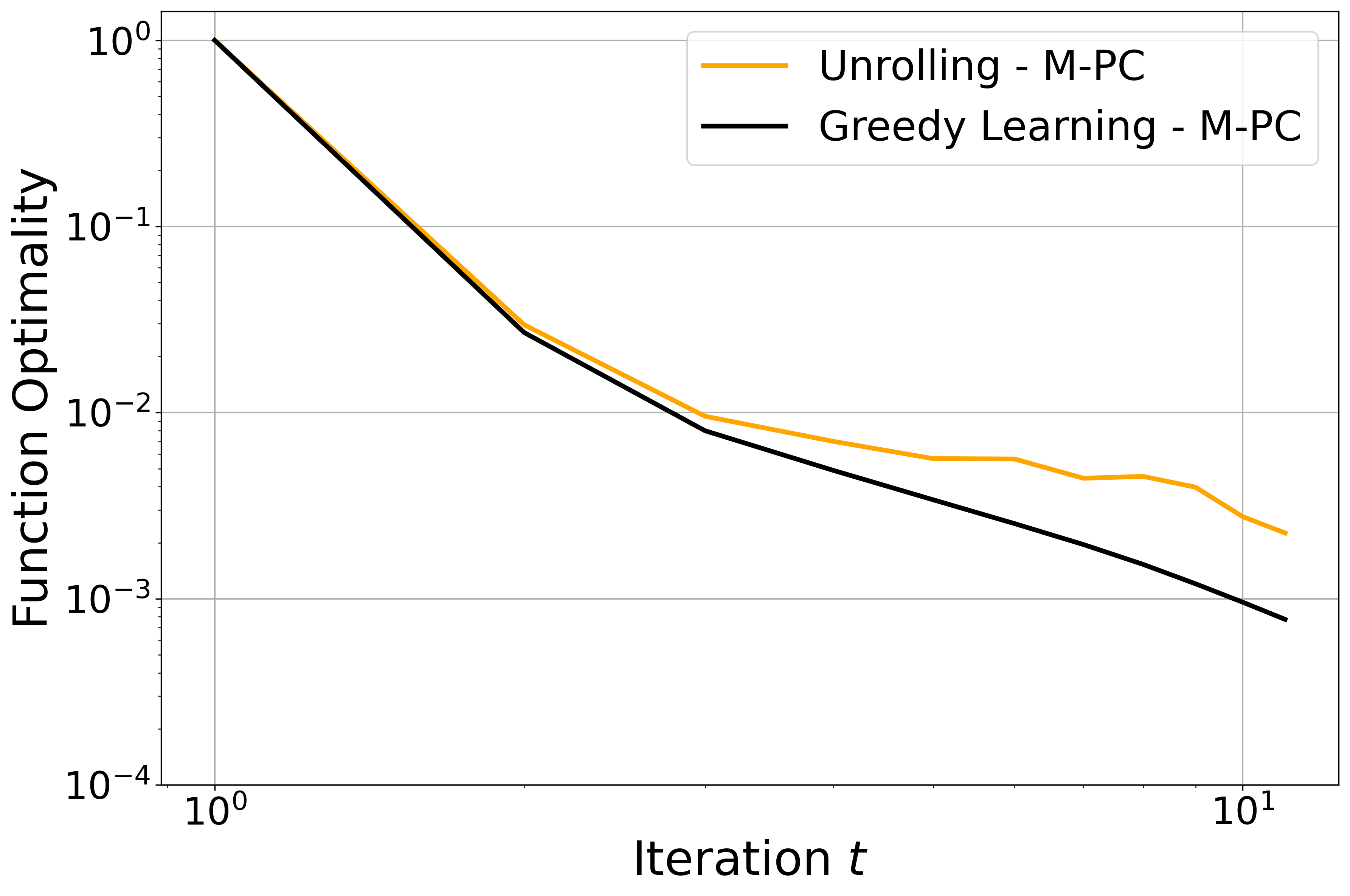} % Replace with your image file
        
        %\caption{}
        \label{unroll-perf}
    \end{subfigure}
    \hspace{15px}
    \begin{subfigure}[t]{0.38\textwidth}
        \centering
        \includegraphics[width=\textwidth]{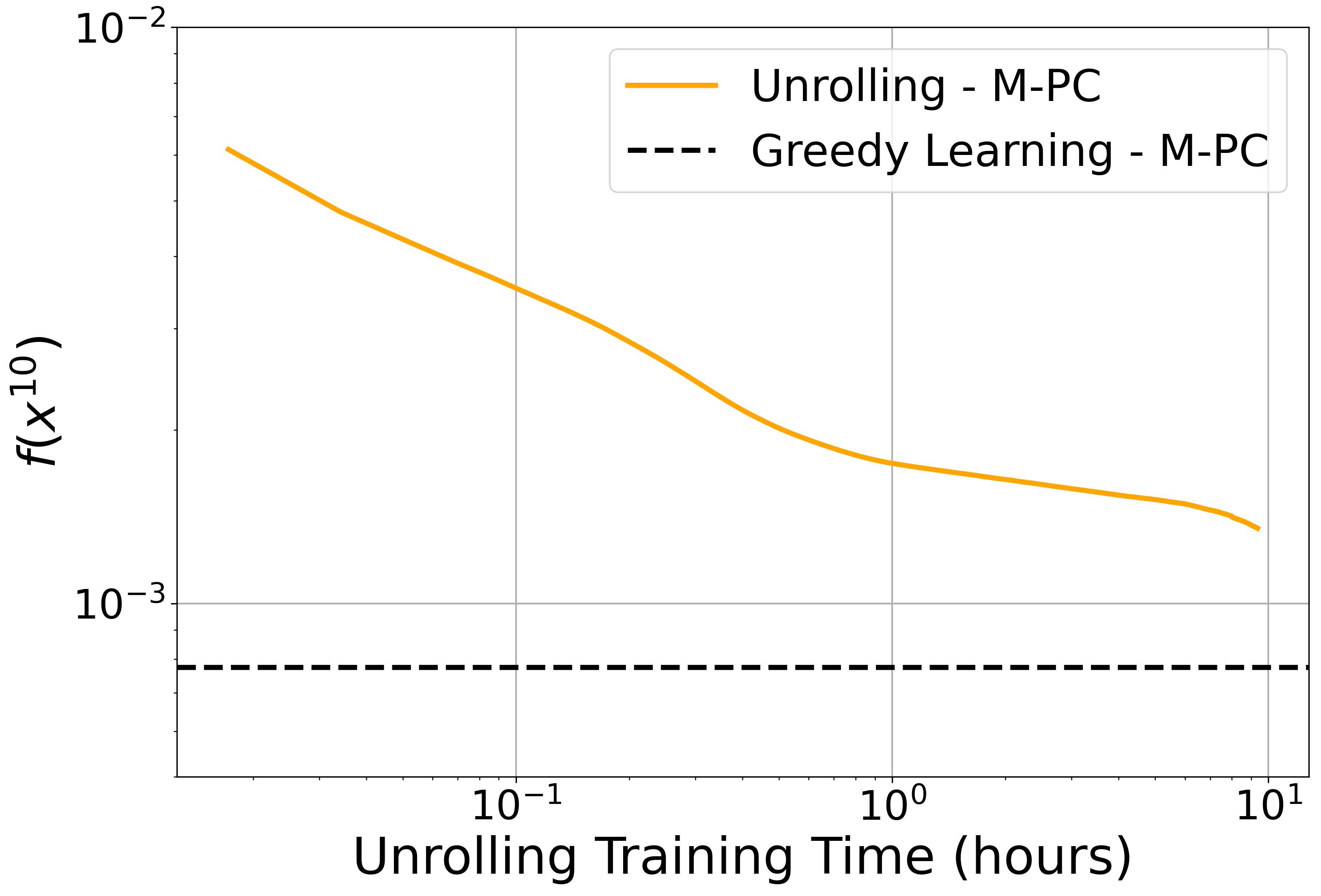} % Replace with your image file
        
        %\caption{}
        \label{unroll-train-time}
    \end{subfigure}
    
    \caption{Left: Performance of the learned unrolled algorithm versus the greedy learned algorithm on test data for Problem 1 over the first $10$ iterations. Right: The algorithm learned using unrolling evaluated on the test data at iteration $10$ versus training time, with the Greedy L2O learned algorithm value at iteration $10$ for comparison. We see that using unrolling takes a significant amount of time to reach a mediocre performance.}
    \label{fig--unroll}
\end{figure}

\paragraph{Greedy learning vs unrolling} We also compare the time taken for training with the greedy learning approach versus unrolling. For unrolling, we fix $T=10$ iterations and jointly learn the parameters $(\theta_0, \phi_0), \dots, (\theta_{T}, \phi_t)$ (all initialized as $(\tilde{\theta}, \tilde{\phi})$) in the update rule \cref{eq--general-update} with the M-PC parametrization. The same training dataset as the greedy method is used with the loss function defined in \cref{eq--finite-unroll-l2o} with $\omega_t=1$ for all $t$. Parameters are learned using Adam \cite{kingma2014adam}, with the learning rate equal to $2\times 10^{-5}$ selected via grid search and a batch size of $4$. The unrolling method was trained for approximately $10$ hours, whereas our method learns these parameters in just over $2$ hours. \cref{fig--unroll} shows that the learned parameters using unrolling achieve significantly worse performance than those learned using our method, despite a greater training time.

\begin{figure}
    \centering
    \centering
    \begin{subfigure}[t]{0.43\textwidth}
        \centering
        \includegraphics[width=\textwidth]{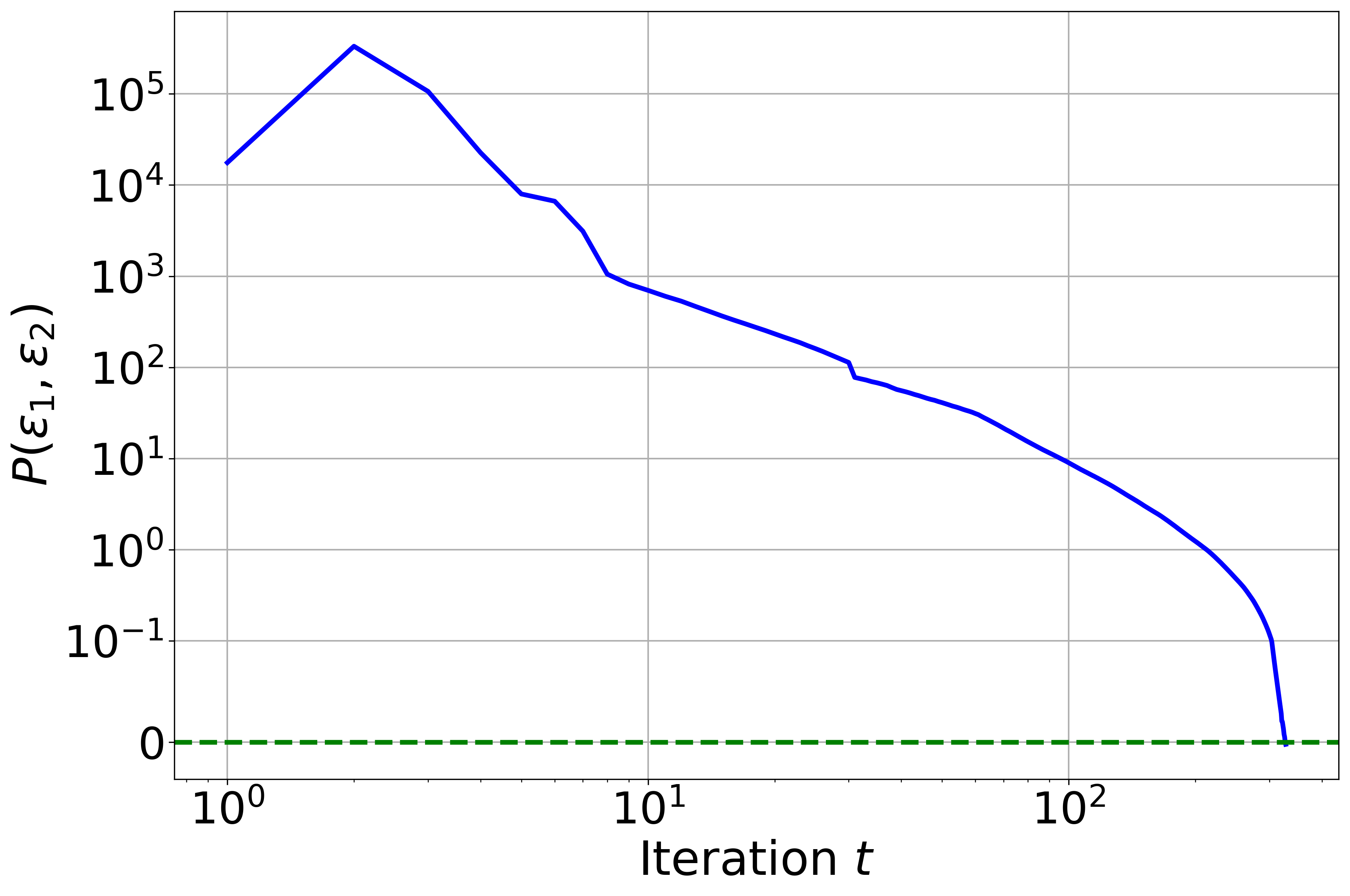} % Replace with your image file
        
        %\caption{}
        \label{fig:P-conv}
    \end{subfigure}
    \hspace{15px}
    \begin{subfigure}[t]{0.49\textwidth}
        \centering
        \includegraphics[width=\textwidth]{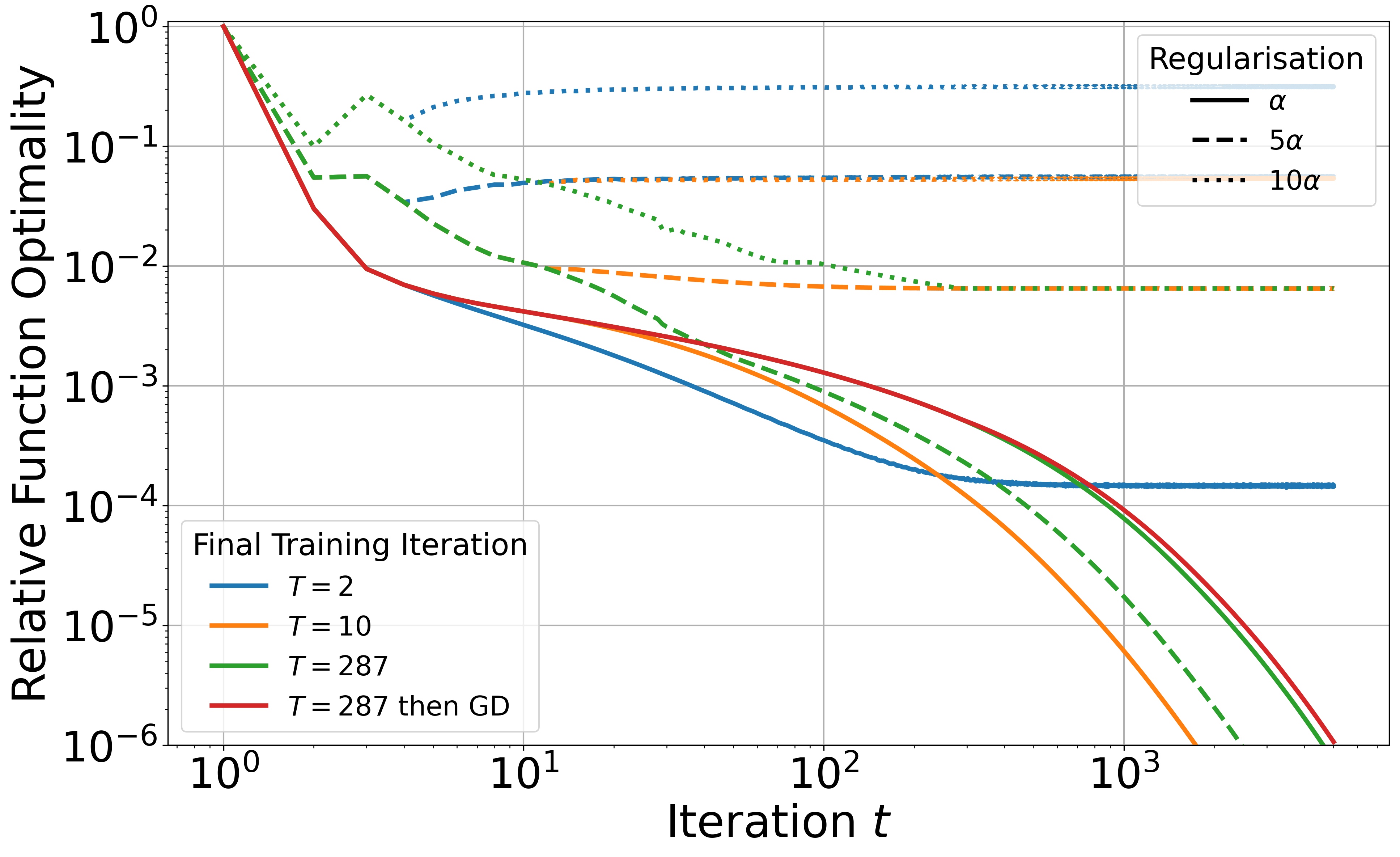} % Replace with your image file
        
        %\caption{}
        \label{fig:reg-plot}
    \end{subfigure}
    
    \caption{Left: The value of the polynomial $P$ \cref{eq::P} in training for $t \in \{0, \dots, 287\}$ for Problem 1 with $\lambda_t = 10^{-4}, \mu_t = 5 \times 10^{-2}$. Right: Objective function values in testing for \cref{alg--test-algo} with final training iteration $T=2$, $T=10$, and $T=287$. We see that the condition $P<0$ leads to guaranteed convergence (solid green curve), but we also see an example of convergence when $P > 0$ (solid orange curve). We also show the change in performance if the regularization parameter $\alpha$ used during training is scaled, with $5\alpha$ with corresponding $L=2.2 < 2L_{\text{train}} = 2.48$ leading to convergence only for $T=287$, and $10\alpha$ with corresponding $L=4.3 > 2L_{\text{train}}$ leading to non-convergence for all learned algorithms.}
    \label{fig:reg-plot-all}
\end{figure}

\paragraph{Regularization} To show how the convergence guarantees can be obtained from \cref{thm--reg-conv}, we consider the M-PC parametrization learned with $\lambda_t = 10^{-4}, \mu_t = 5 \times 10^{-2}$. \cref{fig:reg-plot-all} shows how the inequality $P < 0$, for $P$ defined in \cref{eq::P}, is eventually satisfied during training at iteration $T=287$. \cref{fig:reg-plot-all} illustrates that when the condition on $P$ is satisfied, then convergence holds as $t \to \infty$ ($T = 287$). We also see that if this condition is not satisfied, we can observe both convergence ($T = 10$) and non-convergence ($T = 2$). Note that the test performance when $T=10$ is greater than when $T=287$, this may be due to the fact that at $T=10$, the regularization has less effect on parameter learning, so parameters are more data-adaptive. The learned preconditioner when $T=10$ is more data-adaptive due to the increasing effect of constant regularization over training iterations, as shown in \cref{lemma--theta-t-conv}. We also see that switching to gradient descent iterations with step size $1/L$ achieves slightly slower empirical convergence than the learned algorithm, meaning that the parameters learned at iteration $T=287$ are still data-adaptive. \cref{fig:reg-plot-all} also illustrates that when $L < 2L_{\text{train}}$ we can obtain convergence of \cref{alg--test-algo} as seen in \cref{thm--reg-conv}. In this instance, we only see convergence when training is stopped at $T=287$, for which we inherit the convergence guarantees. We also see that choosing a regularization parameter so that $L > 2L_{\text{train}}$ may not converge in this instance.

% \rcoment{The learned scalar quickly diverges as the final step size learned is greater than $2/L$, despite decreasing for values that fluctuate above and below $2/L$. Similarly, the learned pointwise parametrization diverges even more quickly, due to the wide range of parameter values.}

% \textbf{Ablation study: size of learned kernels} \cref{fig--kernel-sizes} shows that many of the learned convolutional algorithms outperform NAG for Problem 1. We see that the $5 \times 5$ kernels significantly outperform the NAG kernels and perform similarly to the $7 \times 7$ kernels.Furthermore, we see similar performance for the $11 \times 11$ kernels and the $96 \times 96$ kernels. 
% \begin{figure}[h!]
%     \centering
%     \centering
%     \includegraphics[width=.4\textwidth]{blur_images/kernel_sizes.png} % Replace 
    
%     \caption{Test performance of different kernel sizes in the convolutional parametrisation, averaged over the test dataset in Problem 1. Tested kernel sizes are $3 \times 3$, $5 \times 5$, $7 \times 7$, $11 \times 11$, $3 \times 3$, $96 \times 96$.}
%     \label{fig--kernel-sizes}
% \end{figure}

\subsection{Results for Problem 2}
\label{sec::ct}

% \textbf{Training details.} Greedy training was performed with $\lambda_t, \mu_t = 0$ for all iterations $t$. For the M-PC parametrization, parameters were learned up to iteration $T=50$, and for the M-PS parametrization, parameters were learned up to iteration $T=100$. The total time for training M-PS was about $5$ hours, and M-PC took approximately $33$ hours.

\begin{figure}[h!]
    \centering
    % \begin{subfigure}[t]{.49\textwidth}
    %     \centering
    %     \includegraphics[width=\textwidth]{images/ct/ct_kernels.png}
        
    %     \caption{}
    %     \label{fig:ct-kernels}
    % \end{subfigure}
    % \begin{subfigure}[t]{.45\textwidth}
    %     \centering
    %     \includegraphics[width=\textwidth]{images/ct/ct_pw.png}
        
    %     \caption{}
    %     \label{fig:ct-diags}
    % \end{subfigure}
    % \newline 
    % \begin{subfigure}[t]{0.33\textwidth}
    %     \centering
    %     \includegraphics[width=\textwidth]{images/ct/ct_scalars.png} % Replace with your image file
        
    %     \caption{}
    %     \label{fig:ct-scalars}
    % \end{subfigure}
    \begin{subfigure}[t]{.45\textwidth}
        \centering
        \includegraphics[width=\textwidth]{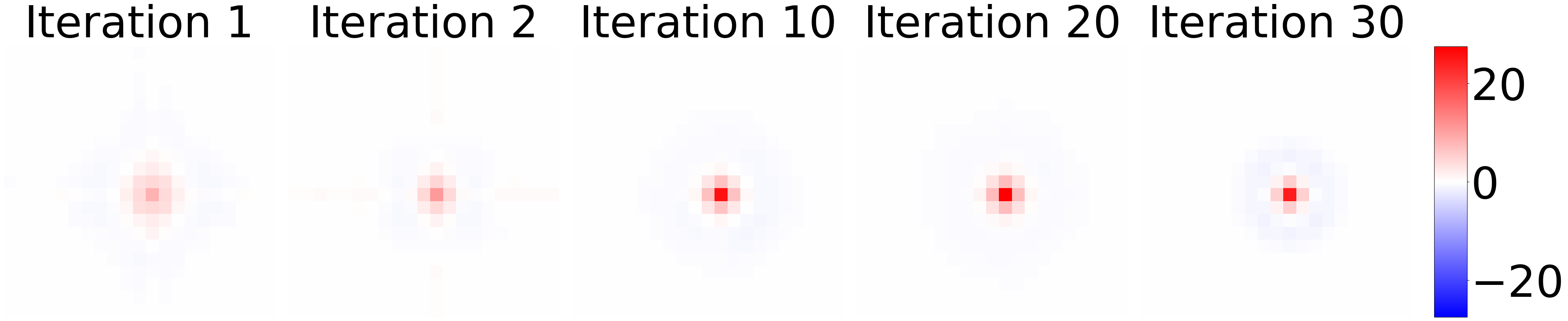}
        
        \caption{$21 \times 21$ crop of $\theta_t$ M-PC kernels.}
        \label{fig:ct-kernels-p}
    \end{subfigure}
    \begin{subfigure}[t]{.45\textwidth}
        \centering
        \includegraphics[width=\textwidth]{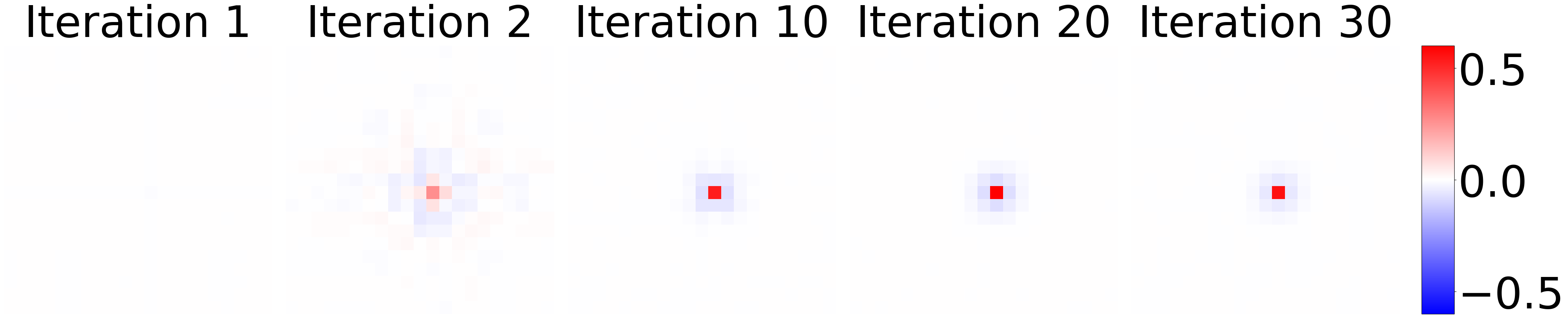}
        
        \caption{$21 \times 21$ crop of $\phi_t$ M-PC kernels.}
        \label{fig:ct-kernels-q}
    \end{subfigure}
    \begin{subfigure}[t]{.4\textwidth}
        \centering
        \includegraphics[width=\textwidth]{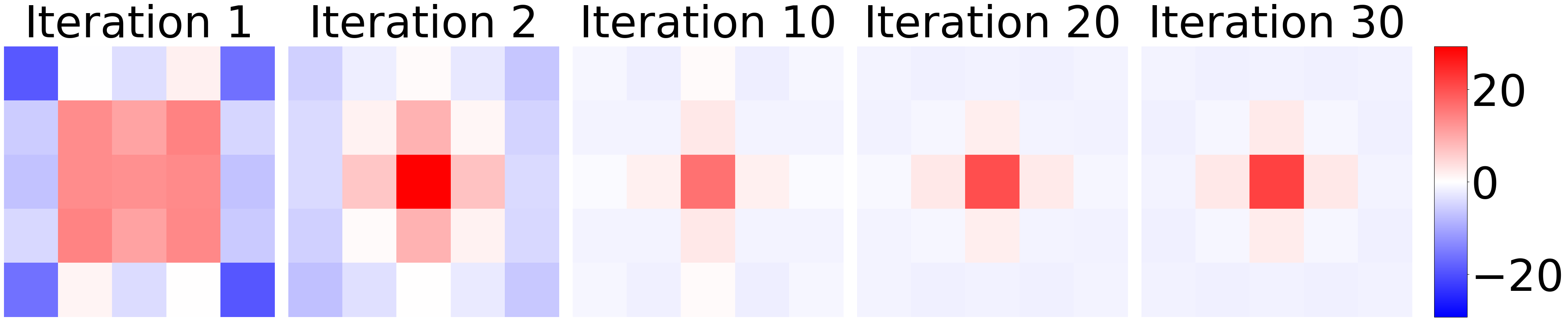}
        
        \caption{Learned $5 \times 5$ $\theta_t$ M-PC kernels.}
        \label{fig:ct-kernels-p-5}
    \end{subfigure}
    \begin{subfigure}[t]{.4\textwidth}
        \centering
        \includegraphics[width=\textwidth]{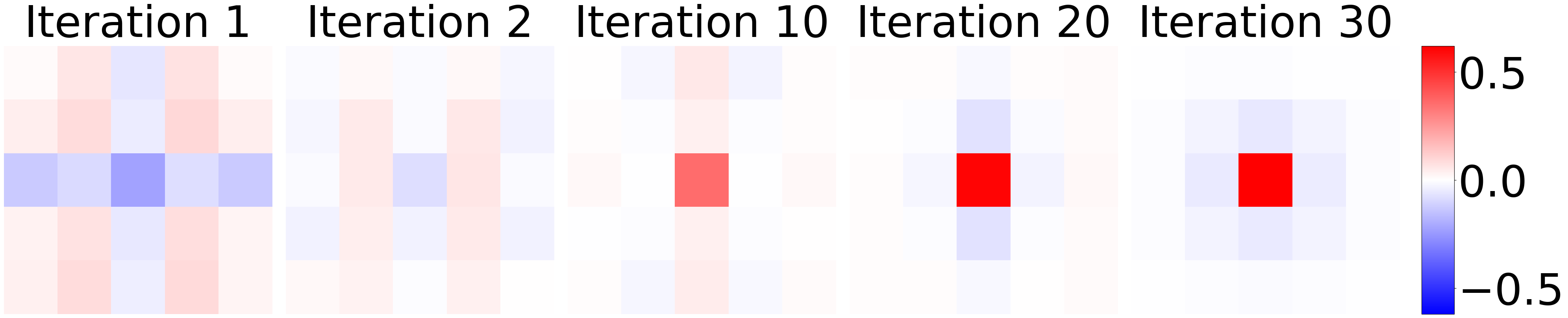}
        
        \caption{Learned $5 \times 5$ $\phi_t$ M-PC kernels.}
        \label{fig:ct-kernels-q-5}
    \end{subfigure}
    \begin{subfigure}[t]{0.4\textwidth}
        \centering
        \includegraphics[width=\textwidth]{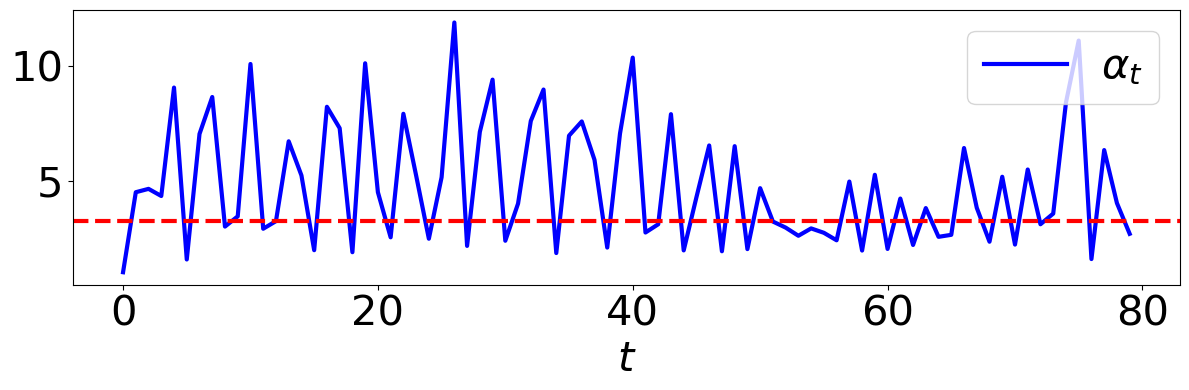} % Replace with your image file
        
        \caption{Learned $\alpha_t$ M-PS scalars.}
        \label{fig:ct-alpha-p}
    \end{subfigure}   
    \begin{subfigure}[t]{0.4\textwidth}
        \centering
        \includegraphics[width=\textwidth]{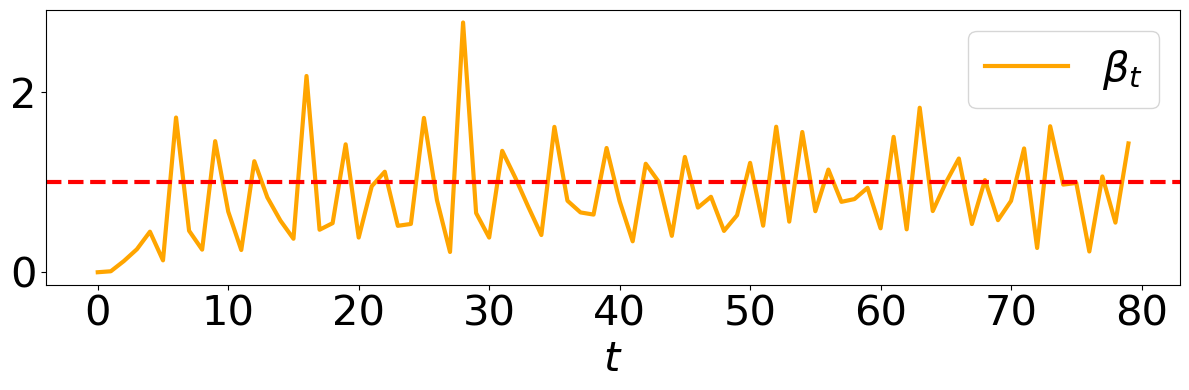} % Replace with your image file
        
        \caption{Learned $\beta_t$ M-PS scalars.}
        \label{fig:ct-alpha-q}
    \end{subfigure}   
    \caption{Learned preconditioners for Problem 2.}
    \label{fig--CT-obs-recon-scalars}
\end{figure}

\paragraph{Visualizing learned preconditioners} \cref{fig:ct-kernels-p} and \cref{fig:ct-kernels-q} consider M-PC, where the learned kernels are weighted heavily towards the center, and with $\phi_t$ having pixel values less than $0.5$, as seen in Problem 1. The $5 \times 5$ learned kernels in \cref{fig:ct-kernels-p-5} and \cref{fig:ct-kernels-q-5} lie within roughly the same range as for the $256 \times 256$ learned kernels. \cref{fig:ct-alpha-p} and \cref{fig:ct-alpha-p} show that the learned scalars oscillate around the Heavy Ball values $\alpha^*$ and $\beta^*$, rather than approach these values monotonically.

\begin{figure}[h!]
    \centering
    \begin{subfigure}[t]{0.42\textwidth}
        \centering
        \includegraphics[width=\textwidth]{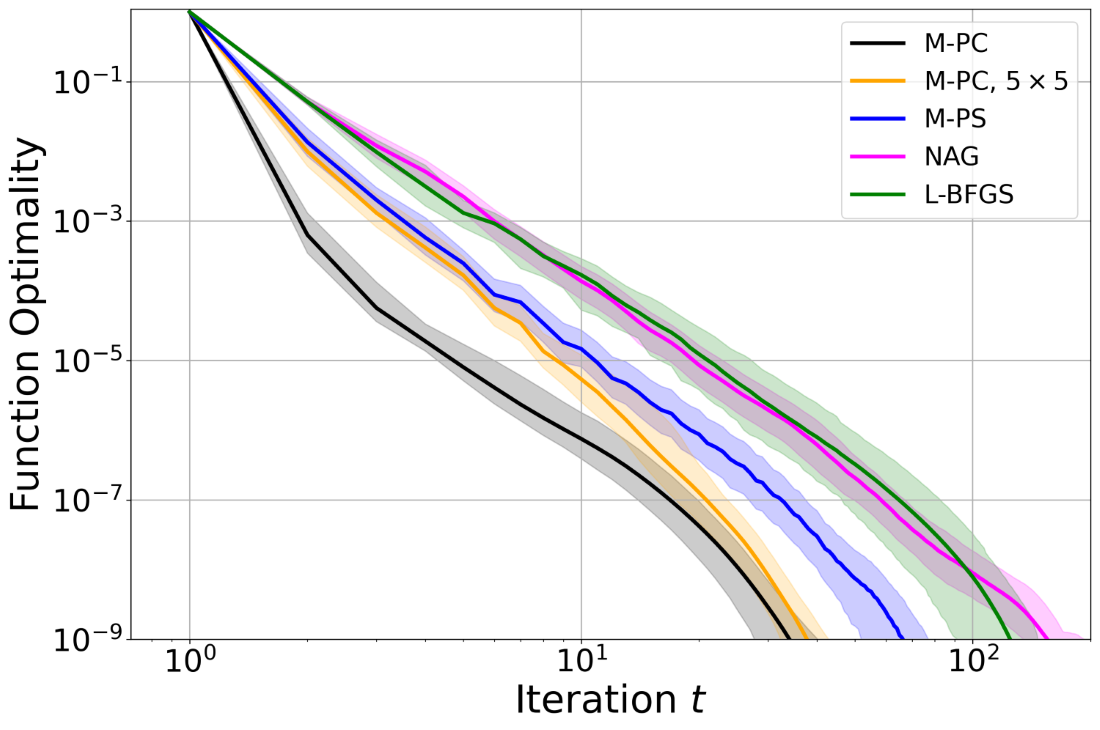} % Replace with your image file
        
        %\caption{}
        \label{ct-test-iter}
    \end{subfigure}
    \hspace{10px}
    \begin{subfigure}[t]{0.42\textwidth}
        \centering
        \includegraphics[width=\textwidth]{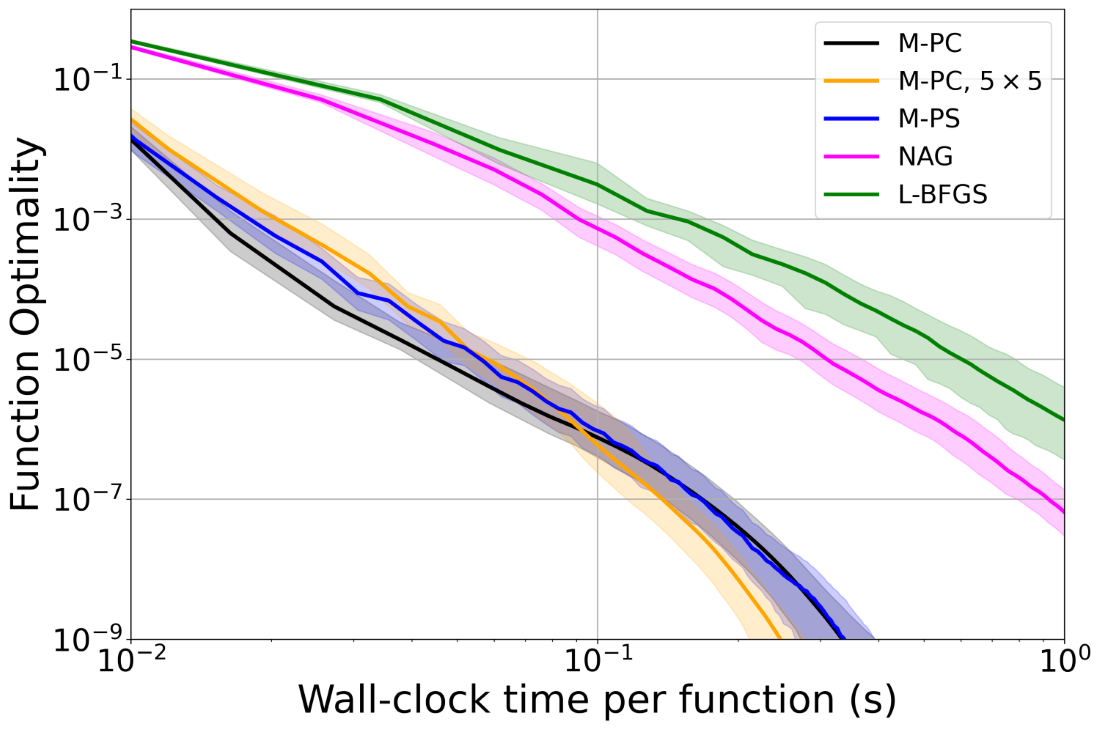} % Replace with your image file
        
        %\caption{}
        \label{ct-test-time}
    \end{subfigure}
    
    \caption{Performance of the proposed method with M-PC and M-PS parametrizations versus benchmark non-learned algorithms for Problem 2. We see that the learned algorithms significantly outperform the benchmark algorithms. Left: Test performance versus benchmark optimization algorithms. Right: Wall-clock time test performance versus benchmark optimization algorithms.}
    \label{fig--CT-test-big}
\end{figure}

\begin{figure}[h!]
    \centering
    \begin{subfigure}[t]{0.27\textwidth}
        \centering
        \includegraphics[width=\textwidth]{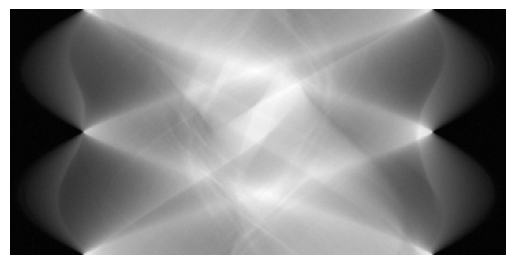} % Replace with your image file
        
        \caption{Observation.}
        \label{fig--CT-obs-big}
    \end{subfigure}
    \begin{subfigure}[t]{0.2\textwidth}
        \centering
        \includegraphics[width=.75\textwidth]{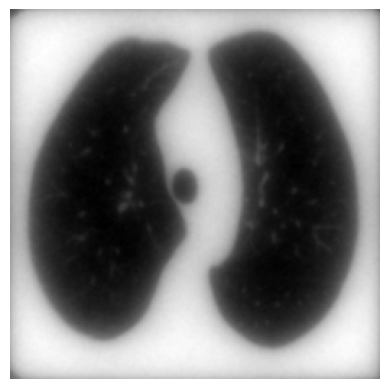} % Replace with your image file
        
        \caption{NAG iteration $15$.}
        \label{ct_recon_nag-big}
    \end{subfigure}
    \begin{subfigure}[t]{0.2\textwidth}
        \centering
        \includegraphics[width=.75\textwidth]{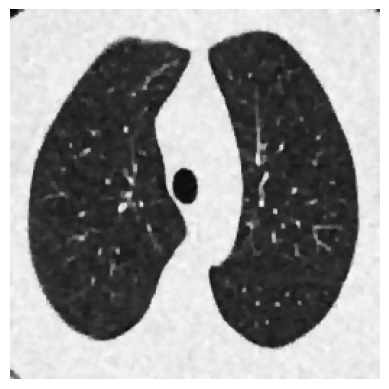} % Replace with your image file
        
        \caption{M-PC iteration $15$.}
        \label{ct_recon_conv-big}
    \end{subfigure}
    \begin{subfigure}[t]{0.2\textwidth}
        \centering
        \includegraphics[width=.75\textwidth]{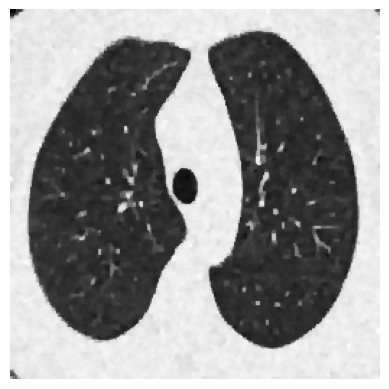} % Replace with your image file
        
        \caption{Final reconstruction.}
        \label{ct_gound_truth}
    \end{subfigure}
    \caption{A comparison of image reconstructions at iteration $15$ for Problem 2. The learned M-PC algorithm provides significantly better reconstruction at iteration $15$.}
    \label{fig:ct-recon}
\end{figure}

\paragraph{Learned algorithm performance} Similar to Problem 1, \cref{fig--CT-test-big} shows that the learned M-PC parametrization outperforms NAG and L-BFGS on the CT test data, reaching a tolerance of $10^{-7}$ in an average of approximately $20$ iterations, compared with about $90$ for both L-BFGS and NAG. Using learned kernels of size $5 \times 5$ in M-PC leads to slower convergence initially compared to the $256 \times 256$ kernel, but this difference decreases as iterations increase. Furthermore, the performance with respect to wall-clock time of the learned algorithms greatly outperforms NAG and L-BFGS. However, due to the computational cost of convolution with $256 \times 256$ kernels, as iterations increase, the learned scalars achieve similar performance with respect to time for larger iterations. However, we see that learning smaller kernels maintains an increased performance per unit time over the scalar parametrization. \cref{fig:ct-recon} show qualitatively that the learned convolutional algorithm achieves a good reconstruction faster than NAG.

\begin{figure}[h!]
    \centering
    \begin{subfigure}[t]{0.38\textwidth}
        \centering
        \includegraphics[width=\textwidth]{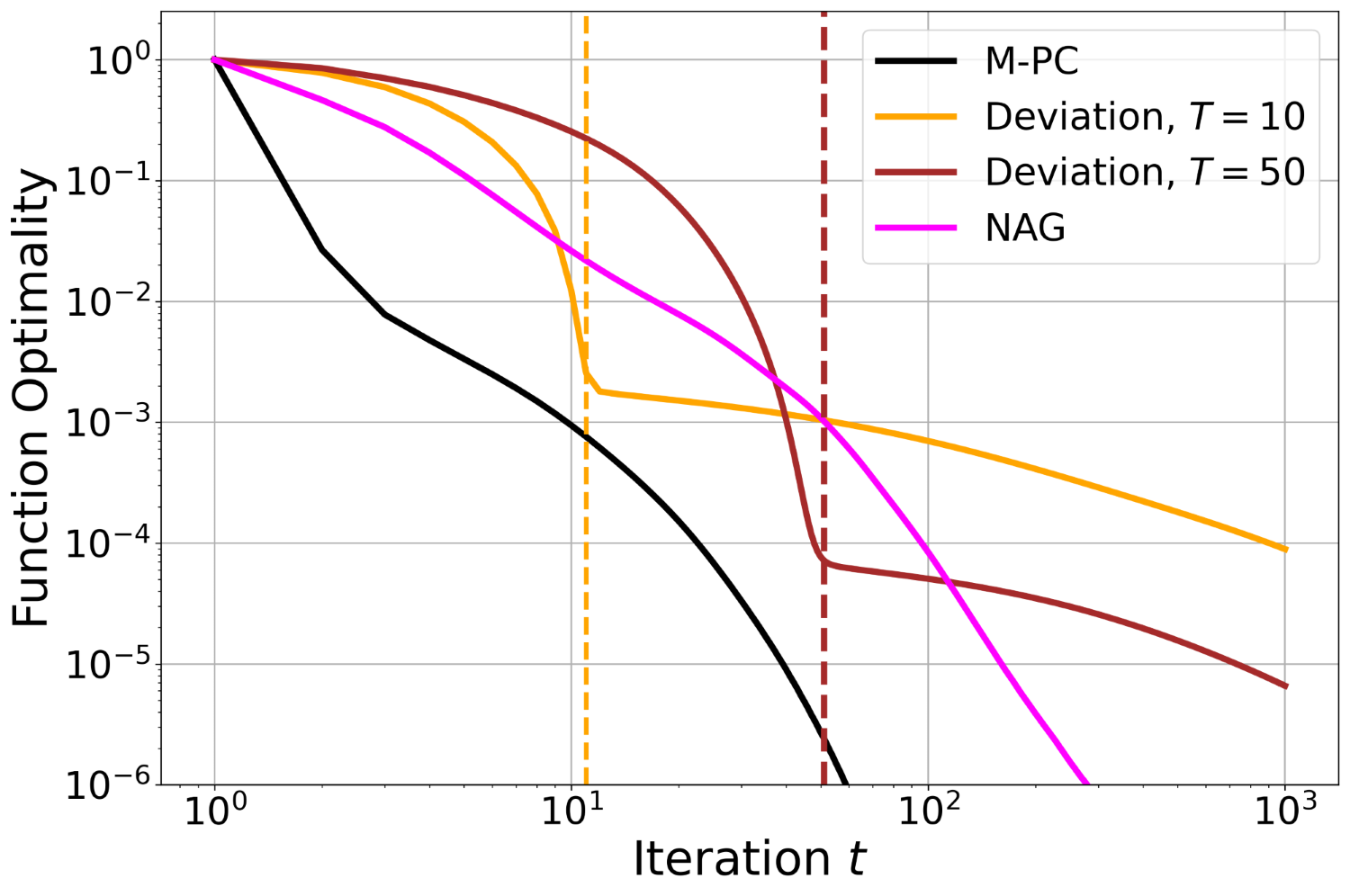} % Replace with your image file
        
        %\caption{}
        \label{fig5b}
    \end{subfigure}
    \hspace{15px}
    \begin{subfigure}[t]{0.38\textwidth}
        \centering
        \includegraphics[width=\textwidth]{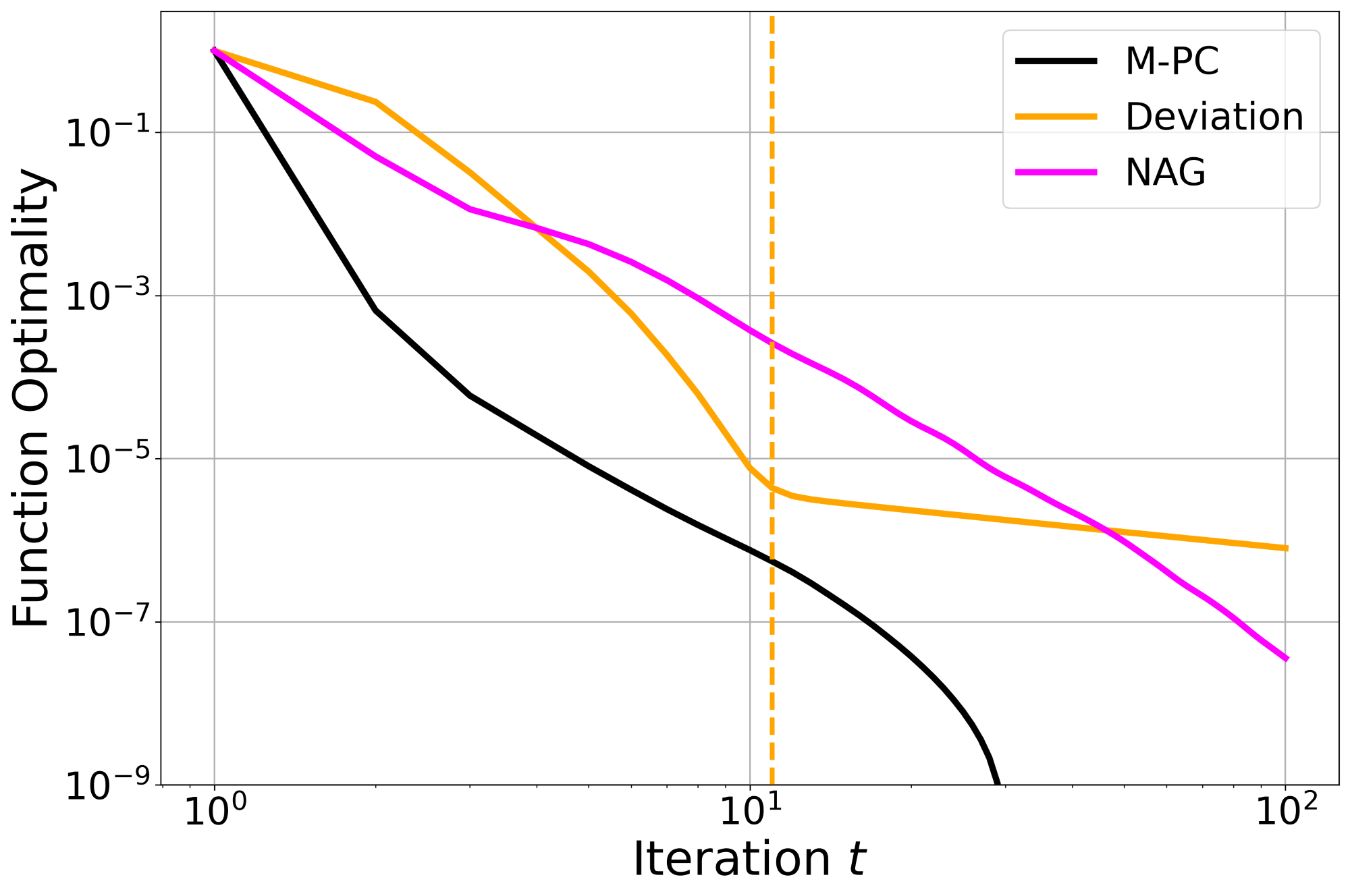} % Replace with your image file
        
        %\caption{}
        \label{fig5c}
    \end{subfigure}
    
    \caption{A comparison to the deviation-based approach \cite{banert2024accelerated}. Left: Performance on Problem 1, for two different unroll lengths, $T=10$ (orange dashed line), $T=50$ (brown dashed line). For $T=10$ we trained for $1$ hour, and for $T=50$, we trained for $4.5$ hours. Right: Performance on Problem 2 with an unroll length of $T=10$ (orange dashed line). This is the maximum unroll length considered for CT due to memory constraints. Training took $4$ hours. We see the learned M-PC algorithm outperforms the deviation-based approach on test data for both Problem 1 and Problem 2.}
    \label{fig--banert}
\end{figure}

\paragraph{Comparison to an existing L2O method} We also compare to the deviation-based Learning to Optimize approach \cite{banert2024accelerated} for Problem 1 and Problem 2. In \cref{fig--banert} we see that our learned M-PC parametrization outperforms the deviation-based approach within training iterations, as well as being able to train over more iterations. It is also worth noting that the performance of the learned deviation-based approach slows for iterations after the final training iteration $T$.

\subsection{Results for Problem 3} \label{sec::nonconvex-exp}

% \textbf{Training details.} Using the greedy method, we learn the convolution parametrization with momentum M-PC  up to iteration $T=100$ and the scalar parametrization with momentum up to iteration $T=150$, using $\lambda_t, \mu_t = 0$ for all $t$. The total training time for M-PS was about $4.5$ hours, and for M-PC it was $31$ hours.

\begin{figure}[h!]
    \centering
    \begin{subfigure}[t]{0.38\textwidth}
        \centering
        \includegraphics[width=\textwidth]{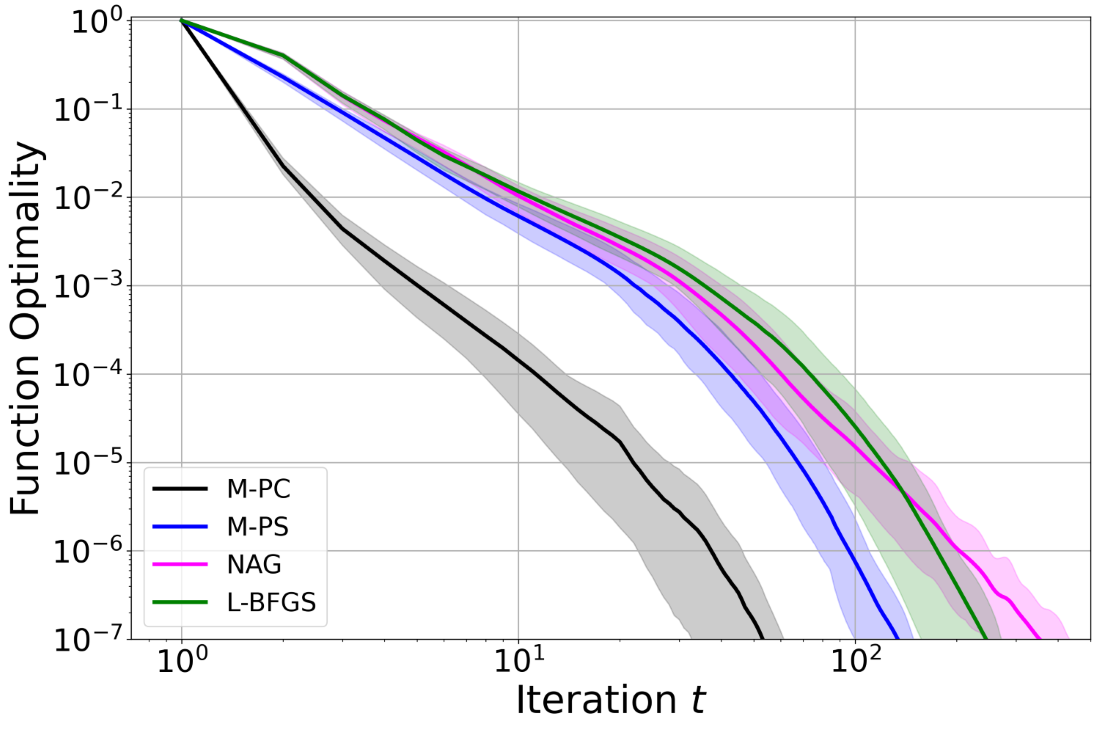} % Replace with your image file
        %\caption{}
        \label{wcrr-iter}
    \end{subfigure}
    \hspace{15px}
    \begin{subfigure}[t]{0.38\textwidth}
        \centering
        \includegraphics[width=\textwidth]{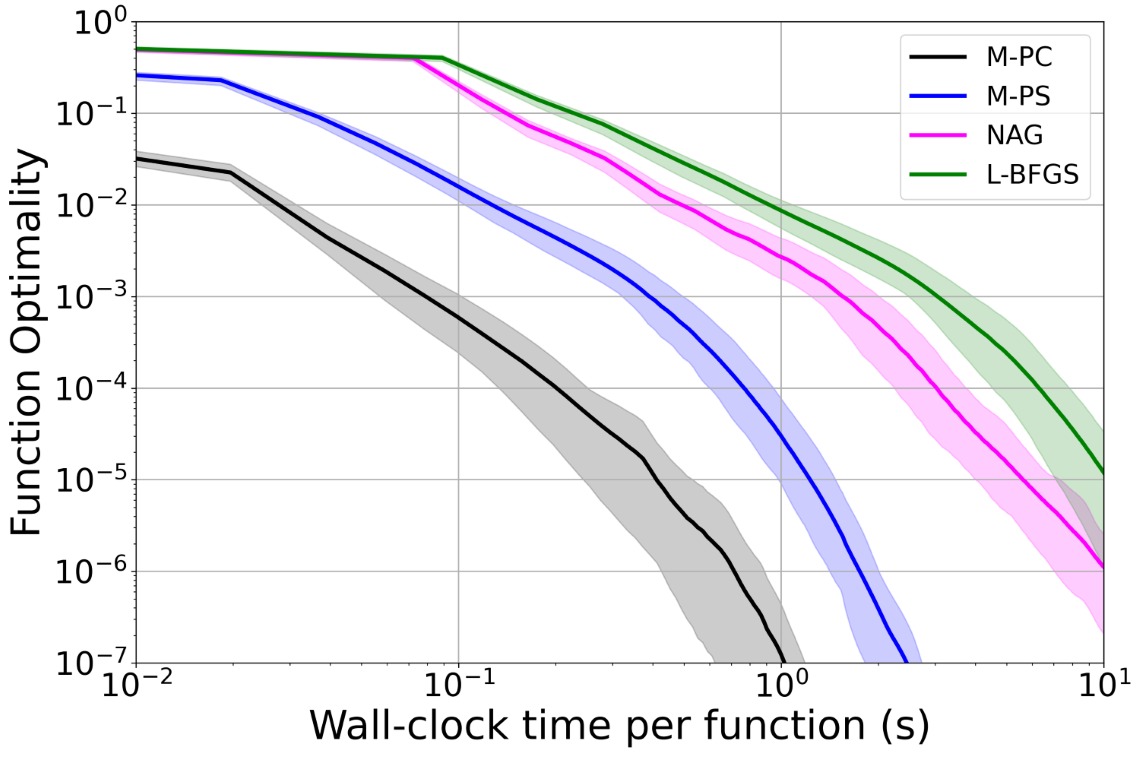} % Replace with your image file
        %\caption{}
        \label{wcrr-time}
    \end{subfigure}
    \hspace{10px}

    \caption{Test performance of the proposed method with M-PC and M-PS parametrizations versus benchmark non-learned algorithms for Problem 3 with a nonconvex regularizer. We see that the learned algorithms significantly outperform the benchmark algorithms. Intervals around each mean represent maximum and minimum values over the dataset. }
    \label{fig--wcrr-blur-test}
\end{figure}

\cref{fig--wcrr-blur-test} shows that the learned M-PC parametrization outperforms NAG and L-BFGS on the test data. Furthermore, the performance with respect to wall-clock time of the learned algorithms greatly outperforms NAG and L-BFGS. Due to the increased cost of convolution, the learned M-PS algorithm performs more similarly to M-PC, but the learned kernels still outperform the other algorithms. 

\section{Conclusions} \label{section--future}

Our contribution is a novel L2O approach for minimizing unconstrained problems with differentiable objective functions. Our method employs a greedy strategy to learn linear operators at each iteration of an optimization algorithm, so that device-memory requirements are constant with respect to the number of training iterations and the training horizon \(T\) does not need to be selected a priori. Parameter learning in our framework corresponds to solving convex optimization problems whenever the objective functions are convex, enabling the use of efficient algorithms. These factors allow training over a large number of iterations, which would otherwise be prohibitively expensive with standard unrolling-based approaches. Furthermore, we prove convergence to stationarity on the training set even when the learned preconditioners are neither symmetric nor positive definite, and obtain convergence guarantees on a more general function class through regularization. Numerical results on imaging inverse problems demonstrate that our approach with a novel convolutional parameterization outperforms competing methods, including NAG and L-BFGS, on convex and nonconvex problems.

\bibliography{references}

@article{armijo1966minimization,
  author  = {Armijo, Larry},
  title   = {Minimization of functions having {L}ipschitz continuous first partial derivatives},
  journal = {Pacific Journal of Mathematics},
  volume  = {16},
  number  = {1},
  pages   = {1--3},
  year    = {1966}
}

@book{briggs1995dft,
  author    = {Briggs, William L. and Henson, Van Emden},
  title     = {The {DFT}: An Owner's Manual for the Discrete Fourier Transform},
  publisher = {SIAM},
  address   = {Philadelphia, PA},
  year      = {1995}
}

@article{qu2024optimal,
  author  = {Qu, Zhaonan and Gao, Wenzhi and Hinder, Oliver and Ye, Yinyu and Zhou, Zhengyuan},
  title   = {Optimal diagonal preconditioning},
  journal = {Operations Research},
  volume  = {73},
  number  = {3},
  pages   = {1479--1495},
  year    = {2025},
  doi     = {10.1287/opre.2022.0592}
}

@article{duchi2011adaptive,
  author  = {Duchi, John and Hazan, Elad and Singer, Yoram},
  title   = {Adaptive subgradient methods for online learning and stochastic optimization},
  journal = {Journal of Machine Learning Research},
  volume  = {12},
  number  = {61},
  pages   = {2121--2159},
  year    = {2011}
}

@article{mcmahan2010adaptive,
  author  = {McMahan, H. Brendan and Streeter, Matthew},
  title   = {Adaptive bound optimization for online convex optimization},
  journal = {arXiv preprint arXiv:1002.4908},
  year    = {2010}
}

@article{wolfe1969convergence,
  author  = {Wolfe, Philip},
  title   = {Convergence conditions for ascent methods},
  journal = {SIAM Review},
  volume  = {11},
  number  = {2},
  pages   = {226--235},
  year    = {1969}
}

@article{beck2009fast,
  author  = {Beck, Amir and Teboulle, Marc},
  title   = {A fast iterative shrinkage-thresholding algorithm for linear inverse problems},
  journal = {SIAM Journal on Imaging Sciences},
  volume  = {2},
  number  = {1},
  pages   = {183--202},
  year    = {2009}
}

@inproceedings{orabona2016coin,
  author    = {Orabona, Francesco and P{\'a}l, D{\'a}vid},
  title     = {Coin betting and parameter-free online learning},
  booktitle = {Advances in Neural Information Processing Systems 29},
  pages     = {577--585},
  year      = {2016},
  publisher = {Curran Associates, Inc.}
}

@article{hazan2016introduction,
  author  = {Hazan, Elad},
  title   = {Introduction to online convex optimization},
  journal = {Foundations and Trends{\textregistered} in Optimization},
  volume  = {2},
  number  = {3--4},
  pages   = {157--325},
  year    = {2016}
}

@book{polyak1987,
  author    = {Polyak, Boris T.},
  title     = {Introduction to Optimization},
  series     = {Translations Series in Mathematics and Engineering},
  publisher = {Optimization Software, Inc.},
  address   = {New York},
  year      = {1987}
}

@article{nesterov1983,
  author  = {Nesterov, Yurii},
  title   = {A method for solving the convex programming problem with convergence rate {$O(1/k^2)$}},
  journal = {Soviet Mathematics Doklady},
  volume  = {27},
  number  = {2},
  pages   = {372--376},
  year    = {1983}
}

@article{gower2020variancereduced,
  author  = {Gower, Robert M. and Schmidt, Mark and Bach, Francis and Richt{\'a}rik, Peter},
  title   = {Variance-reduced methods for machine learning},
  journal = {Proceedings of the IEEE},
  volume  = {108},
  number  = {11},
  pages   = {1968--1983},
  year    = {2020}
}

@article{banert2020,
  author  = {Banert, Sebastian and Ringh, Axel and Adler, Jonas and Karlsson, Johan and {\"O}ktem, Ozan},
  title   = {Data-driven nonsmooth optimization},
  journal = {SIAM Journal on Optimization},
  volume  = {30},
  number  = {1},
  pages   = {102--131},
  year    = {2020}
}

@article{monga2021algorithm,
  author  = {Monga, Vishal and Li, Yuelong and Eldar, Yonina C.},
  title   = {Algorithm unrolling: Interpretable, efficient deep learning for signal and image processing},
  journal = {IEEE Signal Processing Magazine},
  volume  = {38},
  number  = {2},
  pages   = {18--44},
  year    = {2021}
}

@incollection{nocedal2006quasi,
  author    = {Nocedal, Jorge and Wright, Stephen J.},
  title     = {Quasi-{N}ewton methods},
  booktitle = {Numerical Optimization},
  publisher = {Springer},
  address   = {New York},
  edition   = {2},
  pages     = {135--163},
  year      = {2006}
}

@article{2016RLL2O,
  author  = {Li, Ke and Malik, Jitendra},
  title   = {Learning to optimize},
  journal = {arXiv preprint arXiv:1606.01885},
  year    = {2016}
}

@inproceedings{heaton2023safeguarded,
  author    = {Heaton, Howard and Chen, Xiaohan and Wang, Zhangyang and Yin, Wotao},
  title     = {Safeguarded learned convex optimization},
  booktitle = {Proceedings of the AAAI Conference on Artificial Intelligence},
  volume    = {37},
  number    = {6},
  pages     = {7848--7855},
  year      = {2023},
  publisher = {AAAI Press}
}

@inproceedings{andrychowicz2016learning,
  author    = {Andrychowicz, Marcin and Denil, Misha and Gomez, Sergio and Hoffman, Matthew W. and Pfau, David and Schaul, Tom and Shillingford, Brendan and de Freitas, Nando},
  title     = {Learning to learn by gradient descent by gradient descent},
  booktitle = {Advances in Neural Information Processing Systems},
  volume    = {29},
  pages     = {3981--3989},
  year      = {2016}
}

@inproceedings{LISTA,
  author    = {Gregor, Karol and LeCun, Yann},
  title     = {Learning fast approximations of sparse coding},
  booktitle = {Proceedings of the 27th International Conference on Machine Learning ({ICML}-10)},
  pages     = {399--406},
  year      = {2010}
}

@book{Beck2014,
  author    = {Beck, Amir},
  title     = {Introduction to Nonlinear Optimization: Theory, Algorithms, and Applications with {MATLAB}},
  publisher = {SIAM},
  address   = {Philadelphia, PA},
  year      = {2014}
}

@article{tan2023mirror,
  author  = {Tan, Hong Ye and Mukherjee, Subhadip and Tang, Junqi and Sch{\"o}nlieb, Carola-Bibiane},
  title   = {Data-driven mirror descent with input-convex neural networks},
  journal = {SIAM Journal on Mathematics of Data Science},
  volume  = {5},
  number  = {2},
  pages   = {558--587},
  year    = {2023}
}

@article{liu1989limited,
  author  = {Liu, Dong C. and Nocedal, Jorge},
  title   = {On the limited memory {BFGS} method for large scale optimization},
  journal = {Mathematical Programming},
  volume  = {45},
  number  = {1},
  pages   = {503--528},
  year    = {1989}
}

@article{polyak1964,
  author  = {Polyak, Boris T.},
  title   = {Some methods of speeding up the convergence of iteration methods},
  journal = {USSR Computational Mathematics and Mathematical Physics},
  volume  = {4},
  number  = {5},
  pages   = {1--17},
  year    = {1964}
}

@book{golub2013matrix,
  author    = {Golub, Gene H. and Van Loan, Charles F.},
  title     = {Matrix Computations},
  edition   = {4},
  publisher = {Johns Hopkins University Press},
  address   = {Baltimore, MD},
  year      = {2013}
}

@book{nesterov2018lectures,
  author    = {Nesterov, Yurii},
  title     = {Lectures on Convex Optimization},
  volume    = {137},
  series    = {Springer Optimization and Its Applications},
  publisher = {Springer},
  address   = {Cham},
  year      = {2018}
}

@article{rudin1992TV,
  author  = {Rudin, Leonid I. and Osher, Stanley and Fatemi, Emad},
  title   = {Nonlinear total variation based noise removal algorithms},
  journal = {Physica D: Nonlinear Phenomena},
  volume  = {60},
  number  = {1--4},
  pages   = {259--268},
  year    = {1992}
}

@incollection{huberloss1992,
  author    = {Huber, Peter J.},
  title     = {Robust estimation of a location parameter},
  booktitle = {Breakthroughs in Statistics: Methodology and Distribution},
  pages     = {492--518},
  publisher = {Springer},
  address   = {New York},
  year      = {1992}
}

@article{tan2023boosting,
  author  = {Tan, Hong Ye and Mukherjee, Subhadip and Tang, Junqi and Sch{\"o}nlieb, Carola-Bibiane},
  title   = {Boosting data-driven mirror descent with randomization, equivariance, and acceleration},
  journal = {arXiv preprint arXiv:2308.05045},
  year    = {2023}
}

@article{SARS-CoV-2CT-scan,
  author  = {Soares, Eduardo and Angelov, Plamen and Biaso, Sarah and Froes, Michele Higa and Abe, Daniel Kanda},
  title   = {{SARS-CoV-2} {CT}-scan dataset: A large dataset of real patients {CT} scans for {SARS-CoV-2} identification},
  journal = {medRxiv},
  year    = {2020},
  doi     = {10.1101/2020.04.24.20078584}
}

@article{chambolle2016introduction,
  author  = {Chambolle, Antonin and Pock, Thomas},
  title   = {An introduction to continuous optimization for imaging},
  journal = {Acta Numerica},
  volume  = {25},
  pages   = {161--319},
  year    = {2016}
}

@article{banert2024accelerated,
  author  = {Banert, Sebastian and Rudzusika, Jevgenija and {\"O}ktem, Ozan and Adler, Jonas},
  title   = {Accelerated forward-backward optimization using deep learning},
  journal = {SIAM Journal on Optimization},
  volume  = {34},
  number  = {2},
  pages   = {1236--1263},
  year    = {2024}
}

@inproceedings{metz2019understanding,
  author    = {Metz, Luke and Maheswaranathan, Niru and Nixon, Jeremy and Freeman, Daniel and Sohl-Dickstein, Jascha},
  title     = {Understanding and correcting pathologies in the training of learned optimizers},
  booktitle = {Proceedings of the 36th International Conference on Machine Learning},
  series    = {Proceedings of Machine Learning Research},
  volume    = {97},
  pages     = {4556--4565},
  year      = {2019},
  publisher = {PMLR}
}

@inproceedings{mathstructures,
  author    = {Liu, Jialin and Chen, Xiaohan and Wang, Zhangyang and Yin, Wotao and Cai, Hanqin},
  title     = {Towards constituting mathematical structures for learning to optimize},
  booktitle = {Proceedings of the 40th International Conference on Machine Learning},
  series    = {Proceedings of Machine Learning Research},
  volume    = {202},
  pages     = {21426--21449},
  year      = {2023},
  publisher = {PMLR}
}

@article{PACOther,
  author  = {Sambharya, Rajiv and Stellato, Bartolomeo},
  title   = {Data-driven performance guarantees for classical and learned optimizers},
  journal = {Journal of Machine Learning Research},
  volume  = {26},
  number  = {171},
  pages   = {1--49},
  year    = {2025}
}

@article{l2OPACSucker,
  author  = {Sucker, Michael and Fadili, Jalal and Ochs, Peter},
  title   = {Learning-to-optimize with {PAC}-{B}ayesian guarantees: Theoretical considerations and practical implementation},
  journal = {Journal of Machine Learning Research},
  volume  = {26},
  number  = {211},
  pages   = {1--53},
  year    = {2025}
}

@inproceedings{STL10,
  author    = {Coates, Adam and Ng, Andrew and Lee, Honglak},
  title     = {An analysis of single-layer networks in unsupervised feature learning},
  booktitle = {Proceedings of the Fourteenth International Conference on Artificial Intelligence and Statistics},
  pages     = {215--223},
  year      = {2011},
  publisher = {JMLR Workshop and Conference Proceedings}
}

@article{kingma2014adam,
  author  = {Kingma, Diederik P. and Ba, Jimmy},
  title   = {Adam: A method for stochastic optimization},
  journal = {arXiv preprint arXiv:1412.6980},
  year    = {2014}
}

@misc{ODL,
  author    = {Adler, Jonas and Kohr, Holger and {\"O}ktem, Ozan},
  title     = {Operator discretization library ({ODL})},
  year      = {2017},
  note      = {Zenodo},
  publisher = {Zenodo},
  doi       = {10.5281/zenodo.249479}
}

@article{grimmer2024longsteps,
  author  = {Grimmer, Benjamin},
  title   = {Provably faster gradient descent via long steps},
  journal = {SIAM Journal on Optimization},
  volume  = {34},
  number  = {3},
  pages   = {2588--2608},
  year    = {2024}
}

@article{silversteps,
  author  = {Altschuler, Jason and Parrilo, Pablo},
  title   = {Acceleration by stepsize hedging: Multi-step descent and the silver stepsize schedule},
  journal = {Journal of the ACM},
  volume  = {72},
  number  = {2},
  pages   = {Article 12, 1--38},
  year    = {2025},
  doi     = {10.1145/3708502}
}

@article{learningalgohyperparams,
  author  = {Sambharya, Rajiv and Stellato, Bartolomeo},
  title   = {Learning algorithm hyperparameters for fast parametric convex optimization},
  journal = {arXiv preprint arXiv:2411.15717},
  year    = {2024}
}

@article{WCRR,
  author  = {Goujon, Alexis and Neumayer, Sebastian and Unser, Michael},
  title   = {Learning weakly convex regularizers for convergent image-reconstruction algorithms},
  journal = {SIAM Journal on Imaging Sciences},
  volume  = {17},
  number  = {1},
  pages   = {91--115},
  year    = {2024}
}

@article{ASTRA,
  author  = {van Aarle, Wim and Palenstijn, Willem Jan and Cant, Jeroen and Janssens, Eline and Bleichrodt, Folkert and Dabravolski, Andrei and De Beenhouwer, Jan and Batenburg, Kees Joost and Sijbers, Jan},
  title   = {Fast and flexible {X}-ray tomography using the {ASTRA} toolbox},
  journal = {Optics Express},
  volume  = {24},
  number  = {22},
  pages   = {25129--25147},
  year    = {2016},
  doi     = {10.1364/OE.24.025129}
}

@article{borgerding2017amp,
  title={AMP-inspired deep networks for sparse linear inverse problems},
  author={Borgerding, Mark and Schniter, Philip and Rangan, Sundeep},
  journal={IEEE Transactions on Signal Processing},
  volume={65},
  number={16},
  pages={4293--4308},
  year={2017},
  publisher={IEEE}
}

@article{metzler2017learned,
  title={Learned D-AMP: Principled neural network based compressive image recovery},
  author={Metzler, Chris and Mousavi, Ali and Baraniuk, Richard},
  journal={Advances in Neural Information Processing Systems},
  volume={30},
  year={2017}
}

@article{chen2018theoretical,
  title={Theoretical linear convergence of unfolded ISTA and its practical weights and thresholds},
  author={Chen, Xiaohan and Liu, Jialin and Wang, Zhangyang and Yin, Wotao},
  journal={Advances in Neural Information Processing Systems},
  volume={31},
  year={2018}
}

@article{wilson2021lyapunov,
  author  = {Wilson, Ashia C. and Recht, Ben and Jordan, Michael I.},
  title   = {A Lyapunov Analysis of Accelerated Methods in Optimization},
  journal = {Journal of Machine Learning Research},
  volume  = {22},
  pages   = {1--34},
  year    = {2021}
}

@article{drori2014performance,
  title={Performance of first-order methods for smooth convex minimization: a novel approach},
  author={Drori, Yoel and Teboulle, Marc},
  journal={Mathematical Programming},
  volume={145},
  number={1},
  pages={451--482},
  year={2014},
  publisher={Springer}
}

@article{taylor2017smooth,
  title={Smooth strongly convex interpolation and exact worst-case performance of first-order methods},
  author={Taylor, Adrien B and Hendrickx, Julien M and Glineur, Fran{\c{c}}ois},
  journal={Mathematical Programming},
  volume={161},
  number={1},
  pages={307--345},
  year={2017},
  publisher={Springer}
}

@article{taylor2017exact,
  title={Exact worst-case performance of first-order methods for composite convex optimization},
  author={Taylor, Adrien B and Hendrickx, Julien M and Glineur, Fran{\c{c}}ois},
  journal={SIAM Journal on Optimization},
  volume={27},
  number={3},
  pages={1283--1313},
  year={2017},
  publisher={SIAM}
}

@article{condat2026nesterov,
  title={A Nesterov-Accelerated Primal-Dual Splitting Algorithm for Convex Nonsmooth Optimization},
  author={Condat, Laurent and Sadiev, Abdurakhmon and Richtarik, Peter},
  journal={arXiv preprint arXiv:2604.09245},
  year={2026}
}

@article{ochs2014ipiano,
  title={iPiano: Inertial proximal algorithm for nonconvex optimization},
  author={Ochs, Peter and Chen, Yunjin and Brox, Thomas and Pock, Thomas},
  journal={SIAM Journal on Imaging Sciences},
  volume={7},
  number={2},
  pages={1388--1419},
  year={2014},
  publisher={SIAM}
}

@article{attouch2013convergence,
  title={Convergence of descent methods for semi-algebraic and tame problems: proximal algorithms, forward--backward splitting, and regularized Gauss--Seidel methods},
  author={Attouch, Hedy and Bolte, J{\'e}r{\^o}me and Svaiter, Benar Fux},
  journal={Mathematical Programming},
  volume={137},
  number={1},
  pages={91--129},
  year={2013},
  publisher={Springer}
}

@article{donoho2009message,
  title={Message-passing algorithms for compressed sensing},
  author={Donoho, David L and Maleki, Arian and Montanari, Andrea},
  journal={Proceedings of the National Academy of Sciences},
  volume={106},
  number={45},
  pages={18914--18919},
  year={2009},
  publisher={National Academy of Sciences}
}
\bibliographystyle{siamplain}

\appendix

% If instead we chose to learn $T$ sets of parameters $\theta_t$ for $t \in \{ 0, 1, 2, \dots, T-1\}$ simultaneously, such that
% \begin{equation}
%     (\theta_0, \dots, \theta_{t}) = \min_{\tilde{\theta}_0, \dots, \tilde{\theta}_{T-1}} f(x^t(\tilde{\theta}_0, \dots, \tilde{\theta}_{T-1})),
% \end{equation}
% we have obtained a nonconvex optimization problem for $T>1$, where
% \begin{equation}
%     x^{t+1}(\tilde{\theta}_0, \dots, \tilde{\theta}_{t}) = x^t(\tilde{\theta}_0, \dots, \tilde{\theta}_{t-1}) - G_{\theta_t}(\nabla f (x^t(\tilde{\theta}_0, \dots, \tilde{\theta}_{t-1}))).
% \end{equation}

% As $G_{\theta}\nabla f_k (x_k^t)$ is linear in $\theta$ and each function $f_k$ is convex, this optimization problem is convex. In this optimization problem, given the current iterates $x_k^t$ for $k \in \{1, \dots, N\}$, we seek to choose the optimal greedy parameters $\theta_t$ at iteration $t$ such that we minimize the mean over every $f_k(x_k^{t+1})$. In L2O, learning parameters is often a non-convex optimization problem. Therefore, the performance of learned optimizers is highly dependent on the optimization algorithm used and its hyperparameters. However, as our problem is convex, one can use any convex optimization algorithm with convergence guarantees.

\section{Further Notation}

% A function $f: \mathcal{X} \rightarrow \mathbb{R}$ is strongly convex with parameter $\mu > 0$ if $f - \mu\| \cdot \|^2/2$ is convex. We say $f \in \mathcal{F}_{L, \mu}$ if $f \in \mathcal{F}_{L}$ and $f$ is $\mu$-strongly convex. 
 For elements $x, y, z \in \mathcal{X}$, define the linear operator $x \otimes y$, the outer product, by $(x \otimes y)z := \langle y, z \rangle x$, with the property that $[x \otimes y]_{q,i} = \langle y, e_i \rangle \langle x, e_q \rangle = x_q y_i$.
% For two linear operators $A, B \in \mathcal{L}(\mathcal{X})$, define the Kronecker product $A \otimes B$ by
% \begin{equation} \label{eq--AotimesB-identity}
%     [A \otimes B]_{ij,kl} = A_{i,k}B_{j,l}.
% \end{equation}
% \noindent\textbf{\cref{prop--immediate-full-conv}.} 
%     For $k \in \{1,\dots, N\}$, assume that $f_k: \mathcal{X} \to \mathbb{R}$ is convex, continuously differentiable and bounded below, with any initial point $x_k^0 \in \mathcal{X}$. Assume that the set of gradients $\{ \nabla f_1(x_1^0), \dots,  \nabla f_N(x_N^0)\}$ is linearly independent. Then if $N \leq n$, there exists an operator $P \in \mathcal{L}(\mathcal{X})$ such that $x_k^0 - P  \nabla f_k(x_k^0) \in \arg\min_x f_k(x)$, for all $k \in \{1, \dots, N\}$.  

% \Alternative ? \begin{proof}
%     Find $P$ such that $PA=B$, denote by $C_i$ the $i$th row of a matrix C and $C_{\cdot, j}$ the $j$th column. Therefore $\langle P_i, A_{\cdot, j} \rangle = B_{i,j}$ for $i \in \{1, \dots, n\}, j \in \{1, \dots, N\}$. Meaning that $A^TP_i = B_i$ for all $i \in \{1, \dots, n\}$, and such a vector $P_i$ exists if $B_i$ lies in the range of $A^T$, which is equivalent to saying that $B_i$ lies in the column space of $A$. So $B_i$ can be written as as a linear combination of the rows of $A$, equivalently the row rank of $B$ is at most the row rank of $A$ which is equivalent to $\operatorname{rank}(B) \leq \operatorname{rank}(A)$. In particular. this can occur when $N \geq 2$.
% \end{proof}}

\section{Proofs for \Cref{section--convergence}} \label{app:convergence}
\\

\begin{proof}[Proof of \cref{prop:F-properties}]
Firstly, define the norm on $\mathcal{X}^N$ for $x=(x_1, \dots, x_N) \in \mathcal{X}^N$  by $\| x \| = \sqrt{\sum_{k=1}^N \|x_k\|^2}$. For any $x, y \in \mathcal{X}^N$, we have $\nabla F(x) = 1/N \left(\nabla f_1(x_1), \dots, \nabla f_N(x_N) \right)$ and therefore 
\begin{align*}
    \| \nabla F(x) - \nabla F(y) \| &= \frac1N\sqrt{ \sum_{k=1}^N \left \| \nabla f_k (x_k) - \nabla f_k(y_k) \right \|^2}\\
&\leq \frac1N\sqrt{  \sum_{k=1}^N L_k^2 \left \| x_k - y_k \right \|^2} %\quad \text{($L_k$-smoothness of $f_k$.)}\\
\leq  \frac{\max\{L_1, \dots, L_N\}}N \| x-y \|.
\end{align*}
\end{proof}

%\begin{proof}
% \textbf{Proof of Remark \cref{remark-L-range}.}
\begin{proof}[Proof of \cref{remark-L-range}]
\textbf{Case 1} - $L \ge 1/\tau$.
Using condition \cref{eq::P} in this case gives
\begin{equation*}
        L < \frac{ (2\tau + 2\varepsilon_1 + \varepsilon_2) - \varepsilon_2 (\tau + \varepsilon_1)^2 - 2 \varepsilon_2^2(\tau + \varepsilon_1) - \varepsilon_2^3 }{ (\tau + \varepsilon_1)(\tau + \varepsilon_1 + \varepsilon_2) } =  \frac{1 - (\tau+\varepsilon_1)\varepsilon_2 - \varepsilon_2^2}{\tau+\varepsilon_1} + \frac{1}{\tau+\varepsilon_1+\varepsilon_2}.
\end{equation*}
Taking the partial derivative with respect to $\varepsilon_2$ of this upper bound gives:
\begin{equation*}
   -1 - \frac{2\varepsilon_2}{\tau+\varepsilon_1} - \frac{1}{(\tau+\varepsilon_1+\varepsilon_2)^2} < 0,
\end{equation*}
as $\tau > 0$, $\varepsilon_1 \ge 0$, $\varepsilon_2 \ge 0$. Therefore for any fixed $\varepsilon_1 \ge 0$, the maximum occurs at $\varepsilon_2 = 0$. With $\varepsilon_2 = 0$, the upper bound becomes $2/(\tau+\varepsilon_1)$, which is maximised at $\varepsilon_1 = 0$. Therefore, the maximum value of the upper bound is $2/\tau = 2L_{\text{train}}$, at $(\varepsilon_1, \varepsilon_2) = (0,0)$.

\textbf{Case 2} - $L < 1/\tau$, then define $h_2(\varepsilon_1, \varepsilon_2)$ as
\begin{equation*}
\frac{\left(\tau - \varepsilon_1 - \frac{\varepsilon_2}{2}\right)\left(1 - \varepsilon_2(\varepsilon_1 + \varepsilon_2 + \tau)\right) - \frac{\varepsilon_2}{2} (\tau + \varepsilon_1)(\tau + \varepsilon_1 + \varepsilon_2)}{\frac{\varepsilon_2}{2} (\varepsilon_1 + \varepsilon_2 - \tau) (\tau + \varepsilon_1)(\tau + \varepsilon_1 + \varepsilon_2) + \left(\frac{\tau^2}{2} - \varepsilon_1\tau - \frac{\varepsilon_2\tau}{2} + \frac{\varepsilon_1^2}{2} + \frac{\varepsilon_1\varepsilon_2}{2}\right) \left(1 - \varepsilon_2(\varepsilon_1 + \varepsilon_2 + \tau)\right)},
\end{equation*}
then in this case, the condition \cref{eq::P} is equivalent to $L < h_2 (\varepsilon_1, \varepsilon_2)$. Note that $h_2(\varepsilon_1, 0) = 2(\tau - \varepsilon_1)/(\tau - \varepsilon_1)^2 = 2/(\tau - \varepsilon_1) $ if $\varepsilon_1 \neq \tau$. Then
\begin{equation*}
    \lim_{\varepsilon_1 \to \tau^+} h_2(\varepsilon_1, 0) = \lim_{\varepsilon_1 \to \tau^+} \frac2 {\tau - \varepsilon_1} = - \infty.
\end{equation*}
\end{proof}

\begin{lemma} \label{lemma::convergence_grad_fk}
Suppose that $F(x^{t+1}) \leq F\left(x^t\right) - 1/2L_F \| \nabla F\left(x^t\right) \|^2$. Then \\ $\sum_{n=0}^t \| \nabla F (x^n) \|^2  \leq (F\left(x^0\right) - F^*)/(2 L_F )$ and $\nabla f_k (x_k^t) \to 0$ as $t \to \infty$ for all $k \in \{1, \dots, N\}$.
\begin{proof}
    Rearranging and taking a summation from $0$ to $t$, we get
    \begin{align*}
        \sum_{n=0}^t \| \nabla F (x^n) \|^2  \leq \sum_{n=0}^t\frac{F(x^n) - F(x^{n+1})}{2L_F} 
       =\frac{F\left(x^0\right) - F(x^{t+1})}{2 L_F} \leq\frac{F\left(x^0\right) - F^*}{2 L_F },
    \end{align*}
    therefore $\lim_{t \to \infty} \sum_{n=0}^t \| \nabla F (x^n) \|^2  \leq (F\left(x^0\right) - F^*)/(2 L_F )$, meaning that $ \| \nabla F\left(x^t\right) \|^2 \to 0$ as $t \to \infty$. Finally, $\|\nabla F\left(x^t\right) \|^2 = 1/N^2 \sum_{k=1}^N \| \nabla f_k(x_k^t)\|^2 \to 0$ as $t \to \infty$, which implies that $\nabla f_k (x_k^t) \to 0$ as $t \to \infty$ for all $k \in \{1, \dots, N\}$.
\end{proof}
    
\end{lemma}

\begin{corollary} [Example error bounds for $L=L_{\text{train}}$] \label{cor::ex-error-bounds}
Let \cref{assumption::liminfbgdcts} hold. Then, there exists a final training iteration $T$ such that for all $f \in \mathcal{F}_{L_{\text{train}}}$ and any starting point $x^0$, using \cref{alg--test-algo} we have $\nabla f\left(x^t\right) \rightarrow 0$ as $t \rightarrow \infty$, and the convergence rate in \cref{eq:nonconvex-conv-rate} holds. In particular, this rate follows if at iteration $T$, $\| G_{\theta_T} - \tau I \| \le \varepsilon_1 $  and $\| H_{\phi_T} \| \le \varepsilon_2 $, where $ \varepsilon_1 = \tau / 2$, and  $\varepsilon_2 = \tau/(3(\tau^2 + 1))$.
\end{corollary}

\begin{proof}
Due to $G_{\theta_t} \to \tau I$ and $H_{\phi_t} \to 0$ as $t \to \infty$ there exists a final training iteration $T$ such that $\| G_{\theta_T} - \tau I \| \le \varepsilon_1 $  and $\| H_{\phi_T} \| \le \varepsilon_2 $, where $\varepsilon_1 = \tau/2$, and $\varepsilon_2 = \tau/(3(\tau^2 + 1))$. Using \cref{eq::general-ineq-errors}, in the case $L = L_{\text{train}}$, for $t > T$, we have
\begin{align*}
    f\left(x^{t+1}\right) & \le  f\left(x^t\right) - \frac 1 {2 \tau}\left(\tau^2 - \varepsilon_1^2 - \varepsilon_1 \varepsilon_2\right) \| \nabla f\left(x^t\right) \|^2 
    + \frac {\varepsilon_2} {2\tau } (\varepsilon_1 + \varepsilon_2)\ \|  x^t - x^{t-1} \|^2.
\end{align*}
Then the constants $C_1, C_2$ are now given by
% \begin{align*}
%     C_1(\varepsilon_1, \varepsilon_2, \gamma) &=  \tau - \varepsilon_1 - \frac {\varepsilon_2} 2 - \frac L 2  (\tau + \varepsilon_1)\left(\varepsilon_1 +  \varepsilon_2 + \tau \right) - \gamma  (\tau + \varepsilon_1)\left(\varepsilon_1 +  \varepsilon_2 + \tau \right)\\
%     C_2(\varepsilon_1, \varepsilon_2, \gamma) &=  \gamma \left( 1 -  \varepsilon_2\left( \varepsilon_1 + \varepsilon_2 + \tau \right) \right) -\frac {\varepsilon_2} 2 ( 1 + L( \varepsilon_1 + \varepsilon_2 + \tau)),
% \end{align*}
\begin{align*}
    C_1(\varepsilon_1, \varepsilon_2, \gamma) &=  \frac 1 {2 \tau}\left(\tau^2 - \varepsilon_1^2 - \varepsilon_1 \varepsilon_2 - 2\gamma \tau (\tau + \varepsilon_1)\left(\varepsilon_1 +  \varepsilon_2 + \tau\right)\right), \\
    C_2(\varepsilon_1, \varepsilon_2, \gamma) &= \frac 1 {2 \tau}\left(-\varepsilon_2 (\varepsilon_1 + \varepsilon_2) + 2 \gamma \tau \left( 1 -  \varepsilon_2\left( \varepsilon_1 + \varepsilon_2 + \tau \right) \right) \right).
\end{align*}
For $C_1(\varepsilon_1, \varepsilon_2, \gamma), C_2(\varepsilon_1, \varepsilon_2, \gamma) > 0$, with $P$ defined as in \cref{eq::P}, we now have the condition
\begin{equation*}
   P(\varepsilon_1, \varepsilon_2) = 2 \tau \varepsilon_2 (\tau + \varepsilon_1)(\tau + \varepsilon_1 + \varepsilon_2) ( \varepsilon_1 +\varepsilon_2) - 2 \tau\left(\tau^2 - \varepsilon_1^2 - \varepsilon_1 \varepsilon_2\right)\left( 1 -  \varepsilon_2\left( \varepsilon_1 + \varepsilon_2 + \tau \right) \right) < 0.
\end{equation*}
Then, $P\left(\tau/2, \varepsilon_2\right) = \tau^2/2 \left(4\varepsilon_2^3 + 12\tau\varepsilon_2^2 + (9\tau^2 + 2)\varepsilon_2 - 3\tau \right)$,
% and 
% \begin{equation}
%      \frac{2}{\tau^2} P\left(\frac{\tau}{2}, \varepsilon_2\right) = 4\varepsilon_2^3 + 12\tau\varepsilon_2^2 + (9\tau^2 + 2)\varepsilon_2 - 3\tau,
% \end{equation}
and 
\begin{equation*}
     \frac{2 }{\varepsilon_2^3} P\left(\frac{\tau}{2}, \frac \tau {3(1 + \tau^2)}\right) = - 27 \tau ^4 - 90 \tau^2 - 63 < 0.
\end{equation*}
Therefore, as in \cref{thm--reg-conv}, $C_1, C_2 > 0$, as the convergence results are obtained.
\end{proof}

\begin{lemma} \label{lemma-nesterov-lemma}
    Let $\Delta_t \in [0, \infty)$ be a sequence such that $\Delta_{t+1} \leq \Delta_t- c\Delta_t^2$ for some constant $c > 0$ and $\Delta_t \ge \Delta_{t+1}$. Then $\Delta_t \leq 1/(ct)$ for all $t \geq 0$.
\end{lemma}
\begin{proof}
    If for some $t$, $\Delta_t = 0$ then $\Delta_t = 0$ for all $t > T$ and we are done. Otherwise assume $\Delta_t, \Delta_{t+1} \neq 0$. As in \citep{nesterov2018lectures}, by dividing by $\Delta_t\Delta_{t+1}$ both sides,
    \begin{align*}
        &\frac1 {\Delta_t} \leq \frac1 {\Delta_{t+1}} - c \frac{\Delta_t}{\Delta_{t+1}} \leq \frac1 {\Delta_{t+1}} - c\implies c + \frac1 {\Delta_t} \leq \frac1 {\Delta_{t+1}}.
    \end{align*}
    Taking a summation gives
    \begin{align*}
        &\sum_{k=0}^{t-1} c \leq \sum_{k=0}^{t-1} \left(\frac1 {\Delta_{k+1}} - \frac1 {\Delta_k} \right)  \implies ct \leq \frac1 {\Delta_t} - \frac1 {\Delta_0}.
    \end{align*}
    Therefore
    \begin{align*}
        \Delta_t \leq \frac{1}{\frac1 {\Delta_0} + ct} = \frac{ \Delta_0}{1 + c \Delta_0 t} \leq \frac{\Delta_0}{c \Delta_0 t} = \frac{1}{ct}.
    \end{align*}
\end{proof}

\section{Proofs for \Cref{section--linear-params}}

\subsection{Approximating optimal linear parameters} \label{appendix:approx-linear-params}

Using the general result in \cref{prop:approx-general}, we can calculate $\nabla_{\theta} g_{t, \lambda_t, \mu_t}$ and $\nabla_{\phi} g_{t, \lambda_t, \mu_t}$ and their associated Lipschitz constants for specific parametrizations of $G$ and $H$.
\begin{corollary}
    \textbf{Pointwise parametrization.} For the pointwise parametrization, $\theta \in \mathcal{X}$, the adjoints $(B_k^t)^* = B_k^t$, $(C_k^t)^* = C_k^t$, and $\|B_k^t \| \le \|\nabla f_k(x_k^t)\|_\infty, \| C_k^t \| \le \|x_k^t - x_k^{t-1}\|_\infty$.
            
    \textbf{Full operator parametrization.} In this case we have $\theta \in \mathcal{L}(\mathcal{X})$. Then the adjoint is given by $(B_k^t)^* (w) = w \otimes \nabla f_k (x_k^t)$, $(C_k^t)^* (w) = w \otimes(x_k^t - x_k^{t-1})$, and $\| B_k^t \| \le  \left\| \nabla f_k(x_k^t)\right\|$,  $\| B_k^t \| \le  \left\| x_k^t - x_k^{t-1}\right\|$. 

    \textbf{Scalar step size.} We now take $\theta \in \mathbb{R}$, then the adjoints are given by $(B_k^t)^*w = \langle w, \nabla f_k (x_k^t) \rangle$ and $(C_k^t)^*w = \langle w, (x_k^t - x_k^{t-1}) \rangle$, and $\|B_k^t \| =  \| \nabla f_k(x_k^t)  \|$, and $\|C_k^t \| =  \| x_k^t - x_k^{t-1}  \|$.

    \textbf{Convolution.} In this case we have $\theta \in \mathbb{R}^{n_1 \times n_2}$. Define the discrete Fourier transform by $\mathcal{F}$, then $B_k^t = \ft^{-1} \operatorname{diag}(\ft\nabla f_k(x_k^t))\ft$, $C_k^t = \ft^{-1} \operatorname{diag}(\ft(x_k^t - x_k^{t-1}))\ft$, the adjoints are given by $(B_k^t)^* = B_k^t$ and $(C_k^t)^* = C_k^t$, $\|B_k^t \| \le \|\mathcal{F} \nabla f_k(x_k^t)\|_\infty, \| C_k^t \| \le \| \mathcal{F}(x_k^t - x_k^{t-1})\|_\infty$.
    
    % In this case we have $\theta \in \mathbb{R}^{n_1 \times n_2}$. The adjoints are given by $(B_k^t)^*w = w\ast \overline{\nabla f_k (x_k^t)}$ and $(C_k^t)^*w = w\ast \overline{(x_k^t - x_k^{t-1})}$. 
\end{corollary}

\begin{proof}
    \textbf{Pointwise parametrization.} In this case, we have $\theta \in \mathcal{X}$ and $B_k^t(x) = \nabla f_k(x_k^t) \odot x$. For the adjoints, $\langle B_k^t(x), w \rangle = \langle x \odot \nabla f_k(x_k^t), w \rangle = \langle x, w \odot \nabla f_k(x_k^t) \rangle = \langle x,(B_k^t)^* w \rangle$, and so $(B_k^t)^* = B_k^t$ and similarly, $(C_k^t)^* = C_k^t$. Furthermore,
\begin{align*}
    \|B_k^t\| &= \max_{x  \neq 0} \frac{\|x \odot \nabla f_k(x_k^t))\|}{ \| x \|} = \max_{x  \neq 0} \sqrt{\frac{\sum_{i=1}^n x_i^2 [\nabla f_k (x_k^t)]_i^2}{ \sum_{i=1}^n x_i^2}} \\
    &\leq  \max_q |[\nabla f_k (x_k^t)]_q| \max_{x  \neq 0} \sqrt{\frac{ \sum_{i=1}^n x_i^2 }{ \sum_{i=1}^n x_i^2}}   = \| \nabla f_k(x_k^t) \|_\infty.
\end{align*}
\textbf{Full operator parametrization.} $\theta_t \in \mathcal{L}(\mathcal{X}) $ and for the adjoint of $B_k^t$,  
\begin{align*}
    \langle B_k^t(P), w \rangle &= \sum_{i=1}^n [P \nabla f_k (x_k^t)]_i w_i 
    = \sum_{i=1}^n \sum_{j=1}^n [P]_{i,j} [\nabla f_k (x_k^t)]_j w_i \\
    &= \langle P, (B_k^t)^* w \rangle 
    = \sum_{i=1}^n \sum_{j=1}^n [P]_{i,j} [(B_k^t)^* w]_{i,j},
\end{align*}
and therefore $[(B_k^t)^* w]_{i,j} = w_i [\nabla f_k (x_k^t)]_j$, which means $(B_k^t)^* (w) = w \otimes \nabla f_k (x_k^t)$. Similarly, $(C_k^t)^* (w) = w \otimes (x_k^t - x_k^{t-1})$. For the Lipschitz constants, note that $\| B_k^t (P) \| = \| P \nabla f_k(x_k^t) \| \leq \| P \| \| \nabla f_k(x_k^t) \|$, and therefore 
\begin{align*}
    \| B_k^t \| &= \max_{P  \neq 0} \frac{\| B_k^t (P) \|}{ \| P \|} \leq \| \nabla f_k(x_k^t) \|.
\end{align*}
\textbf{Scalar step size.}  $B_k^t(\alpha) = \alpha \nabla f_k (x_k^t)$. For $\alpha \in \mathbb{R}$, $\langle B_k^t(\alpha), w \rangle = \langle \alpha \nabla f_k (x_k^t), w \rangle = \alpha \langle \nabla f_k (x_k^t), w \rangle$, and therefore $(B_k^t)^*(w) = \langle \nabla f_k (x_k^t), w \rangle$, and similarly $(C_k^t)^*(w) = \langle (x_k^t - x_k^{t-1}), w \rangle$. Furthermore, $\| B_k^t(\alpha)\| = \| \alpha \nabla f_k (x_k^t) \| = | \alpha | \| \nabla f_k (x_k^t) \|$, and so 
    \begin{align*}
        \| B_k^t \| &= \max_{\alpha \neq 0} \frac{\| B_k^t(\alpha)\|}{|\alpha|} = \| \nabla f_k (x_k^t) \|.
    \end{align*}
    \textbf{Convolution.} For the operators $B_k^t$ and $C_k^t$, for $v \in \mathcal{X}$, as $\mathcal{F}$ is bijective, and $\ft(\kappa \ast v) = \ft\kappa \odot \ft v$, we have
\begin{align*}
    \kappa \ast v = \ft^{-1} \ft(\kappa \ast v) =  \ft^{-1} (\ft\kappa \odot \ft v) = \ft^{-1} \operatorname{diag}(\ft v)\ft\kappa
\end{align*}
Then as $\ft^{*} = c\ft^{-1}$ for some constant $c$ \cite{briggs1995dft} (this constant depends on the implementation of the discrete Fourier transform. Often $c=n_1n_2$ or $c=1$), $(\ft^{-1})^* = \tfrac1c \ft$. $(B_k^t)^* = \ft^{*} \operatorname{diag}(\ft\nabla f_k(x_k^t))(\ft^{-1})^* = \ft^{-1} \operatorname{diag}(\ft\nabla f_k(x_k^t)) \ft = B_k^t$. Lastly,
\begin{align*}
\|B_k^t\| &= \left\| \mathcal{F}^{-1} \operatorname{diag}(\mathcal{F}\nabla f_k(x_k^t))\mathcal{F} \right\| = \left\| \operatorname{diag}(\mathcal{F}\nabla f_k(x_k^t)) \right\| \le \left\| \mathcal{F}\nabla f_k(x_k^t) \right\|_\infty,
\end{align*}
where the second equality follows as the operator norm is unitarily invariant.
\end{proof}

One also has the option to approximate the operator norms of $B_k^t$ and $C_k^t$ using the power method \cite{golub2013matrix}.

\section{Small-scale CT problem}

\textbf{Problem details.} We use $90$ projection angles and extract $40 \times 40$ pixel crops from the center of each ground-truth image in the dataset. The noise $\varepsilon$ is modeled with a standard deviation of $ 10^{-2}$, and we set $\alpha = 10^{-4}$, resulting in $L = 1.08$.

\textbf{Learning Parameters.}  We use a training set of $25$ functions for parametrizations PS, PP, and PC. 
For PF, the model is trained using $1000$ functions. Testing is performed on a separate set of $100$ functions. 

\textbf{Training details.} Greedy training was performed up to iteration $T=200$ with $\lambda_t = 0$ for all iterations $t$ for the PS, PP, and PC parametrizations. The total time for training PS was about $10$ minutes, for PP was about $67$ minutes, and PC took approximately $10$ hours. For the PF parametrization, training was performed up to iteration $T=11$ with $\lambda_t = 0$ for all $t$. Furthermore, the PF parametrization was trained with regularization such that $\lambda_t = 10^{-10}$ for $t < T=101$ iterations. At iteration $101$, the learned operator $G_{\theta_T}$ satisfied $\| G_{\theta_T} - \tau I \| < \tau$, guaranteeing convergence on iterations $t \geq T$. For each iteration $t$, solving the optimization problem \cref{eq--opt-problem} with the PF parametrization took one hour.

\begin{figure}[h!]
    \centering
    \begin{subfigure}[t]{\textwidth}
        \centering
        \includegraphics[width=.5\textwidth]{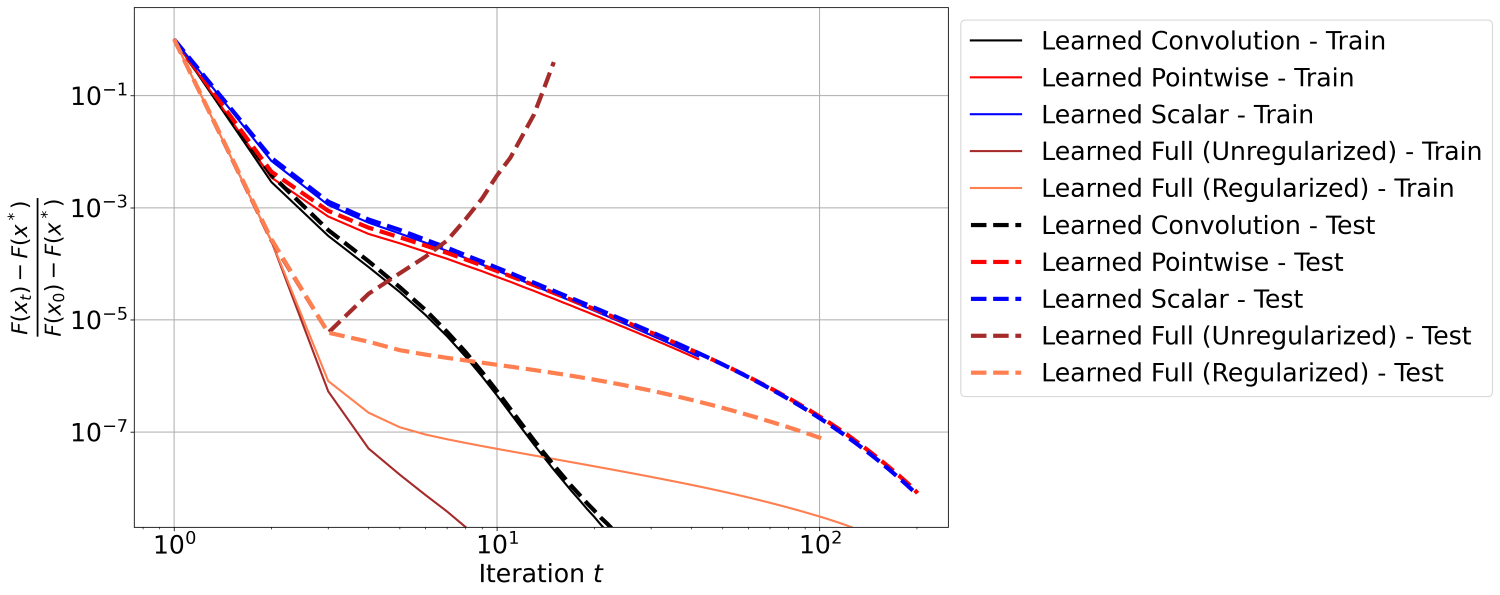} % Replace with your image file
        
        \caption{}
        \label{6a}
    \end{subfigure}
    \newline
    \begin{subfigure}[t]{0.275\textwidth}
        \centering
        \includegraphics[width=\textwidth]{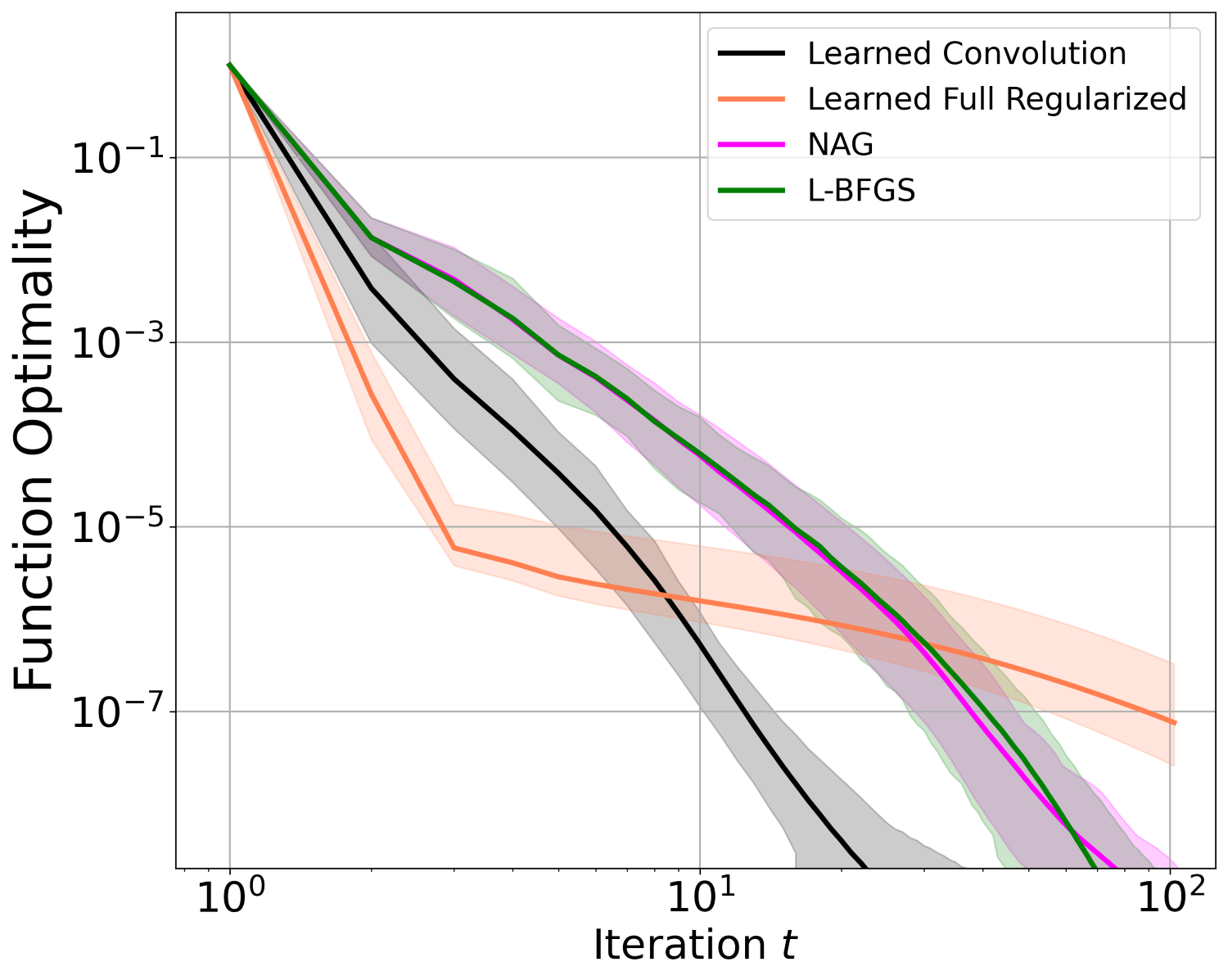} % Replace with your image file
        
        \caption{}
        \label{6b}
    \end{subfigure}
    \begin{subfigure}[t]{0.275\textwidth}
        \centering
        \includegraphics[width=\textwidth]{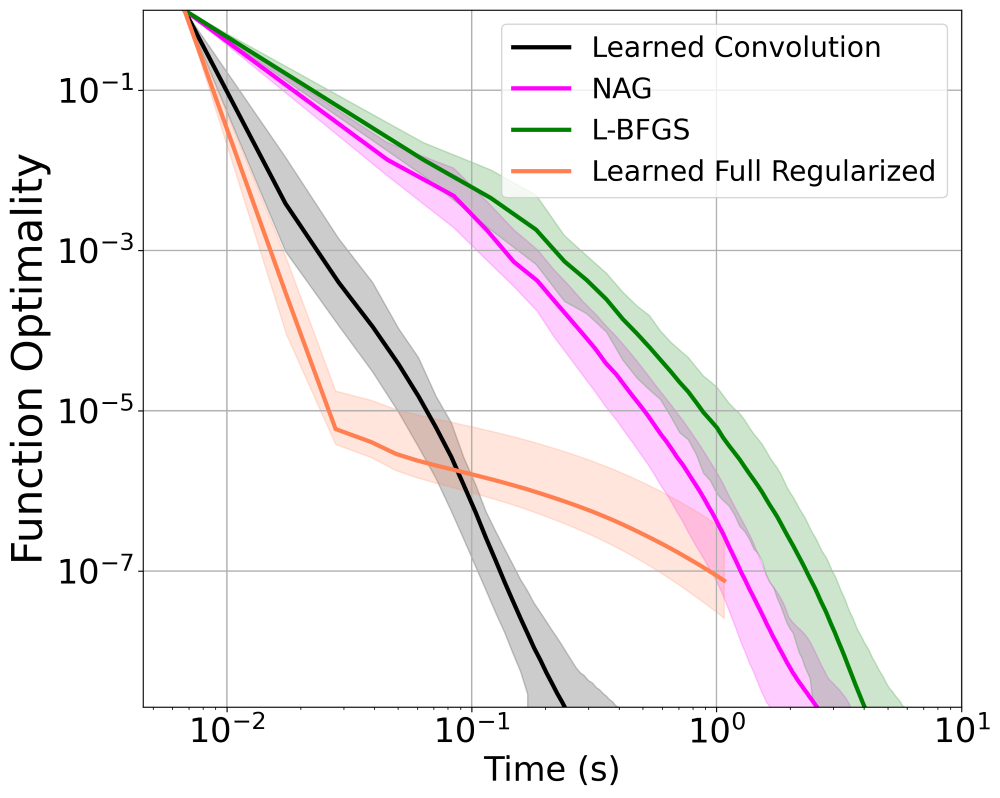} % Replace with your image file
        
        \caption{}
        \label{6c}
    \end{subfigure}
    
    \caption{Performance of learned algorithms for the small-scale CT problem. (a) Train versus test set performance of the learned parameterizations. (b) Test performance versus benchmark optimization algorithms. (c) Wall-clock test performance versus benchmarks.}
    \label{fig--CT-gen-test}
\end{figure}

\textbf{Learned algorithm performance.} \cref{6a} shows that the learned parametrizations PS, PP, and PC generalize well to test data for Problem 2. Again the learned PC parametrization outperforms NAG and L-BFGS on the CT test data as shown by \cref{6b}, reaching a tolerance of $10^{-10}$ in approximately $30$ iterations, compared with about $ 80$ for L-BFGS and NAG. \cref{6c} shows that PC also outperforms in terms of wall-clock time.

\textbf{Full operators.} The full parametrization PF shows signs of overfitting, as it does not generalize well to test data. It performs well in the first two iterations, but then diverges. The PF parametrization with regularization mitigates this issue, as the generalization performance is seen to improve. \cref{6b} shows it initially converges quickly but its speed decreases later due to regularization. This is because, with increasing iterations, the learned update gets closer to gradient descent.
% \begin{figure}[h!]
%     \centering
%     \begin{subfigure}[t]{0.27\textwidth}
%         \centering
%         \includegraphics[width=\textwidth]{CT_images/generalisation.png} % Replace with your image file
%     \end{subfigure}\hspace{20px}
%     \begin{subfigure}[t]{0.27\textwidth}
%         \centering
%         \includegraphics[width=\textwidth]{CT_images/ct_test_all.png} % Replace with your image file
%     \end{subfigure}
%     \caption{ Left: Train versus test performance of learned methods for Problem 2. Right: Test performance versus benchmark optimization algorithms.}
%     \label{fig--CT-gen-test}
% \end{figure}

\end{document}